# Foundations of Geometry

Ruth Moufang
with the assistance of H. Rixecker

Translated and edited by John Stillwell

# Contents

# Introduction

1. The following investigations deal with the unprovable assertions on which geometry is founded. Every mathematical theory is based on certain unprovable assertions; the separation between provable and unprovable is made by formulating axioms or postulates. Euclid's *Elements* (325 BCE) are the first example of an axiomatic foundation for elementary geometry, particularly admirable in the case of the theory of proportions. The main goals of the developments that have flowed from Euclid are:

   (a) Investigation of the parallel axiom. In the 19th century it was first recognized that the parallel axiom is unprovable from the remaining axioms of Euclid, and that it is not reducible to simpler assumptions. However, there are other assertions equivalent to the parallel axiom. We shall not go further into these investigations, which concern the so-called noneuclidean geometry.

   (b) Investigations of continuity, particularly the axiom of Archimedes. The existence of incommensurable segments shows that the axiom of Archimedes does not suffice to introduce the number concept in geometry via coordinates. From the standpoint of geometric intuition, it is unsatisfactory to assume from the outset that geometric space is a number manifold. Continuity plays a crucial role in the foundation of the theory of proportions and the theory of planar and spatial content. We shall not go further into these investigations either. Here we wish to introduce coordinates into synthetic geometry independently of continuity assumptions. The Archimedean axiom then gives important insights into the structural properties of synthetic geometry via the detour through coordinate geometry.

   (c) In projective geometry, investigations of incidence properties are not supposed to involve the concept of magnitude. This raises the question of "pure foundations of projective geometry", that is, foundations without limit processes and congruence assumptions. Here it turns out that planar and spatial geometry are essentially different.

   The investigations that follow are concerned mainly with point (c). We leave aside all problems normally treated in basic courses in noneuclidean geometry or "elementary mathematics from an advanced standpoint". Our goal is to study ways of establishing projective geometry on the basis of the Hilbert axioms. The algebraic construction of geometry as coordinate geometry over a generalized number system plays an essential role.

2. We conclude this introduction with a few remarks on the axiomatic method.

   The axioms involve various properties of geometric figures: incidence (for example, two points determine exactly one line), order (for example, when three points lie on a line, exactly one of them is between the other two), congruence, continuity, and parallelism. Hilbert's axioms are correspondingly divided into groups I to V. Hilbert's division is not the only one possible. For example, Pasch based his investigations on another division, into so-called "core theorems". The essence of mathematics lies in its freedom.[1] But for each system one requires the axioms to be

   $\alpha$) independent,

   $\beta$) consistent,

   $\gamma$) complete.

   Requirement $\alpha$)—eliminating superfluous assumptions—is not only for the sake of elegance. Deciding what needs proof and what does not is the basic problem. To prove the independence of a given axiom system one has to show that an assertion $A$ does not follow from assertions $A_1$ to $A_n$. For a geometric axiom system, independence is proved, for example, by giving a model geometry in which assertions $A_1$ to $A_n$ hold but assertion $A$ is false. This idea is often used in what follows. The model need not be minimal, that is, it may have properties $B_1$ to $B_m$ as well as properties $A_1$ to $A_n$. If $A$ does not follow from $A_1$ to $A_n$ and $B_1$ to $B_m$ then it certainly does not follow from $A_1$ to $A_n$ alone.

   Consistency is a very deep problem, which we shall not expand on here. The consistency of geometry may be reduced to the consistency of arithmetic. The requirement of completeness lies at the foundation of all geometric disciplines. For example, a part of the axiom system may suffice to found a subdiscipline of geometry, say projective geometry, but not the whole of elementary geometry.

3. The investigations below are based on Hilbert's axiom system. The method he discovered and used, namely to connect general number systems with geometric theorems and thereby reduce geometric problems to algebraic problems, is systematically extended.

   The circle of problems in projective geometry is not confined to showing the hypotheses and paths leading to the fundamental theorem of projective geometry. It also includes structural questions, for example, the previously unsolved problem of classifying plane configuration theorems: given configuration theorems $S_1$ and $S_2$, one wants a procedure to decide whether $S_2$ follows from $S_1$ or not.

   The claim that $S_2$ does not follow from $S_1$ may be proved, for example, by giving a model geometry in which $S_1$ holds but $S_2$ is false. It certainly appears more natural to systematically generate all consequences of $S_1$ and present them so directly that the absence of $S_2$ is clear, but this path is connected with severe difficulties of principle.

---

[1]Translator's note. Here Moufang is quoting a remark of Cantor (almost word for word), from *Math. Ann.* 21, (1883), p. 564.

4. As prerequisites we assume familiarity with elementary geometry, analytic and projective geometry, and some basic concepts from the algebra of quaternions. Beyond that, foundational investigations require a certain mathematical maturity, an appreciation of structural questions rather than the simple desire to acquire knowledge. Naturally, proofs become longer when they are based on more restricted hypotheses.

# Chapter 1

# The axioms of geometry

## 1.1 Incidence axioms

1. Geometric figures consist of three types of things: points, lines and planes. We denote points by large Roman letters, lines by small Roman letters, and planes by large German letters. It is not necessary to attach an intuitive meaning to these concepts. Only the relations between them are important, for example, the concepts of "lying upon", "between", "congruent", "parallel", "continuous". These relations satisfy the following axioms.

   I.1 Through two distinct points $A$ and $B$ there is always a line $g$, the so-called connecting line of $A$ and $B$.

   I.2 Two distinct points have no more than one connecting line.

   I.3 On any line $g$ there are always at least two distinct points $A$ and $B$. There are at least three points not in the same line.

   I.4 Through any three distinct points not in the same line there is always a plane $\mathfrak{E}$, the so-called connecting plane of the three points. In each plane there is always at least one point.

   I.5 Three distinct points not in a line have no more than one connecting plane.

   I.6 Axiom of planes: if $\mathfrak{E}$ contains the points $A$ and $B$ ($A \neq B$) then any point on the connecting line of $A$ and $B$ also lies in $\mathfrak{E}$.

   I.7 If two planes have a point in common then they have at least one more point in common.

   I.8 There are at least four points not in the same plane.

   Axioms 1 to 3 are called the plane incidence axioms, and 4 to 8 are the space incidence axioms.

   In what follows, the expressions $A$ lies in $\mathfrak{E}$, $A$ is a point of $\mathfrak{E}$, $A$ belongs to $\mathfrak{E}$, $A$ is incident with $\mathfrak{E}$, will be used synonymously. Likewise, the expressions $A$ lies on $g$, $g$ goes through $A$, $A$ belongs to $g$, $A$ is incident with $g$, mean one and the same thing. When $A$ lies on both $g$ and $h$, $A$ is called the *intersection point* of g and h. The axioms are formulated so that each axiom contains only *one* requirement, and these requirements say as little as possible. For example, Axiom I.1 does not require that there be exactly one line through two points. Uniqueness is demanded only after existence has been postulated. Also, it is not required, for example, that each line contain many points. It suffices to demand the existence of two points, because the existence of further points can then be proved.

2. **Consequences of the incidence axioms.**

   (a) Two distinct lines of a plane have either no common point or at least one. In the latter case they have exactly one common point, by Axiom I.2.

   (b) A line $g$ and a point $A$ outside it determine exactly one plane $\mathfrak{E}$. Indeed, by Axiom I.3 there are at least two different points $B$ and $C$ on $g$. By I.4 and I.5 there is exactly one plane through $A$, $B$, and $C$, and by I.6 it contains the whole line $g$.

   (c) A plane $\mathfrak{E}$ and a line $g$ that does not lie in $\mathfrak{E}$ have either no point in common or—when a common point exists—only one, otherwise Axiom I.6 would be contradicted.

   (d) Two distinct lines $g$ and $h$ with one common point $A$ determine exactly one plane. Indeed, by I.3 there is a point $B$ on $g$ different from $A$, and a point $C$ on $h$ different from $A$. By I.4 there is at least one plane $\mathfrak{E}$ connecting the three points $A, B, C$ not in a line, and by I.5 no more than one. By I.6, all points of $g$ and all points of $h$ lie in $\mathfrak{E}$.

   (e) Two distinct planes $\mathfrak{E}_1$ and $\mathfrak{E}_2$ have either no point in common or else a line in common, because the existence of a common point $A$ implies the existence of a further common point $B$, by I.7. Then by I.6 the connecting line $g$ of $A$ and $B$ lies entirely in $\mathfrak{E}_1$ and $\mathfrak{E}_2$. If $\mathfrak{E}_1$ and $\mathfrak{E}_2$ also have a common point $P$ not in $g$, then $\mathfrak{E}_1$ and $\mathfrak{E}_2$ must coincide by consequence (b).

## 1.2 Order axioms

1. The axioms of this section describe the "betweenness" concept.

   II.1 When the point $B$ lies between $A$ and $C$, then $A, B, C$ are three different points on a line, and $B$ also lies between $C$ and $A$.

   The relation "between" is therefore symmetric with respect to the endpoints. We write $(ABC) \overset{\leftarrow}{\to} (CBA)$.

   II.2 For any two points $A, B$ on a line $g$ there is a point $C$ with $(ABC)$.

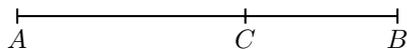
$A \qquad C \qquad B$

- II.3 Of any three points on a line, at least one lies between the others.

    In what follows, the "segment $AB$" means the set of all points $C$ with $(ABC)$, together with the points $A, B$ which we call the endpoints of the segment. All the other points of the connecting line of $A$ and $B$ are said to be "outside the segment $AB$". Three points on the same line are called "collinear"; three points not on a line are called "noncollinear".

- II.4 The plane order axiom. Let $A, B, C$ be noncollinear points and suppose that $g$ is a line in the plane they determine, but not passing through any of $A, B, C$. If $g$ goes through a point of the segment $AB$, then $g$ also goes through a point of the segment $AC$ or a point of the segment $BC$.

The idea of axiom II.4 is that a line which enters a triangle must also leave. It is provable that it cannot meet both the segments $AC$ and $BC$.
[]$^1$

2. **Consequences of the order axioms.**

    (a) Suppose we have a triangle$^2$ $A, B, C$ and a line with points $E, F, D$.

    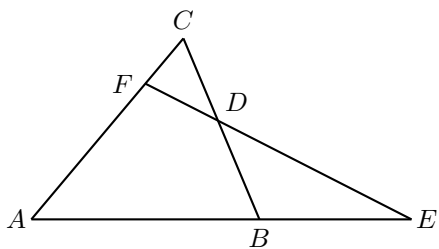

    Then $(ABE)$ and $(BDC)$ imply $(AFC)$, otherwise Axiom II.4 would be contradicted. This theorem is also known as the *axiom of Pasch*. It says that, if a line meets a side of a triangle at an interior point, and another side at an exterior point, then it necessarily meets the third side at an interior point.

    (b) **On Axiom II.2:**

    *The segment $AB$ is by definition the set of all points $C$ on the line connecting $A$ and $B$ for which $(ACB)$, together with the endpoints. We claim that segment $AB$ contains at least one point $C$ with $(ACB)$.*

---

$^1$Translator's note. Here I omit Moufang's sentence "In this chapter we denote the connecting line of points $P$ and $Q$ by $\overline{PQ}$, and the segment with endpoints $P$ and $Q$ by $PQ$" because I think it best to dispense with the bar notation. The *opposite* notation—$\overline{PQ}$ for the line segment and $PQ$ for the line—is used by Hilbert and also by Hartshorne (2000), and in any case Moufang drops the bar notation after a while. I have tried to minimize confusion by using $PQ$ for the line through $P$ and $Q$ and saying "segment $PQ$" explicitly when the segment is meant.

$^2$Translator's note. From now on, Moufang denotes the triangle with vertices $A$, $B$, $C$ by the usual notation $\triangle ABC$.

Proof: An $E$ outside segment $AB$ exists by I.3;
$F$ on $EA$ exists with $(AEF)$ by II.2;
$G$ on $FB$ exists with $(FBG)$ by II.2.

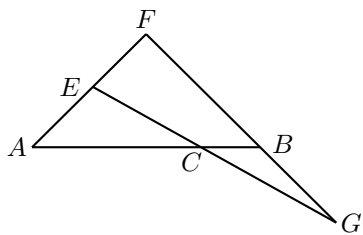

Now consider $\Delta ABF$ and the line $EG$. By Pasch's axiom (consequence (a)), $C$ exists as the intersection point of $AB$ with $EG$, and $(ACB)$.  Q.E.D.

Axiom II.2 ensures that a given segment may be extended beyond its endpoints; consequence (b) ensures that the segment itself is not empty of interior points. It has at least one, and hence infinitely many. It is thereby proved that the line connecting two points contains infinitely many further points.

(c) **On Axiom II.3:**

*Suppose that $A, B, C$ are collinear and that $(ACB)$ and $(CAB)$ are false. Then $(ABC)$ is true, that is, of three points on a line $g$, exactly one lies between the other two.*

Proof: A point $D$ outside $g$ exists by I.3;
an $E$ on $BD$ with $(BDE)$ exists by II.2.

Now consider $\Delta ABE$ and $DC$. Since $(ACB)$ is false and $(BDE)$ is true, we know by II.4 that there is an $H$ on segment $AE$ with $(AHE)$.

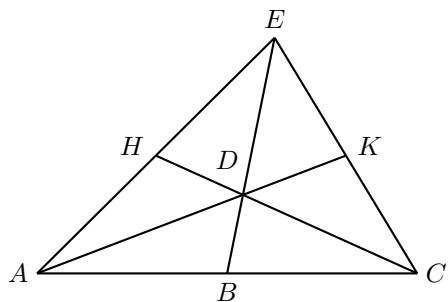

Next we consider $\Delta BCE$ and $AD$ and similarly conclude that there is a $K$ with $(EKC)$, because $(CAB)$ is false and $(BDE)$ is true. Then for $\Delta AKE$ and the line $HC$ it follows from $(EKC)$ and $(AHE)$ by Pasch's axiom (consequence (a)) that $(ADK)$. Finally we consider $\Delta AKC$ and $EB$. From $(ADK)$ and $(EKC)$ it follows, by Pasch's axiom, that $(ABC)$.  Q.E.D.

(d) **On Axiom II.4:**

*A line that meets all three sides of a triangle, two of them in the interior, necessarily meets the third side in the exterior. In other words: if $(AEB)$ and $(BFC)$ then $(AGC)$ is false.*

Proof. Let $E, F, G$ be collinear on $g$. Then either $(EFG)$ or $(FEG)$ or $(FGE)$. Suppose, for example, that $(EFG)$ holds.

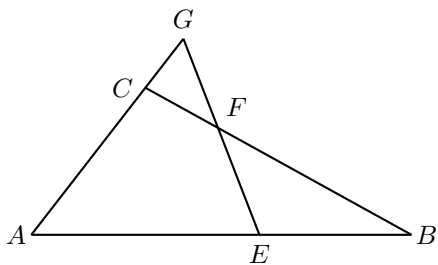

I consider $\Delta AEG$ and $BC$. It follows from the assumptions $(AEB)$ and $(BFC)$, by Pasch's axiom, that $(ACG)$ also holds. Then $(AGC)$ is false by (c).

The remaining cases are handled similarly. Q.E.D.

3. **Ordering of points in the plane.**

*Each line $g$ divides the points of the plane outside $g$ into two classes $\{A\}$ and $\{B\}$. Two points in the same class have a connecting segment that contains no point of $g$. On the other hand, the segment connecting two points not in the same class contains a point of $g$.*

Proof. Outside $g$ there is at least one point $A$, by I.3. Let $S$ be a point on $g$, so there is an $A_1$ on $SA$ such that $(SAA_1)$. Hence $(ASA_1)$ is false. Thus $A_1$ and $A$ belong to the class of $A$, since $(A_1SA)$ is likewise false.

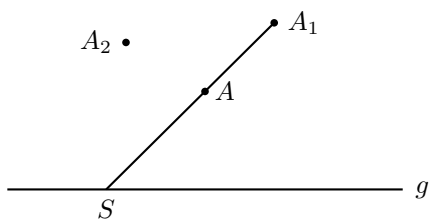

This "common membership" relation between $A$ and $A_1$ is symmetric and transitive. The symmetry is evident.

Proof of transitivity.

Suppose $A_2 \in \{A\}$ and $A_1 \in \{A\}$. Then segment $A_1A_2$ contains no point of $g$, because $g$ meets $\Delta AA_1A_2$ outside the segment $AA_1$ by the proof above. Likewise for segment $AA_2$. Hence $g$ also cannot meet segment $A_1A_2$ in its interior, by II.4, so $A_2 \in \{A_1\}$.

Conversely: by hypothesis there is a $U$ on $g$ with $(AUB)$. To show that $B \notin \{A_1\}$ we consider $\Delta AA_1B$ and $g$. Then $g$ meets segment $AB$ in the interior, but not segment $AA_1$. Hence $(A_1VB)$ by Pasch's axiom. That is, $A_1 \notin \{B\}$.

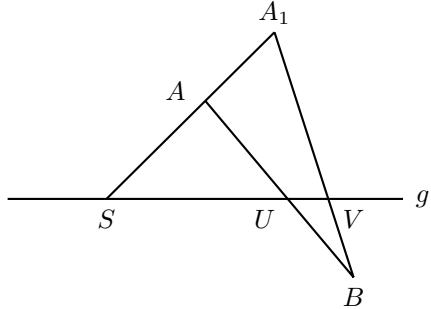

Each point therefore belongs to at most one class, but to at least one. Indeed, if $A$ and $B$ belong to different classes, and if $X$ is an arbitrary point on $AB$, then exactly one of $(BXU)$, $(UXA)$, $(XBU)$, $(XAU)$ holds. If $X$ is not on $AB$, then it follows from Pasch's axiom for $\Delta ABX$ and the line $g$ that $g$ meets either segment $BX$ or segment $AX$ in the interior, since it meets segment $AB$ in the interior.                              Q.E.D.

4. **Ordering of points in space.**

   *Each plane $\mathfrak{E}$ divides the points of space outside $\mathfrak{E}$ into classes, such that the segment connecting two points in different classes contains a point of $\mathfrak{E}$, while the segment connecting two points in the same class does not.*

   Proof. There is a point $A$ outside $\mathfrak{E}$ by I.8. In $\mathfrak{E}$ there is a point $S$ by I.4. On $SA$ there is a point $A_1$ with $(SAA_1)$ by II.2, so $(ASA_1)$ is false, that is, $A_1 \in \{A\}$.

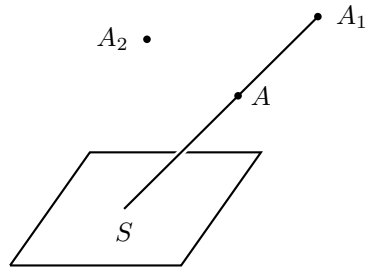

Now suppose that $A_1 \in \{A\}$ and $A_2 \in \{A\}$ and that $A, A_1, A_2$ are not collinear. It follows that $A_2 \in \{A_1\}$, that is, $A_1A_2$ contains no point of $\mathfrak{E}$. Indeed the plane $AA_1A_2 = \mathfrak{E}_1$ is uniquely determined, and since $\mathfrak{E}$ and $\mathfrak{E}_2$ have $S$ in common they have a common line $s$ by I.7 and I.6. The line $s$ meets $\Delta AA_1A_2$ in the interior of neither segment $AA_1$ nor segment $AA_2$, and hence also does not meet segment $A_1A_2$. Thus segment $A_1A_2$ contains no point of $\mathfrak{E}$.

Conversely, it remains to show that, when segment $AB$ contains a point of $\mathfrak{E}$, then so does segment $A_1B$.

Well, the plane through $A, A_1, B$ has the point $U$ in common with $\mathfrak{E}$, and hence a whole common line $t$. The line $t$ meets segment $AB$ in the interior, but not segment $AA_1$, hence it meets segment $A_1B$ in the interior, say at $V$, so $(A_1VB)$ holds. Thus $V$ also lies on $\mathfrak{E}$, and hence segment $A_1B$ contains a point of $\mathfrak{E}$.                                Q.E.D.

5. **Transitivity of the ordering of points on a line.**

    *Let $A, B, C, D$ be four points on a line. Then $(ABC)$ and $(BCD)$ imply $(ABD)$ and $(ACD)$.*

    We give the proof only for $(ACD)$.

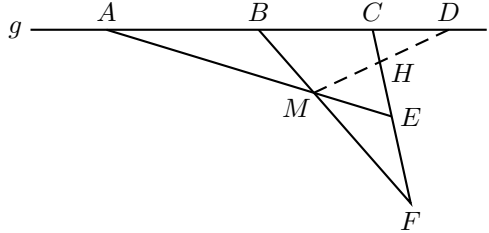

Outside $g$ there is a point $E$. On segment $CE$ there is an $F$ with $(CEF)$ by II.2. I apply Pasch's theorem to $\Delta ACE$ and $BF$. We have $(ABC)$ and $(CEF)$, hence there is an $M$ with $(AME)$. By the same argument in $\Delta BCF$ and the line $AE$, $(CEF)$ and $(ABC)$ give $(BMF)$. Similarly with $\Delta BDM$ and $FC$: $(BCD)$ and $(BMF)$ give $(MHD)$. Finally, in $\Delta ADM$ and $FC$, $(AME)$ and $(MHD)$ give $(ACD)$.           Q.E.D.

To prove $(ABD)$ one makes appropriate changes in the construction.

We also have:

$(ABC)$ *and* $(ACD)$ *imply* $(BCD)$ *and* $(ABD)$.

We prove the first assertion [$(ABC)$ *and* $(ACD)$ *imply* $(BCD)$].

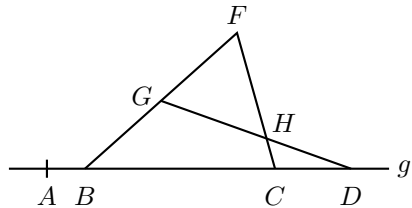

Outside $g$ there is a point $G$, and on $BG$ there is an $F$ with $(BGF)$. Applying Pasch's axiom to $\Delta ABG$ and $CF$, we find that $CF$ does not meet the segment $AG$ in its interior. It follows, by Pasch's axiom for $\Delta AGD$ and $CF$, and $(ACD)$, that $FC$ meets $GD$ at a point $H$ with $(GHD)$. Now Pasch's axiom for $\Delta BGD$ and $FC$, together with $(BGF)$ and $(GHD)$, implies that $(BCD)$.                     Q.E.D.

6. **The concept of "half-line".**

   Just as space is divided into two classes by a plane, and the plane is divided into two classes by a line, so too is a line divided into half-lines or rays.

   On the line $g$ there exists a point $O$, by I.3, and further points $B, A_1, A_2$. Suppose that $(A_1 O A_2)$ is false but $(A_1 O B)$ is true. In other words, $A_1$ and $A_2$ lie on the same side of $O$, but $A_1$ and $B$ lie on different sides. We have to prove that $A_2$ and $B$ also lie on different sides of $O$. We have either the ordering $(O A_1 A_2)$ or the ordering $(O A_2 A_1)$. In both cases it follows from (5) [that is, transitivity]—since $(A_1 O B)$ is true by hypothesis—that $(B O A_2)$ or $(A_2 O B)$.                                                            Q.E.D.

7. **The harmonic scale**

   Finally we show that infinitely many points can be constructed on a line. *First we prove that the fourth harmonic point is different from initial points $A, B, C$.*

   Suppose that the ordering $(CAB)$ holds. There is an $E$ outside $g$; draw $EB$ and $EA$. Then $EB \neq EA$. On $BE$ there is an $F$ such that $(EFB)$, so $F$ does not lie on $AE$. (By Pasch's axiom for $\Delta EBA$ and $CE$, and because $(EGA)$ and $(CAB)$, the segment $CF$ meets the segment $AE$ at $G$ with $(AGE)$. By Pasch's axiom for $\Delta CAF$ and $GB$, and because $(CAB)$ and $(CGF)$, we also have $(AMF)$.)

   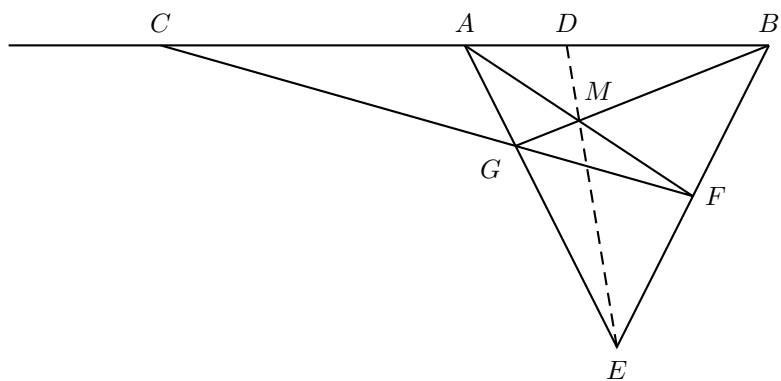

   We now apply Pasch's theorem four times:
   $\Delta CBF$ and $AE$: $(BFE)$ and $(CAB)$ imply $(CGF)$
   $\Delta ABE$ and $CF$: $(BFE)$ and $(CAB)$ imply $(AGE)$
   $\Delta CAF$ and $GB$: $(CGF)$ and $(CAB)$ imply $(AMF)$
   $\Delta ABF$ and $ME$: $(AMF)$ and $(BFE)$ imply $(ADB)$,
   where $D$ is the intersection of $ME$ and $AB$.

   From $(ADB)$ it also follows that $D \neq A$ and $D \neq B$. Since $(CAB)$ and $(ADB)$ hold, it follows from the transitivity of the ordering that $(CDB)$, hence $D \neq C$ also.                                                                        Q.E.D.

   This construction may now be repeated.

Again let $E$ be an arbitrary point outside $AB$, and take $F$ on $AE$ so that $(AFE)$. Then $FC$ and $EB$ intersect at a point $G$ with $(BGE)$ and $(FGC)$. One shows as above that $AG$ and $EC$ meet at a point $R$ with $(ERC)$, that $RB$ and $FC$ meet at a point $G_1$ with $(BG_1R)$, and that $EG_1$ meets the segment $BC$ at a point $B_1$. If one denotes the intersection[3] of the lines $UX$ and $VY$ by $UX \cap VY$, then the construction proceeds recursively as follows:

$$B_i R \cap FC = G_{i+1},$$
$$G_{i+1} E \cap AC = B_{i+1}.$$

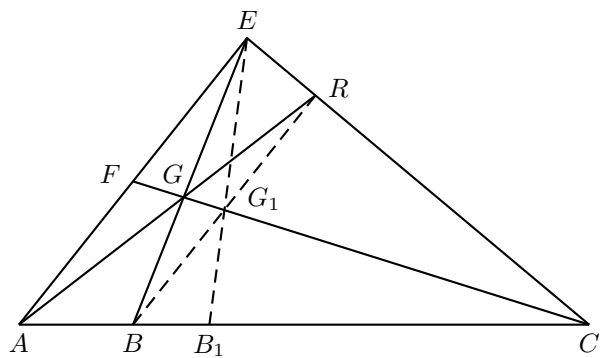

Then all the $B_i$ are different and they are in the order of their indices, that is, $(B_{i-1}B_iB_{i+1})$ holds, and also $(AB_iC)$ for all natural numbers $i$. Thus there are infinitely many distinct points $B_i$ on the segment $BC$. The proof depends on Pasch's axiom and the transitivity of order.

## 1.3   Congruence axioms

The congruence axioms involve relations between segments that are congruent or equal to each other.

III.1 **Transportability of segments.**

> *Suppose two points $A$ and $B$ are given on $g$, together with a point $A'$ on $g'$. Then on each half-line of $g'$ there is at least one point $B'$ such that the segment $AB$ is congruent or equal to the segment $A'B'$: $AB = A'B'$.*

One says that $AB$ and $A'B'$ are "congruent" or "equal" and writes $=$ or $\equiv$ or $\cong$ between them.[4] The order of the endpoints of a segment does not matter because the segment is defined by its set of endpoints. The definition of congruence depends symmetrically on the endpoints, hence $AB \cong A'B'$ and $A'B' \cong AB$. Also, as will be shown below, segment transport is unique.

---

[3]Translator's note. Moufang uses the symbol $\times$ for intersection. We use the ordinary intersection sign to avoid any possible confusion with a product operation.

[4]Moufang uses $=$ for congruence, but I have replaced this by $\cong$, since congruence does not mean absolutely equal. I use $=$ only between identical objects, though I follow ordinary prose usage in saying segments or angles are "equal" when they are equal in size. Moufang uses the $\equiv$ sign for "equals by definition", and also for "identically equals" and I have retained this notation.

**III.2 Transitivity of segment congruence.**

$A_1B_1 \cong AB$ and $A_2B_2 \cong AB$ imply $A_1B_1 \cong A_2B_2$.

It follows that each segment is congruent to itself.[5] By making the special choice $A_2 \cong A_1$ and $B_2 \cong B_1$, we get $A_1B_1 \cong A_1B_1$. The symmetry of segment congruence also follows, that is, $A'B' \cong AB$ implies $AB \cong A'B'$. Indeed III.2 says that $AB \cong AB$ and $A'B' \cong AB$ together imply $AB \cong A'B'$. Q.E.D.

Segment congruence is therefore a symmetric and transitive relation. One can say that two segments are "congruent to each other". The set of all segments is partitioned into classes of congruent segments.

**III.3 Additivity of segments**

*Let $A, B, C$ be three points on a line $g$, and let $A', B', C'$ be three points on a line $g'$.*

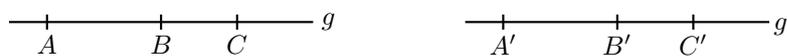

*Suppose that segments $AB$ and $BC$ have only the point $B$ in common, and that segments $A'B'$ and $B'C'$ have only the point $B'$ in common. Then $AB \cong A'B'$ and $BC \cong B'C'$ imply $AC \cong A'C'$.*

**Angle transport.**

We understand an *angle* to be a pair of rays[6] $\overline{g}, \overline{h}$ on two distinct lines $g, h$ with a common point (the vertex of the angle): $\angle(\overline{g},\overline{h}) \cong \angle(\overline{h},\overline{g})$. An angle divides the points of the plane into inner and outer points.

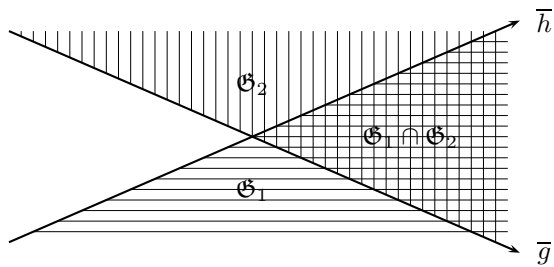

Suppose that $\mathfrak{G}_1$ is a half plane bounded by $\overline{h}$, and $\mathfrak{G}_2$ is a half plane bounded by $\overline{g}$. Then the common points of $\mathfrak{G}_1$ and $\mathfrak{G}_2$ comprise the interior of the angle. The connecting segment of two inner points meets neither $\overline{g}$ nor $\overline{h}$. Conversely: if $H$ is a point of $\overline{h}$ and $G$ is a point of $\overline{g}$, then the segment $HG$ consists entirely of inner points. The proof is analogous to that in Section 1.2, (3).[7]

---

[5]Translator's note. This property (reflexivity) should really be made an axiom, as was done by Euclid. Reflexivity does not follow from transitivity as Moufang claims, since the empty relation is transitive but not reflexive.

[6]Translator's note. I have gone along with Moufang's notation $\overline{g}$ for a half line of $g$. Some authors, such as Hartshorne (2000), use $\vec{g}$, because a half line has a direction. But Moufang also uses the bar notation for a half plane, where an arrow seems undesirable.

[7]Translator's note. Moufang here cites "consequence (e)", which does not exist. It seems likely that she means (3), which follows consequence (d) and is the place where the ordering of points in the plane is introduced.

III.4 **Axiom of angle transport (movability of figures).**

*Suppose $\angle(\overline{g},\overline{h})$ is given in $\mathfrak{E}$, and also $g'$ through $O'$ in $\mathfrak{E}'$. Let $\overline{\mathfrak{E}}'$ be a half plane of $\mathfrak{E}'$ determined by $g'$ and let $\overline{g}'$ be a half line of $g'$. Then there is exactly one half line $\overline{h}'$ originating from $O'$ and in the half plane $\overline{\mathfrak{E}}'$ such that $\angle(\overline{g},\overline{h}) \cong \angle(\overline{g}',\overline{h}')$.*

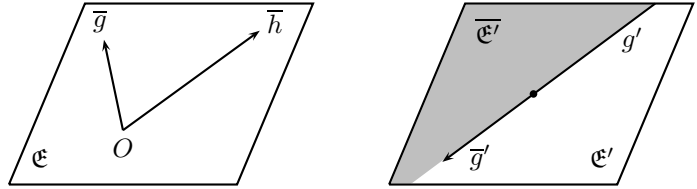

Here the sense of rotation of the angle is disregarded, just as the sense of a line segment was earlier. However, in contrast to segment transport, angle transport demands uniqueness as well as existence.

If $A$ is the vertex of an angle $\angle(\overline{h},\overline{k})$, $B$ point on $\overline{h}$, and $C$ a point on $\overline{k}$, then it is convenient to write $\angle(\overline{h},\overline{k})$ as $\angle BAC \cong \angle CAB$.

III.5 **Weak congruence theorem.**

*Given $\Delta ABC$ and $\Delta A'B'C'$, if $AB \cong A'B'$, $AC \cong A'C'$, and $\angle BAC \cong \angle B'A'C'$, then $\angle ABC \cong \angle A'B'C'$.*

It follows, by exchanging $B$ and $C$, that $\angle ACB \cong \angle A'C'B'$ also.

Now for some important *consequences* of the congruence axioms.

(a) **Uniqueness of segment transport.**

The uniqueness of segment transport follows from the uniqueness of angle transport and axiom III.5 by an indirect argument.

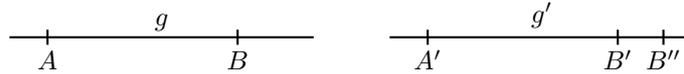

Suppose $AB \cong A'B'$ and $AB \cong A'B''$. Take a point $C$ outside $g$ and consider the triangles $\Delta A'B'C$ and $\Delta A'B''C$. We have

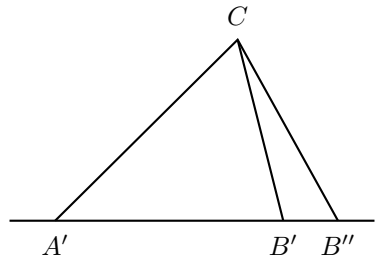

$$\angle CA'B' \cong \angle AA'B'',$$
$$A'C \cong A'C,$$
$$A'B' \cong A'B'',$$

hence by III.5 we have $\angle A'CB' \cong \angle A'CB''$. Then it follows from the uniqueness of angle transport that $B' = B''$. This proves the uniqueness of segment transport.

(b) **The complete congruence theorem for triangles.**[8]
*From $AB \cong A'B'$, $AC \cong A'C'$, $\angle CAB \cong \angle C'A'B'$ it follows that $BC \cong B'C'$.*

Proof. Suppose that $BC \not\cong B'C'$. Then there is a $D'$ on $B'C'$ such that $CB = B'D'$. Consider $\triangle ACB$ and $\triangle A'D'B'$. By III.5, all their corresponding angles are equal, so in particular

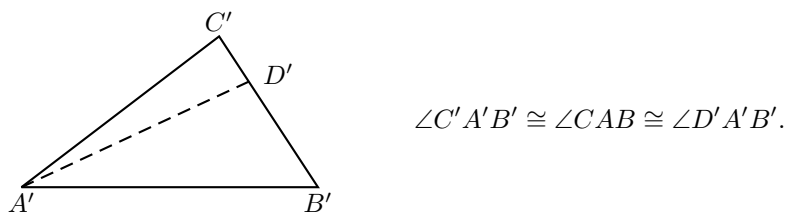

$$\angle C'A'B' \cong \angle CAB \cong \angle D'A'B'.$$

Since $D' \neq C'$, this contradicts the uniqueness of angle transport. Hence we must have $BC \cong B'C'$. Q.E.D.

(c) **Congruent angles have congruent supplementary angles.**
Suppose $AC \cong A'C'$, $AB \cong A'B'$, $AD \cong A'D'$, $\angle BAC \cong \angle B'A'C'$.

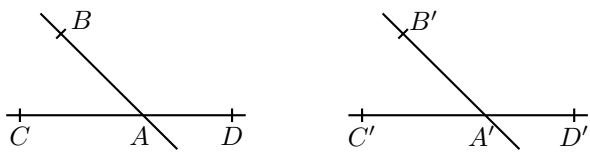

It is required to prove that $\angle BAD \cong \angle B'A'D'$.

Proof. Consequence (b), applied to $\triangle ABC$ and $\triangle A'B'C'$, gives $BC \cong B'C'$ and $\angle BCA \cong \angle B'C'A'$. By III.3, $CD \cong C'D'$.
Consequence (b), applied to $\triangle CBD$ and $\triangle C'B'D'$, now gives $BD \cong B'D'$ and $\angle BDC \cong \angle B'D'C'$. In $\triangle ABD$ and $\triangle A'B'D'$, Axiom III.5 gives $\angle BAD \cong \angle B'A'D'$. Q.E.D.

As an extension, we now prove *congruence of opposite angles:* $\alpha \cong \beta$.

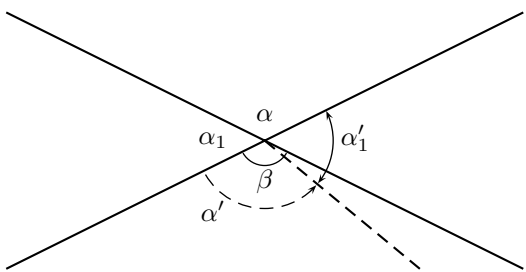

Proof. If $\alpha \not\cong \beta$ there is a unique $\alpha'$ such that $\alpha' \cong \alpha$. By the theorem just proved, applied to $\alpha$ and $\alpha'$, we have $\alpha_1 \cong \alpha'_1$. Applying the same theorem to $\alpha_1$ and $\alpha'_1$ then gives $\beta \cong \alpha'$, whence it follows by transitivity that $\alpha \cong \beta$. Q.E.D.

---
[8]Translator's note. This is the property known in English as SAS or "side-angle-side".

(d) **Theorem of the exterior angle.**

*An exterior angle of a triangle is greater than each non-adjacent interior angle.*

Proof. Make $AD \cong CB$, where $(DAB)$. Then we have to prove that $\angle CAD > \angle ACB$. First we prove that $\angle CAD \not\cong \angle ACB$.

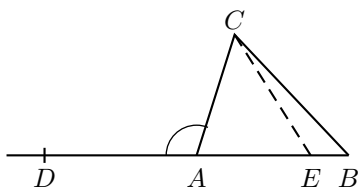

If these angles are equal, then applying Axiom III.5 to $\Delta DAC$ and $\Delta BCA$ gives $\angle DCA \cong \angle CAB$. By (c), $\angle DCA$ must equal the supplementary angle of $\angle BCA$. By the uniqueness of angle transport, $B$ therefore lies on $CD$ or, what is the same thing, $C$ lies on $DB$. Thus $C, A, B$ lie on a line, contrary to hypothesis.

Now we prove that the assumption $\angle DAC < \angle ACB$ also leads to a contradiction. Under this assumption there is an $E$ such that $\angle DAC \cong \angle ACE$, and $E$ lies on the same side of $A$ as $B$. But then $\Delta ACE$ is a triangle in which the exterior angle at $A$ equals the interior angle at $C$, which we have just proved to be impossible.

The theorem of the exterior angle is thereby proved for $\angle ACB$. The proof for $\angle CBA$ is analogous, with $B$ interchanged with $C$ and use of the opposite angle theorem.　　　　　　　　　　　　Q.E.D.

## 1.4 The parallel axiom

*Suppose we are given a plane $\mathfrak{E}$, a line $g$ in it, and a point $P$ in $\mathfrak{E}$ outside $g$. Then $\mathfrak{E}$ contains a line $h$ through $P$ not meeting $g$.*

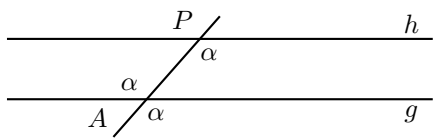

Proof. The existence of an intersection point $g \cap h = B$ in the figure contradicts the theorem of the exterior angle. For each line $PA$ one can find such a line $h$ through $P$ by transportation of the alternate angle $\alpha$.　　　　　　Q.E.D.

IV **Parallel axiom** (also known as the parallel postulate).

　　In $\mathfrak{E}$ there is at most one line through $P$ not meeting $g$.

　　Euclid formulated this postulate as follows. If two lines $g, h$ in the plane are cut by a third line $c$, and if the sum of the interior angles with $c$ on one side is less than two right angles, then $g$ and $h$ will meet on that side if prolonged sufficiently far.

Both formulations have the consequence that lines are parallel if and only if they have equal alternate angles.

The parallel postulate has played a vital role in the development of geometry. Even in ancient times, Euclid's formulation stimulated attempts to prove it from the remaining axioms. The noneuclidean geometries first found by Gauss, Bolyai, and Lobachevsky demonstrated its unprovability.

## 1.5 Axioms of continuity

V.1 **Axiom of Archimedes** (axiom of measure).

Let $AB$ and $CD$ be arbitrary segments, and let $CD$ be transported onto $AB$ repeatedly, giving points $A_1, A_2, \ldots, A_m$ where each $A_{m-1}A_m = CD$. Then for some natural number $m$ one has the ordering $(ABA_m)$.

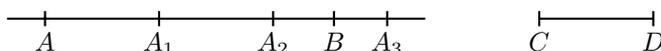

Here we have used the possibility and uniqueness of segment transport. A geometry in which axiom V.1 does not hold is called *non-Archimedean*.

With the axioms enumerated so far it is still not possible to establish the cartesian geometry over the real numbers. One needs to adjoin the following axiom.

V.2 **Axiom of linear completeness.**

The system of points on a line admits no order-preserving extension satisfying the first congruence theorem and the Archimedean axiom, that is, the axioms I.1–I.3, II, III.1, V.1. That is, it is not possible to adjoin further points to this system so that the extended system satisfies all the axioms just listed.

The completeness axiom is satisfiable essentially because the required axioms include the axiom of Archimedes. It possible to show that a system of points on a line satisfying Axioms I.1–I.3, II, III.1 can be extended in an order-preserving way. Thus, in the absence of the axiom of Archimedes, the completeness axiom can be refuted.

The linear completeness axiom implies the **completeness theorem:**

*The points, lines, and planes constitute a system that admits no extension satisfying the axioms I, II, II.1, V.1, and hence no extension satisfying all the axioms.*

Without the completeness axiom, that is, from Axioms I, II, III, IV alone, the identity of the geometry with the usual cartesian analytic geometry does not follow.

Regarding the independence of all these axioms, it has already been remarked that the spatial incidence axioms do not follow from the planar ones. The independence of the order axioms from the incidence axioms is proved with the help of a model geometry over the complex numbers. The congruence axioms are partially independent (namely, III.5 is) of the incidence and order axioms. The independence of the parallel postulate is due to the consistency of noneuclidean geometry, and that of the Archimedean axiom to the consistency of non-Archimedean geometry.

# Chapter 2

# Number systems

## 2.1 Number fields

To construct model geometries one often uses coordinate geometries, and hence algebraic methods. Algebra works with things we call numbers (and denote by small Latin letters), and uses two operations on them. We now consider the axioms of algebra.

1. I′ **Operations and existence.**

    I′.1 Uniqueness of sum: $a + b = c$ is unique.
    Here the order of summands is fixed.

    I′.2 If $a$ and $b$ are given, then there is exactly one $x$ such that $a + x = b$, and exactly one $y$ such that $y + a = b$.

    I′.3 Existence of "zero". There is exactly one element 0 such that $a + 0 = a$ and $0 + a = a$.

    I′.4 Uniqueness of product: $a \cdot b = c$ is unique.
    Also, the order of $a$ and $b$ is fixed.

    I′.5 If $a$ and $b$ are given, and $a \neq 0$, then there is exactly one $x$ such that $ax = b$ and exactly one $y$ such that $ya = b$.

    I′.6 Existence of "one". There is exactly one element 1 such that $a \cdot 1 = a$ and $1 \cdot a = a$.

    II′ **Computation rules.**

    II′.1 Associative law of addition: $(a + b) + c = a + (b + c)$.

    II′.2 Commutative law of addition: $b + a = a + b$.

    II′.3 Associative law of multiplication: $(ab)c = a(bc)$.
    One proves by complete induction that a sum of arbitrarily many summands, and also a product of arbitrarily many factors, may be bracketed arbitrarily.

    II′.4 Left distributive law: $a(b + c) = ab + ac$.

    II′.5 Right distributive law: $(b + c)a = ba + ca$.

II′.6 Commutative law of multiplication: $ab = ba$.

Each system of symbols with sum and product operations satisfying the laws I′ and II′ is called a *number field* (for example, the rational numbers, real numbers, complex numbers). Systems for which II′.3 is not satisfied are called *nonassociative systems*.

2. First we prove some *consequences* of axioms I′ and II′.

    (a) **Equality of right and left inverses.** By I′.5 there are $x$ and $y$ such that $ax = b$ and $ya = b$. Let $b$ be the unit element,[1] which we now denote by $e$.
    Claim: $ax = e, ya = e \to x = y$.
    Proof. $ax = e$ implies $(ax)a = ea = a$.
    $ya = e$ implies $a(ya) = ae = a$ or, by II′.3, $(ay)a = a = (ax)a$.
    Hence $(ay - ax)a = 0$.
    Then, since $a \neq 0$, it follows that $ay - ax = 0 = a(y-x)$ and therefore $y = x$. Q.E.D.
    We denote $y$ and $x$ by $a^{-1}$.
    One shows analogously that $a + x = 0$ and $y + a = 0$ imply $x = y$. We denote $x$ and $y$ by $(-a)$.

    (b) **Inverse of an inverse:** $(a^{-1})^{-1} = a$.
    Proof. Indeed, $x = (a^{-1})^{-1}$ is defined by $a^{-1}x = e$.
    Hence $a(a^{-1}x) = ae = a$, or $(aa^{-1})x = a$, so $ex = a = x$. Q.E.D.
    One shows analogously that $-(-a) = a$.

    (c) **Inverse of a product:** $(ab)^{-1} = b^{-1}a^{-1}$.
    Proof.
    $$\begin{aligned}
    (ab)^{-1} \cdot (ab)\} \cdot b^{-1} &= (ab)^{-1} \cdot \{(ab)b^{-1}\} \\
    &= (ab)^{-1} \cdot \{a(bb^{-1})\} \\
    &= (ab)^{-1} \cdot ae \\
    &= (ab)^{-1} \cdot a.
    \end{aligned}$$

    Hence $b^{-1} = (ab)^{-1}a$ and therefore
    $$\begin{aligned}
    b^{-1} \cdot a^{-1} &= \{(ab)^{-1} \cdot a\} \cdot a^{-1} \\
    &= (ab)^{-1} \cdot (aa^{-1}) \\
    &= (ab)^{-1} \cdot e \\
    &= (ab)^{-1} \quad \text{Q.E.D.}
    \end{aligned}$$

    One notes that the derivation of this relation uses the associative law of multiplication only in the special form where two factors of the product are equal.[2]

3. **The axioms I′ and II′ are not independent.**

---

[1] Translator's note. That is, the element called "one" in Axiom I′.6.
[2] Translator's note. Or rather, inverses of each other.

(a) The existence of zero follows from Axioms I'.1, I'.2, and II'.1.

Proof. By I'.2, the equation $a + x = a$ has exactly one solution $x_a = \bar{a}$, and the equation $y + a = a$ has exactly one solution $y_a = \bar{\bar{a}}$. Then II'.1 gives, for arbitrary $b$,

$$b + (a + x_a) = b + a = (b + a) + x_a.$$

Since $b + a = c$ is arbitrary (that is, for any $a$ we can choose $b$ so that $b + a$ has a prescribed value), it follows that

$$c = c + x_a,$$

that is, $x_a$ is a right zero element.

It also follows from II'.1 that

$$(y_a + a) + b = a + b = y_a + (a + b),$$

that is, $c = y_a + c$, so $y_a$ is a left zero.

It remains to show that $x_a = y_a$. Substituting $c = y_a$ in $c = c + x_a$ gives

$$y_a = y_a + x_a,$$

and it follows from $c = y_a + c$ with $c = x_a$ that

$$x_a = y_a + x_a,$$

hence $x_a = y_a$ follows by I'.1.

This proves the existence of at least one zero element $x_a = y_a = 0$ such that $c + 0 = 0 + c = c$ for all $c$. If there were another zero element $0'$, then for $c = 0'$ we should have

$$0' + 0 = 0 + 0' = 0',$$

and $d + 0' = 0' + d = d$ with $d = 0$ gives

$$0 + 0' = 0' + 0 = 0,$$

hence $0 = 0'$. Q.E.D.

The existence of the unit element is proved analogously from I'.4, I'.5, and II'.3.

(b) The commutative law of addition follows from I', II'1,4,5.

Proof. Let $u$ and $v$ be any two elements with $u \neq 0$.

$$\begin{aligned}(u + v)(u^{-1} + e) &= u(u^{-1} + e) + v(u^{-1} = e) \\ &= uu^{-1} + u + vu^{-1} + v \\ &= e + u + vu^{-1} + v, \\ (u + v)(u^{-1} + e) &= (u + v)u^{-1} + (u + v) \\ &= uu^{-1} + vu^{-1} + u + v \\ &= e + vu^{-1} + u + v.\end{aligned}$$

Now let $\tilde{e} + e = 0$ and $v + \tilde{\tilde{v}} = 0$, and it follows from the equations above by left addition of $\tilde{e}$ and right addition of $\tilde{\tilde{v}}$ that

$$\tilde{e} + e + u + vu^{-1} + v + \tilde{\tilde{v}} = \tilde{e} + e + vu^{-1} + u + v + \tilde{\tilde{v}},$$

hence
$$u + vu^{-1} = vu^{-1} + u.$$

Setting $vu^{-1} = w$ gives

$$u + w = w + u \qquad \text{Q.E.D.}$$

This proof does not carry over to a proof of the commutative law of multiplication, because the analogous distributive law $a + (bc) = (a+b)(a+c)$ does not exist!

## 2.2 Ordered and continuous number systems

III′ **Axioms of order**

III′.1 **Existence and uniqueness of order:** for any two elements $a$, $b$, exactly one of the relations $a \leq b$, $a = b$, $a \geq b$ holds.

III′.2 **Transitivity of order:** $a > b$ and $b > c$ imply $a > c$.

(And similarly for the sign $<$.) The order is therefore *antisymmetric* and *transitive*.

III′.3 **First monotonicity law** (monotonicity of addition).

If $a > b$ then $a + c > b + c$ and $c + a > c + b$ for any $c$.

III′.4 **Second monotonicity law** (monotonicity of multiplication).

If $a > b$ and $c > 0$ then $ac > bc$ and $ca > ba$.

A number system that satisfies Axioms III′.1 to III′.4 is called a linearly ordered number system. Correspondingly, a number field is said to be linearly ordered when it satisfies all the axioms I′, II′, and III′.1 to III′.4.

IV′ **Axioms of continuity**

IV′.1 **Axiom of Archimedes**

For any positive[3] $a$ and $b$ there is a natural number $n$ such that $na > b$.

IV′.2 **Completeness axiom:** It is impossible to add further elements to the system in an order-preserving manner so that the Axioms I′, II′, III′, and IV′.1 remain valid.

An *ordered Archimedean number field* satisfies Axioms I′ to IV′.1.
We remark that II′.6 follows from I′, II′.1–II′.5, III′, and IV′.1, but not without use of IV′.1!

---

[3]Translator's note. That is, $a, b > 0$.

## 2.3   Structure of ordered Archimedean skew fields

**Definition:** A *skew field* is a number system that satisfies Axioms I′ and II′ *without* II′.6. Thus $ab = ba$ does not necessarily hold in a skew field.

**Main theorem:** *Each ordered Archimedean skew field is a field.*
In other words, the commutative law of multiplication is a consequence of the remaining axioms, without the completeness axiom.

  The proof of the main theorem divides into four steps.

1. **Lemma:** *Each linearly ordered skew field $\mathfrak{K}$ contains a subsystem $\mathfrak{R}'$ isomorphic to the field of rational numbers.*

   Proof:

   $\alpha$) In $\mathfrak{K}$ there is a unit element $1'$ and a zero element $0'$. We denote the skew field unit by $1'$ to distinguish it from the 1 of the rational numbers, and the skew field zero by $0'$ similarly to distinguish it from the 0 of the rational numbers. One constructs
   $$1' + 1' = 2',$$
   and so on, recursively:
   $$n' + 1' = (n+1)'.$$
   Each natural number is thereby associated with an element of $\mathfrak{R}'$.
   $1' + x = 0'$ has exactly one solution, namely $x = (-1)'$. We set
   $$(-1)' + (-1)' = (-2)' \quad, \ldots, \quad (-n)' + (-1)' = (-n-1)'.$$
   It will now be proved that the elements $n'$ are isomorphic to the ring of integers, that is [4]
   $$(mn)' = m'n' \quad \text{and} \quad (m+n)' = m' + n'.$$
   First one shows by complete induction that for $n \geq 1$ we also have
   $$(-n)' + 1' = (-n+1)'.$$
   This relation is correct for $n = 1$. Induction shows that, for each natural number $n$,
   $$\begin{aligned}(-(n+1))' + 1' = (-n)' + (-1)' + 1' &= (-n)' \\ &= ((-n-1)+1)' = (-(n+1)+1)'\end{aligned}$$
   We similarly show that
   $$n' + (-1)' = (n-1)',$$

---
[4]Translator's note. A pencilled note in the manuscript at this point reads:
  Remark: These relations prove only a homomorphism! To prove that the homomorphism is an isomorphism requires use of the linear ordering!

because
$$(n+1)' + (-1)' = n' + 1' + (-1)' = n' = (n+1-1)' = ((n+1)-1)'.$$

In summary, we can say that the definitions
$$n' + 1' = (n+1)'$$
$$(-n)' + (-1)' = -(n+1)'$$

imply the relations
$$(-n)' + 1' = (-n+1)'$$
$$n' + (-1)' = (n-1)'.$$

β) For integers $m$ we have
$$(m \cdot 1)' = m' = m' \cdot 1'.$$

We wish to prove that also
$$(m \cdot (-1))' = m'(-1)' = -m'.$$

We have
$$m'(1' + (-1)') = m' \cdot 0' = 0' = m'1' + m'(-1)' = m' + m'(-1)',$$

hence
$$m'(-1)' = (-m)'.$$

If we now assume that $(-m)' = -m'$ holds for some fixed $m > 0$, then it follows from
$$(m+1)' + (-m-1)' = m' + 1' + (-m)' + (-1)'$$
$$= 1' + m' + (-m)' + (-1)'$$
$$= 1' + (-1)' + m' + (-m)' = 0'$$

that
$$(-m-1)' = -(m+1)'.$$

The induction hypothesis
$$(-m)' = -m'$$

is correct for $m = 1$ by definition of $(-1)'$, and hence in general, so
$$m'(-1)' = -m' = (-m)' = (m \cdot (-1))'.$$

γ) Now we show that all integers $m, n$ satisfy the relations
$$n' + m' = (n+m)' \qquad (1)$$
$$n' \cdot m' = (n \cdot m)' \qquad (2)$$

(1) is correct for each $n$ and $m = 1$.
For $m > 1$ we make the following induction step:

$$n' + (m+1)' = n' + m' + 1' = (n' + m') + 1'$$
$$= (n+m)' + 1 = (n + (m+1))'.$$

For $m < -1$ one proceeds from the true relation for $m = -1$,

$$n' + (-1)' = (n-1)'$$

and the induction hypothesis

$$n' + (-m)' = (n-m)'$$

to

$$n' + (-(m+1))' = n' + (-m)' + (-1)' = (n-m)' + (-1)'$$
$$= (n-m-1)' = (n-(m+1))'.$$

(2) is proved analogously.
From $(mn)' = m'n'$, which is valid for each integer $m$ and $n = 1$, it follows that

$$(m(n+1))' = (mn+m)' = (mn)' + m' = m'n' + m'$$
$$= m'(n' + 1') = m'(n+1)'.$$

From $(m(-n))' = m'(-n)'$, which is valid for each integer $m$ and $n = 1$, it follows that

$$(m(-n-1))' = (-mn-m)' = (-mn)' + (-m)'$$
$$= m'(-n)' + (-m)' = m'(-n)' + m'(-1)'$$
$$= m'((-n)' + (-1)') = m'(-n-1)'. \quad \text{Q.E.D.}$$

The isomorphism between the elements $n'$ and the integers is thereby proved.

δ) I now apply division to $m'$ and $n'$.
Suppose $a$ and $m'$ belong to $\mathfrak{K}$.
It follows from the distributive law that

$$am' = a(1' + 1' + \cdots + 1') = a + a + \cdots + a$$
$$= 1'a + 1'a + \cdots + 1'a = (1' + 1' + \cdots + 1')a = m'a,$$

hence
$$am' = m'a.$$

Now $m' \cdot x = 1'$ has exactly one solution, namely

$$x = (m')^{-1}.$$

Then we also have $xm' = 1'$. We denote $(m')^{-1}$ by $\frac{1}{m'}$.

Define fractions by $(m')^{-1} \cdot n' = n'(m')^{-1} = \frac{n'}{m'}$.

We must first show that the latter "fraction" can be simplified:

$$\frac{(n\lambda)'}{(m\lambda)'} = \frac{n'\lambda'}{m'\lambda'} = n'\lambda'(m'\lambda')^{-1} = n' \cdot \lambda' \cdot \lambda'^{-1} \cdot m'^{-1} = n'm'^{-1} = \frac{n'}{m'}.$$

The set of all elements $\frac{n'}{m'}$ is a subdomain $\mathfrak{R}'$ isomorphic to the rational numbers. This is proved by showing that the correspondence

$$\frac{n'}{m'} \mapsto \frac{n}{m}$$

satisfies

$$\frac{n'_1}{m'_1} + \frac{n'_2}{m'_2} \mapsto \frac{n_1}{m_1} + \frac{n_2}{m_2}, \qquad \frac{n'_1}{m'_1} \cdot \frac{n'_2}{m'_2} \mapsto \frac{n_1}{m_1} \cdot \frac{n_2}{m_2}.$$

Indeed, if one thinks of

$$\frac{n'_1}{m'_1} \quad \text{and} \quad \frac{n'_2}{m'_2}$$

with a common denominator $m$ then

$$\frac{n'_1}{m'} + \frac{n'_2}{m'} = n'_1 \cdot m'^{-1} + n'_2 \cdot m'^{-1}$$

$$= (n'_1 + n'_2)m'^{-1} = (n_1 + n_2)'m'^{-1} = \frac{(n_1 + n_2)'}{m'}.$$

And

$$\frac{n'_1}{m'_1} \cdot \frac{n'_2}{m'_2} = n'_1 \cdot (m'^{-1}_1) \cdot n'_2 \cdot (m'^{-1}_2)$$

$$= n'_1 n'_2 \cdot m'^{-1}_1 \cdot m'^{-1}_2 = (n_1 n_2)' \cdot (m_2 m_1)'^{-1} = \frac{(n_1 n_2)'}{(m_1 m_2)'}$$

This completes the proof of the lemma. One may therefore remove the accents from the elements of $\mathfrak{R}'$.

Since $\mathfrak{K}$ is assumed to be linearly ordered, the elements of $\mathfrak{R}'$ within it have the same linear ordering as the rational numbers.

2. **Lemma:** *For any element $\rho$ of the subsystem $\mathfrak{R}'$ (isomorphic to the field $\mathfrak{R}$ of rational numbers) and any element $a$ of $\mathfrak{K}$ one has*

$$\rho a = a\rho.$$

Proof. Suppose $\rho = \frac{n'}{m'}$. Then

$$a\rho = a\frac{n'}{m'} = a(n'm'^{-1}) = (an')m'^{-1} = (n'a)m'^{-1} = n'(am'^{-1})$$

$$= n'(m'a^{-1})^{-1} = ((m'a^{-1})n'^{-1})^{-1} = ((a^{-1}m')n'^{-1})^{-1}$$

$$= (a^{-1}(m'n'^{-1}))^{-1} = (m'n'^{-1})^{-1}a = (n'm'^{-1})a$$

$$= \frac{n'}{m'}a = \rho a \qquad\qquad\qquad\qquad \text{Q.E.D.}$$

3. **Lemma:** *For any two elements a and b of $\mathfrak{K}$ with*

$$a < b$$

*there is an element $\rho$ of $\mathfrak{R}'$ such that*

$$a < \rho < b.$$

*In other words, the rational numbers are dense in $\mathfrak{K}$.*

Proof. Suppose $0 < a < b$. Then the Archimedean axiom gives an $N$ such that

$$N \cdot 1 > a.$$

This means there is an $n$ with $n < a < n+1$.

We investigate two cases:

(1) $n + 1 < b$. Then $n < a < n + 1 < b$, so $n + 1 = \rho$.
(2) $n + 1 > b$. Then $n < a < b < n + 1$, so $0 < b - a < 1$ since $b - a < n + 1 - n = 1$.

In the latter case the Archimedean axiom gives an $m$ such that

$$m(b - a) > 1 \quad \text{or} \quad b - a > \frac{1}{m} \quad \text{for a natural number } m.$$

Thus

$$b > a + \frac{1}{m}.$$

Also, there is a natural number $k$ such that

$$\frac{k}{m} < a < \frac{k+1}{m},$$

from which it follows that

$$b > \frac{k}{m} + \frac{1}{m} = \frac{k+1}{m}, \quad \text{whence} \quad a < \frac{k+1}{m} < b.$$

Taking $\rho = \frac{k+1}{m}$, Lemma 3 is proved. Q.E.D.

**Proof of the main theorem.** Let $a, b > 0$ be elements of $\mathfrak{K}$. We shall derive a contradiction from $ab \neq ba$. By the monotonicity law for multiplication, $ab$ and $ba$ have the same sign, so we can assume that both are positive.

Suppose that $ab < ba$.

Then by Lemma 3 there is a $\rho$ such that $ab < \rho < ba$. Right multiplication of $ab < \rho$ by $a$ gives

$$(ab)a < \rho a, \quad \text{or} \quad aba < \rho a.$$

By Lemma 2, $\rho a = a\rho$, hence $aba < a\rho$.

Now multiplying $\rho < ba$ on the left by $a$ gives

$$a\rho < a(ba) = aba.$$

But this is a contradiction, hence $ab = ba$. Q.E.D.

## 2.4 Ordered number fields: Archimedean and non-Archimedean

1. **The rational numbers** are an Archimedean ordered field. They satisfy all the axioms except the completeness axiom IV′.2. The field is extendible by Dedekind cuts to the *field of real numbers*, in which IV′.2 holds.

2. **A non-Archimedean field: rational functions of one parameter.**

   Let $r_i$ be numbers from the ground field of reals. We construct polynomials with parameter $t$:
   $$P(t) = \sum_{i=0}^{n} r_i t^i.$$

   $P(t)$ is nonzero when not all the coefficients $r_i$ are equal to 0.

   Expressions of the form
   $$f(t) = \frac{P_1(t)}{P_2(t)} = \frac{r_0 + r_1 t + r_2 t^2 + \cdots}{\rho_0 + \rho_1 t + \rho_2 t^2 + \cdots}$$

   are the elements of a system which, as we shall show, is a linearly ordered but non-Archimedean field.

   First we define the zero and unit elements:

   Zero element: $r_0 = \cdots = r_n = 0$, not all $\rho_i = 0$

   Unit element: $n = m$, $r_i = \rho_i$.

   **Operations and computation rules:**

   Addition is performed by termwise addition of numerators, after the denominators have been made equal.

   Polynomials are added and multiplied in the usual way. Thus all the usual rules for computation with fractions hold, that is, the elements of the system constitute a field.

   $$\frac{P_1}{P_2} + \frac{Q_1}{Q_2} = \frac{R_1}{R_2}, \qquad \frac{P_1 \cdot Q}{P_2 \cdot Q} = \frac{P_1}{P_2}, \qquad \frac{P_1}{P_2} \cdot \frac{Q_1}{Q_2} = \frac{P_1 Q_1}{P_2 Q_2}$$

   $$\frac{P_1}{P_2} \cdot x = \frac{Q_1}{Q_2} \quad \longrightarrow \quad x = \frac{Q_1}{Q_2} \cdot \frac{P_2}{P_1}.$$

   **Ordering:**

   A rational function $f(t)$ has constant sign for sufficiently large $t > |t_0|$. If $f_1(t)$ and $f_2(t)$ are elements of the form $\frac{P(t)}{Q(t)}$ then $f_1(t) - f_2(t)$ has constant sign for all sufficiently large positive values of $t$. One now sets:

   III′.1  $f_1(t) >, =, < f_2(t)$ according as $f_1(t) - f_2(t) >, =, < 0$ [for $t$ sufficiently large.]. This orders the rational functions in one indeterminate according to their ultimate behavior.

III$'$.2 This ordering is *transitive*, that is, $f_1 > f_2$ and $f_2 > f_3$ implies $f_1 > f_3$, because

$$(f_1(t) - f_3(t)) = (f_2(t) - f_3(t)) + (f_1(t) - f_2(t)).$$

The left hand side is $> 0$ for sufficiently large $t > |t_0|$ because each term on the right hand side is.

III$'$.3 **Monotonicity theorem for addition:**

If $f_1 - f_2 > 0$ for sufficiently large $t$ then so is

$$(f_1 + f_3) - (f_2 + f_3) = f_1 - f_2,$$

hence $f_1 > f_2$ implies $f_1 + f_3 > f_2 + f_3$.

III$'$.4 **Monotonicity theorem for multiplication:**

If $f_1 > f_2$ and $f > 0$ then $f_1 f > f_2 f$ or $(f_1 - f_2)f > 0$.

This is because

$$f_1 - f_2 = \frac{r_0 + r_1 t + \cdots + r_n t^n}{\rho_0 + \rho_1 t + \cdots + \rho_m t^m} > 0 \quad \text{implies} \quad \frac{r_n}{\rho_m} > 0$$

and

$$f = \frac{r'_0 + r'_1 t + \cdots + r'_k t^k}{\rho'_0 + \rho'_1 t + \cdots + \rho'_l t^l} > 0 \quad \text{implies} \quad \frac{r'_k}{\rho'_l} > 0.$$

In $(f_1 - f_2)f$ the highest power of $t$ in the numerator is in the term $r_n r'_k t^{k+n}$, while the highest power of $t$ in the denominator is in the term $\rho_m \rho'_l t^{m+l}$. Here we have

$$\frac{r_n r'_k}{\rho_m \rho'_l} > 0, \quad \text{that is,} \quad (f_1 - f_2)f > 0. \qquad \text{Q.E.D.}$$

**The Archimedean axiom does not hold.** Indeed, consider the element $f(t) = n - t$ ($n$ a natural number). Then, for fixed $n$, $f(t) < 0$ for sufficiently large positive $t$. Hence the Archimedean axiom does not hold because $f(t) = n - t$ gives $n \cdot 1 < t$. One says that the indeterminate $t$ is "non-Archimedeanly large relative to 1", that is, it is larger than each integer multiple of 1. One writes $t \gg 1$.

**Remark.** This extension of the real numbers to the field of functions $f(t)$ does not preserve *all* the axioms I$'$, II$'$, III$'$, and IV$'$.1 (to IV$'$.2). This is why the extension does not contradict the completeness axiom.

3. Another representation of the field of rational functions $f(t)$ is in terms of *Laurent series*.

$$f(t) = \sum_{i=N}^{\infty} \alpha_i t^i \quad \text{with } N \text{ an integer and the } \alpha_i \text{ real.}$$

We define the *operations* by

$$f_1 + f_2 = \sum_{i=N}^{\infty}(\alpha_i + \beta_i)t^i$$

$$f_1 \cdot f_1 = \left(\sum_{i=N}^{\infty}\alpha_i t^i\right)\left(\sum_{k=M}^{\infty}\beta_k t^k\right)$$

$$= \sum_{l=N+M}^{\infty}\gamma_l t^l$$

with $\gamma_{j+M+N} = \alpha_N \beta_{M+j} + \alpha_{N+1}\beta_{M+j-1} + \cdots + \alpha_{N+j}\beta_M$ for $j = 0, 1, 2, \ldots$

*Zero element*: $\alpha_i = 0$ for all $i$.

*Unit element*: $\alpha_0 = 1$, $\alpha_i = 0$ otherwise.

**Axioms I$'$ and II$'$** are satisfied, as can be seen immediately. For example, inverses arise as follows. Given

$$a = \sum_{i=N}^{\infty}\alpha_i t^i, \quad \text{we seek} \quad b = \sum_{k=M}^{\infty}\beta_k t^k$$

such that $ab = 1$.

One finds that $N + M = 0$, $M = -N$, $\alpha_N \beta_{-N} = 1$, so $\beta_{-N} = \frac{1}{\alpha_N}$.

Then $0 = \gamma_1 = \alpha_N \beta_{-N+1} + \alpha_{N+1}\beta_{-N}$ gives $\beta_{-N+1}$, and $\beta_{-N+l}$ is found recursively from

$$0 = \alpha_N \beta_{-N+l} + \alpha_{N+1}\beta_{-N+l-1} + \cdots + \alpha_{N+l}\beta_{-N}.$$

This gives all the $\beta_k$, and hence $b$.

The existence of inverses implies the existence of quotients with the help of the associative law of multiplication.

**Order axioms III$'$.**

Suppose $f(t) = \sum_{i=N}^{\infty}\alpha_i t^i = a$. When is $a >, =, < 0$? Since it is assumed that $a \neq 0$, there is a nonzero coefficient $\alpha_{N+n}$ with $N + n$ minimal.

We define $a >, =, < 0$ according as $\alpha_{N+n} >, =, < 0$.

Also, $a >, =, < b$ when $a - b >, =, < 0$.

All order axioms and the monotonicity laws are then satisfied, as one sees immediately.

The element $t = \sum \alpha_i t^i$ with $\alpha_1 = 1$, and the other $\alpha_i = 0$, is $> 0$.

**The Archimedean axiom is not satisfied.**

We consider $f(t) = 1 - n \cdot t$, which is $> 0$ (since $f(t) = 1 \cdot t^0 - nt^1$) for any natural number $n$. Thus $1 > nt$ and one writes $t \ll 1$. In contrast to the first representation, $t$ is now "non-Archimedeanly *small* relative to 1".

The convergence of the series in question is completely irrelevant. We focus only on the coefficients and define an element of the number system to be a denumerable sequence of real numbers,

$$a = (\alpha_N, \alpha_{N+1}, \alpha_{N+2}, \ldots),$$
$$b = (\beta_N, \beta_{N+1}, \beta_{N+2}, \ldots), \quad \text{or} \quad (\beta_M, \beta_{M+1}, \beta_{M+2}, \ldots)$$

and define

$$a + b = (\alpha_N + \beta_N, \alpha_{N+1} + \beta_{N+1}, \ldots)$$
$$a \cdot b = (\gamma_{N+M}, \gamma_{N+M+1}, \ldots)$$

where $\gamma_{N+M+l} = \alpha_{N+l}\beta_M + \alpha_{N+l-1}\beta_{M+1} + \cdots + \alpha_N \beta_{M+l}$

By carrying out the representation in this manner, the parameter $t$ is eliminated.

We mention that while the Archimedean axiom implies commutativity of multiplication, the converse does not hold, because—as we now see—there are non-Archimedean ordered fields with commutative multiplication.

4. **The complex numbers** constitute a field which cannot be linearly ordered, so *a fortiori* they cannot be Archimedeanly ordered.

   To prove that they cannot be ordered one proceeds indirectly, using the monotonicity laws to show that the assumptions $i > 0$ and $i < 0$ lead to contradictions.

   Suppose that $i > 0$ and take $i = a$, $0 = b$ in the second monotonicity law ($a > b$, $c > 0 \to ac > bc$). For $c = i$ we get

   $$ii > 0 \cdot i = 0,$$
   $$-1 > 0 \quad \text{False!}$$

   If we suppose that $i < 0$ then $-i > 0$ (otherwise addition of $i < 0$ and $-i < 0$ would give $i - i < 0$, that is, $0 < 0$, which is false).

   But $-i > 0$ multiplied by $-i > 0$ gives $-1 = (-i)(-i) > 0(-i) = 0$ which is false. This completes the proof that the complex numbers cannot be ordered.

5. **Coordinate geometry over a field** is constructed as usual by taking points to be number pairs (or number triples), lines (or planes, respectively) to be given by linear equations, and the axioms of geometry are then investigated algebraically.

   I   To verify the incidence axioms, it suffices to solve linear equations, which is possible in a field.

   II  The order axioms hold when the field is linearly ordered.

   III The congruence axioms are valid when square roots exist (since the length of $P_1 P_2 = \sqrt{(x_1 - x_2)^2 + (y_1 - y_2)^2}$) and motions of translation, reflection, and rotation are representable. Here we confine ourselves to the plane.

Unique transport of a segment onto a given half line $\overline{g}$ says that the circle around $A$ with given radius meets the segment $\overline{g}$ of the diameter in exactly two points. To determine their coordinates requires the square root operation.

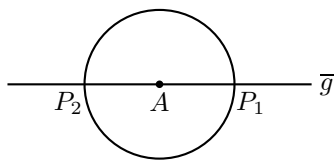

Angle congruence is defined analytically with the help of rotations and reflections. The transformations of the plane are represented as follows:

Translation: $x' = a + x, \quad y' = b + y$.

Reflection in the $x$-axis: $x' = x, \quad y' = -y$.

Rotation: $x' = \frac{a}{|\sqrt{a^2+b^2}|}x - \frac{b}{|\sqrt{a^2+b^2}|}y, \quad y' = \frac{b}{|\sqrt{a^2+b^2}|}x + \frac{a}{|\sqrt{a^2+b^2}|}y$

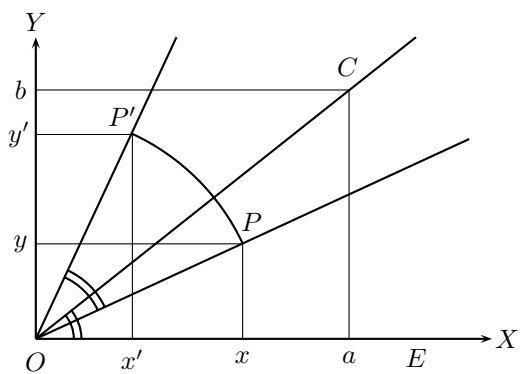

The latter formula represents the rotation through the angle $\angle COE$ that carries $P = (x,y)$ to $P' = (x',y')$. Then $\angle COE$ is called $\angle P'OP$.

Thus the introduction of the congruence axioms is possible in a number domain in which $\sqrt{a^2+b^2} = a\sqrt{1+\frac{a^2}{b^2}}$ is present along with $a$ and $b$. The field of rational functions of $t$ previously considered does not satisfy this condition until extended to the domain of special algebraic functions of $t$, obtained from the real numbers and $t$ by repeated rational operations and the operation $w \mapsto \sqrt{1+w^2}$. This extension may also be linearly ordered by ordering the algebraic functions according to their ultimate behavior. One sees in this way that the Axiom of Archimedes is not a consequence of Axioms I, II, III, IV.

  IV The continuity axioms of the geometry are equivalent to those of the number system.

By means of the plane geometries over the various number systems one obtains the following table, in which we do not further investigate the complex geometry.

| Field | Valid axiom groups |
|---|---|
| Real numbers | I, II, III, IV, V |
| Complex numbers | I; II false |
| Rational functions of $t$ | I, II, IV; III, V false |
| Algebraic functions of $t$ | I, II, III, IV; V false |

Consequences: II is independent of I, V.1 is independent of I to IV.

6. We conclude with **a further independence assertion:** *Axiom III.5 is independent of Axioms I, II, III.1–III.4, IV, V.1.*

   The proof is obtained through a suitable *model geometry* over the real numbers, in which points, lines, and planes are defined as usual but length is defined artificially as follows.

   $$\text{length } P_1 P_2 = \left| \sqrt{(x_1 - x_2 + y_1 - y_2)^2 + (y_1 - y_2)^2 + (z_1 - z_2)^2} \right|$$

   First we prove the uniqueness of segment transport.

   The parametric equation of a line is

   $$x = x_0 + at$$
   $$y = y_0 + bt$$
   $$z = z_0 + ct$$

   so

   $$(\text{length } P_0 P)^2 = l^2 = (at + bt)^2 + b^2 t^2 + c^2 t^2$$

   and therefore

   $$t = \pm \frac{l}{\sqrt{(a+b)^2 + b^2 + c^2}}.$$

   That is, the position of a point $P_0$ on a half line from $P$ is uniquely determined by its distance.

   With III.2 there is nothing to prove, since equality of numbers is transitive.

   Segment addition (Axiom III.3) is derived as follows. Suppose $(P_1 P_2 P_3)$ is the order of three points on a line. Then their parameter values satisfy $t_1 > t_2 > t_3$ (or else $t_1 < t_2 < t_3$). This is because

   $$\text{length } P_1 P_2 = \left| \sqrt{(a+b)^2 + b^2 + c^2} \right| (t_1 - t_2)$$
   $$\text{length } P_2 P_3 = \left| \sqrt{(a+b)^2 + b^2 + c^2} \right| (t_2 - t_3),$$

   hence

   $$\text{length } P_1 P_2 + \text{length } P_2 P_3 = \left| \sqrt{(a+b)^2 + b^2 + c^2} \right| (t_1 - t_3)$$
   $$= \text{length } P_1 P_3.$$

   Hence III.3 holds.

   The proof that III.5 does not hold is carried out by considering the triangle $ABC$ with $A = (1, 0)$, $B = (-1, 0)$, $C = (0, \frac{1}{\sqrt{2}})$.

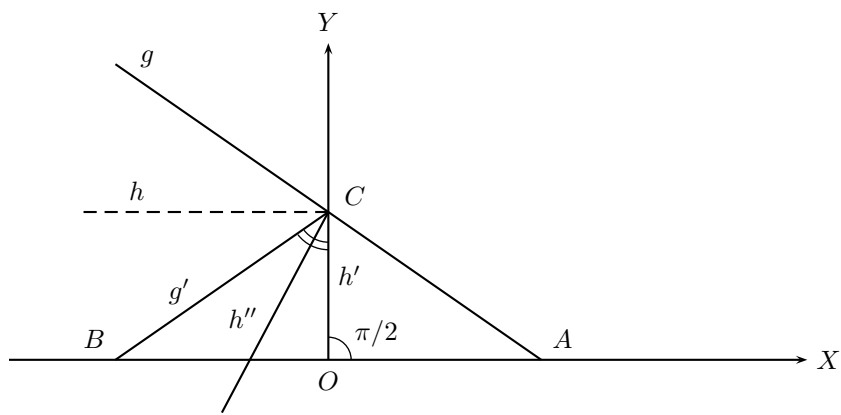

We have length $OA$ = length $OB$ = 1 and

$$\text{length } OC = \sqrt{\left(0 - 0 + 0 - \frac{1}{\sqrt{2}}\right)^2 + \left(0 - \frac{1}{\sqrt{2}}\right)^2 + (0-0)^2} = 1$$

$$\text{measure } \angle AOC = \text{measure } \angle BOC = \frac{\pi}{2}.$$

Thus triangles $OAC$ and $OBC$ satisfy the hypotheses of III.5. Nevertheless $\angle OAC \not\cong \angle OCB$, because if one translates angle $\alpha = \angle CAO$ to the apex $C$ and rotates it about $C$ to bring $g$ to the position of $h$, then $g'$ does not go to $h'$ but to $h''$.[5] We mention that the pseudolength we have introduced agrees with the ordinary length for all point pairs in which $y_1 = y_2$.

## 2.5 Skew fields and Desarguesian geometry

1. A skew field is defined by the operations and computation rules for a field, without the commutative law of multiplication. When $AB \neq BA$ for at least one pair of elements one speaks of a proper skew field. We first treat Hamilton's *quaternions* as an example of a proper skew field.

    The ground field is the real numbers and the basis elements, or units, are $e_0 = 1$, $e_1 = i$, $e_2 = j$, $e_3 = k$. The numbers of the system are

    $$a = r_0 + r_1 e_1 + r_2 e_2 + r_3 e_3.$$

---

[5]Translator's note. Essentially the same example occurs in Hilbert's 1930 *Grundlagen*, but without the transport and rotation of angle. Hilbert observes directly that $\angle OAC \not\cong \angle OCB$, despite SAS being satisfied in triangles $AOC$ and $COB$.

Hilbert also points out that

$$\text{length } AC = \sqrt{2 - \frac{2}{\sqrt{2}}} \neq \sqrt{2 + \frac{2}{\sqrt{2}}} = \text{length } BC$$

and that the isosceles triangle theorem fails for $\triangle AOC$ and $\triangle COB$.

*Addition* is defined by

$$a + b = (r_0 + \rho_0) + (r_1 + \rho_1)e_1 + (r_2 + \rho_2)e_2 + (r_3 + \rho_3)e_3,$$

and *multiplication* by

$$ab = \sum r_i e_i \rho_k e_k = \sum r_i \rho_k e_i e_k$$

and

$$e_i e_k = \alpha_{ik0} + \alpha_{ik1} e_1 + \alpha_{ik2} e_2 + \alpha_{ik3} e_3,$$

where the $\alpha_{ike}$ are given constants of the ground field. In our case we have in particular

$$\begin{array}{lll} e_1 e_2 = e_3 & e_2 e_3 = e_1 & e_3 e_1 = e_2 \\ e_2 e_1 = -e_3 & e_3 e_2 = -e_1 & e_1 e_3 = -e_2 \\ e_1 e_1 = -1 & e_2 e_2 = -1 & e_3 e_3 = -1. \end{array}$$

Quaternion multiplication *is* noncommutative; for example $e_1 e_2 = -e_2 e_1$. The distributive law of multiplication holds by definition. And the associative law holds because $(e_i e_j) e_k = e_i (e_j e_k)$ for all choices of indices, as one may verify by trying all combinations.

The *neutral elements* among the quaternions are:

"Zero": $r_0 = r_1 = r_2 = r_3 = 0$

"One", $e_0$: $r_0 = 1$, $r_1 = r_2 = r_3 = 0$, because $ae_0 = e_0 a = a$.

To determine the reciprocal of a quaternion, we must solve the equations

$$ax = e_0, \quad ya = e_0.$$

To do this one introduces the *conjugate* of $a$,

$$\overline{a} = r_0 - r_1 e_1 - r_2 e_2 - r_3 e_3,$$

for which

$$\begin{aligned} a\overline{a} &= r_0^2 + r_1^2 + r_2^2 + r_3^2 + (r_1 e_1)(-r_2 e_2) + (r_2 e_2)(-r_1 e_1) + \cdots \\ &= r_0^2 + r_1^2 + r_2^2 + r_3^2 - r_1 r_2 (e_1 e_2 + e_2 e_1) + \cdots \\ &= r_0^2 + r_1^2 + r_2^2 + r_3^2 \\ &= \text{"norm" of } a = N(a). \end{aligned}$$

Then the reciprocal of $a$ is

$$a^{-1} = \frac{\overline{a}}{N(a)} \quad \text{because} \quad aa^{-1} = e_0.$$

This establishes the validity of Axioms I′, II′.1–II′.5. II′.6 does not hold, so we have a proper skew field. Since the complex numbers $a = r_0 + r_1 i$ form a subsystem of the quaternions, the quaternions—like the complex numbers—cannot be linearly ordered.

The quaternions are an example of a so-called *hypercomplex algebra*. Such a structure is defined in general by a ground field of the $r_i$ (for example, the real or rational numbers) and basis elements $e_0, e_1, \ldots, e_n$. The general element of the algebra is $\sum_0^n r_i e_i$. Addition and multiplication of $a = \sum_0^n r_i e_i$ and $b = \sum_0^n r'_i e_i$ are defined by

$$a + b = \sum_0^n (r_i + r'_i) e_i$$

$$ab = \sum_{i,k=0}^n r_i r'_k e_i e_k, \quad \text{where} \quad e_i e_k = \sum_{l=0}^n \alpha_{ikl} e_l$$

for fixed $\alpha_{ikl}$ in the ground field. The $\alpha_{ikl}$ determine the so-called *multiplication table* of the algebra.

The condition for commutativity may be expressed by

$$\alpha_{ikl} = \alpha_{kil} \quad \text{for each } i, k, l.$$

For associativity it is necessary and sufficient that

$$(e_i e_k) e_n = e_i (e_k e_n) \quad \text{for each } i, k, n.$$

Now

$$(e_i e_k) e_m = \left( \sum_l \alpha_{ikl} e_l \right) e_m = \sum_s \left( \sum_l \alpha_{ikl} \alpha_{lms} \right) e_s$$

$$e_i (e_k e_m) = \sum_l \alpha_{kml} e_i e_l = \sum_s \left( \sum_l \alpha_{kml} \alpha_{ils} \right) e_s.$$

Hence the condition for associativity is that

$$\sum_l \alpha_{ikl} \alpha_{lms} = \sum_l \alpha_{kml} \alpha_{ils}.$$

The distributive law always holds. In fact if

$$a = \sum r_i^{(1)} e_i, \quad b = \sum r_i^{(2)} e_i, \quad c = \sum \rho_i e_i,$$

then

$$(a+b)c = \sum_i (r_i^{(1)} + r_i^{(2)}) e_i \cdot \sum_k \rho_k e_k$$

$$= \sum_l \left( \sum_{i,k} (r_i^{(1)} + r_i^{(2)}) \rho_k \cdot \alpha_{ikl} \right) e_l$$

$$= \sum_l \left( \sum_{i,k} r_i^{(1)} \rho_k \alpha_{ikl} + \sum_{i,k} r_i^{(2)} \rho_k \alpha_{ikl} \right) e_l$$

$$= ac + bc.$$

The proof of $a(b + c) = ab + ac$ is analogous.

The zero element of the algebra is given by $r_i = 0$ for all $i$.

The element $e_0$ is the unit element when

$$\alpha_{0kl} = \begin{cases} 0 & \text{for } k \neq l \\ 1 & \text{for } k = l. \end{cases}$$

The algebra is a *division algebra* when, for each nonzero element

$$a = \sum r_i e_i,$$

there is an element

$$a' = \sum \rho_l e_l$$

such that

$$aa' = e_0, \quad \text{the unit element.}$$

This requires solvability of the equations

$$\sum_{i,k=0}^{n} r_i \rho_k \alpha_{ikl} = \begin{cases} 0 & \text{for } l \neq 0 \\ 1 & \text{for } l = 0 \end{cases}$$

We content ourselves with these general remarks.

2. **Ordered skew fields (Desarguesian number systems)**

    The skew fields considered thus far are not ordered.

    The *Hilbert number system* we are about to describe is an example of a skew field that can be linearly ordered. The general number of this system is defined by

    $$a = \sum_{i=N}^{\infty} P_i(s) t^i,$$

    where $N$ is an integer, $P_i(s) = \sum_{k=N_i}^{\infty} \alpha_{ki} s^k$, the $N_i$ are integers, the $\alpha_{ki}$ are real, and $s$ and $t$ are parameters.

    *Addition* is defined by

    $$a + b = \sum_{i=N}^{\infty} (P_i(s) + Q_i(s)) t^i,$$

    where

    $$P_i(s) + Q_i(s) = \sum_{k=N_i}^{\infty} (\alpha_{ki} + \beta_{ki}) s^k.$$

    Because of this, addition is commutative and associative.

    The *zero element* has

    $$P_i(s) = 0, \quad \text{that is, } \alpha_{ki} = 0 \text{ for all } i, k.$$

    *Multiplication* is set up to be noncommutative.

First of all we have

$$ab = \sum_{i=N}^{\infty} P_i t^i \cdot \sum_{k=M}^{\infty} Q_k t^k$$
$$= P_N t^N \cdot Q_M t^M + (P_N t^N \cdot Q_{M+1} t^{M+1} + \cdots) + \cdots$$

We force noncommutativity by setting

$$ts = 2st.$$

We also set $t\alpha = \alpha t$ and $s\alpha = \alpha s$ for each rational $\alpha$, whence it follows by complete induction that

$$t^p s^q = 2^{pq} s^q t^p.$$

Then

$$t^N P_M(s) = t^N \sum_{l=NM}^{\infty} \alpha_{lM} s^l = \left(\sum_l \alpha_{lM} \cdot 2^{Nl} s^l\right) t^N \equiv \overline{P}_M^{(N)}(s) \cdot t^N.$$

In this notation, we get

$$ab = P_N \overline{Q}_M^{(N)} t^{N+M} + \left(P_N \overline{Q}_{M+1}^{(N)} + P_{N+1} \overline{Q}_M^{(M+1)}\right) t^{N+M+1} + \cdots,$$

which again represents an element of the system.

The associative law of multiplication holds. In fact, the equation

$$(ab)c = a(bc)$$

is satisfied termwise:

$$(\alpha_{ki} s^k t^i \cdot \beta_{nm} s^n t^m) \gamma_{pq} s^p t^q = (\alpha_{ki} \beta_{nm} \cdot 2^{in} s^{k+n} t^{i+m}) \gamma_{pq} s^p t^q$$
$$= \alpha_{ki} \beta_{nm} \gamma_{pq} \cdot 2^{in+(i+m)p} s^{k+n+p} t^{i+m+q},$$

and we get the same value for

$$\alpha_{ki} s^k t^i (\beta_{nm} s^n t^m \cdot \gamma_{pq} s^p t^q).$$

The distributive laws also hold, for example, $(a+b)c = ac + bd$ follows thus:

$$\left(\sum P_i t^i + \sum Q_i t^i\right) \cdot \sum R_k t^k = \sum_{i+k} (P_i + Q_i) \overline{R}_k^{(i)} t^{i+k}$$
$$= \sum \left(P_i \overline{R}_k^{(i)} + Q_i \overline{R}_k^{(i)}\right) t^{i+k}$$
$$= \sum_{i+k} P_i \overline{R}_k^{(i)} t^{i+k} + \sum_{i+k} Q_i \overline{R}_k^{(i)} t^{i+k}.$$

The *unit element* $e$ is defined by $P_0 = 1$, with the other $P_i = 0$, where $P_0(s) = 1$ means $\alpha_{00} = 1$ and the other $\alpha_{l0} = 0$.

We obviously have

$$ae = ea = a.$$

From now on we write 1 in place of $e$.

**Determination of inverses.**

We let
$$a^{-1} = \sum Q_k t^k, \quad \text{where} \quad Q_k(s) = \sum \beta_{lk} s^l$$
and the $\beta_{lk}$ are to be determined.

We have
$$1 = aa^{-1} = \sum_{i=N}^{\infty} P_i t^i \cdot \sum_{k=M}^{\infty} Q_k t^k = \sum_{i,k} P_i \overline{Q}_k^{(i)} t^{i+k}.$$

First we must have $N + M = 0$, so $M = -N$.

Then we must have $P_N(s) \cdot \overline{Q}_M^{(N)}(s) = 1$, so
$$\sum_{l=L}^{\infty} \alpha_{lN} s^l \cdot \sum_{k=K}^{\infty} \beta_{k,-N} \cdot 2^{Nk} s^k = 1.$$

Therefore $K = -L$ and $\alpha_{LN} \beta_{-L,-N} = 1$, whence $\beta_{-L,-N}$ is determined.

It now follows in turn that
$$\alpha_{L+1,N} \beta_{-L,N} 2^{-NL} + \alpha_{L,N} \beta_{-L+1,N} 2^{N(-L+1)} = 0,$$

and so on, giving a recursive determination of the $\beta_{k,-N}$.

The condition that the coefficient of $t^1$ equals 0 gives
$$P_N \overline{Q}_{-N+1}^{(N)} + P_{N+1} \overline{Q}_{-N}^{(N+1)} = 0$$

or
$$\sum_L \alpha_{lN} s^l \cdot \sum_{K_1} \beta_{k,-N+1} \cdot 2^{Nk} s^k + \sum_{L_1} \alpha_{l,N+1} s^l \cdot \sum_{-L} \beta_{k,-N} 2^{-N+1} s^k = 0,$$

where the second term is the given power series, beginning with $s^{L_1-L}$.

Thus we must have
$$K_1 + L = L_1 - L, \quad \text{or} \quad K_1 = L_1 - 2L$$
and
$$\alpha_{LN} \beta_{k,-N+1} + \alpha_{L_1,N+1} \beta_{-L,-N} = 0.$$

From this we get $\beta_{k,N+1}$ etc.

Since $ts = 2st$, this number system is a proper skew field. Its basis elements are the denumerably many symbols
$$s^k t^i.$$

Now we investigate the possibility of *ordering*. To this end we consider the general element
$$a = \sum_{i,k}^{\infty} \alpha_{k,i} s^k t^i$$

and the coefficient $\alpha_{MN}$ of the smallest power of $t$ and of $s$.

Then
$$a > 0 \quad \text{or} \quad a < 0$$
according as
$$\alpha_{MN} > 0 \quad \text{or} \quad \alpha < 0.$$
We say $a <, > b$ when $a - b >, < 0$.

The transitivity of the elements is inherited from the transitivity of the numbers
$$\alpha_{MN} \beta_{MN} \gamma_{MN}.$$

Monotonicity likewise follows from the monotonicity of the coefficients. The validity of the first monotonicity law,
$$a > b \quad \longrightarrow \quad a + c > b + c$$
follows from the validity of
$$\alpha_{MN} > \beta_{MN} \quad \longrightarrow \quad \alpha_{MN} + \gamma_{MN} > \beta_{MN} + \gamma_{MN}.$$

The validity of the second monotonicity law ($a > b$, $c > 0 \to ac > bc$) follows because
$$\alpha_{MN} > \beta_{MN}, \quad \gamma_{M_1 N_1} > 0$$
implies
$$\alpha_{MN} \cdot \gamma_{M_1 N_1} \cdot 2^{NM_1} > \beta_{MN} \cdot \gamma_{M_1 N_1} \cdot 2^{NM_1},$$
where the left side of the inequality represents the coefficient of
$$s^{M+M_1} t^{N+N_1} \quad \text{in} \quad ac.$$

This establishes the linear ordering. The ordering is *not Archimedean*:

We have $1 \gg s \gg t$, that is, $1 > ns$ and $s > mt$ for arbitrary natural numbers $n$ and $m$.

In fact
$$a = 1 - ns > 0,$$
because, more precisely,
$$1 \cdot s^0 - n \cdot s^1 > 0.$$
Likewise
$$a = s - mt > 0$$
because $a = st^0 - mt^1$, where $\alpha_{MN} = +1$ in the notation above.

3. **Application of Hilbert's number system to geometry.**

We construct a plane coordinate geometry over this skew field. A *point* is a pair $(x, y)$ of numbers; a *line* is represented by $ax + by + c = 0$.

Two points $P_1$ and $P_2$ are called equal when $x_1 = x_2$ and $y_1 = y_2$. Two lines $g_1$ and $g_2$ are called equal when there is a $\lambda$ such that
$$a_1 = \lambda a_2, \quad b_1 = \lambda b_2, \quad c_1 = \lambda c_2.$$

Then there is also a $\mu$ such that $a_2 = \mu a_1$, $b_2 = \mu b_1$, $c_2 = \mu c_1$. We call $g_1$ and $g_2$ *parallel* when $a_1 = \lambda a_2$, $b_1 = \lambda b_2$, and $c_1 \neq \lambda c_2$.

We now prove analytically the incidence theorems of this geometry. For the sake of later investigations we introduce a case distinction whose purpose will become apparent later.

**Existence of an intersection** $g_1 \cap g_2$ of nonparallel lines $g_1$ and $g_2$.

1) $a_1, a_2 \neq 0$.
   Then
$$a_1 x + b_1 y + c_1 = 0 \times a_1^{-1} \longrightarrow x + a_1^{-1} b_1 y + a_1^{-1} c_1 = 0 \quad (1)$$
$$a_2 x + b_2 y + c_2 = 0 \times a_2^{-1} \longrightarrow x + a_2^{-1} b_2 y + a_2^{-1} c_2 = 0 \quad (2)$$

   whence
$$(a_1^{-1} b_1 - a_2^{-1} b_2) y + (a_1^{-1} c_1 - a_2^{-1} c_2) = 0.$$

   If $a_1^{-1} b_1 - a_2^{-1} b_2 = 0$, or $a_1^{-1} b_1 = a_2^{-1} b_2$, then
$$\mu a_1 = a_2 \longrightarrow \mu b_1 = b_2,$$

   that is, $g_1 \parallel g_2$, which is not the case. Therefore
$$(a_1^{-1} b_1 - a_2^{-1} b_2) \neq 0$$

   and the coordinates of the intersection point are determined from (1) to be
$$\overline{y} = (a_1^{-1} b_1 - a_2^{-1} b_2)^{-1} \cdot (a_2^{-1} c_2 - a_1^{-1} c_1) \quad (3)$$
$$\overline{x} = -a_1^{-1} c_1 - a_1^{-1} b_1 (a_1^{-1} b_1 - a_2^{-1} b_2)^{-1} \cdot (a_2^{-1} c_2 - a_1^{-1} c_1) \quad (4)$$

   To show that $\overline{x}, \overline{y}$ also satisfy equation (2) one sets
$$\overline{\overline{x}} = -a_2^{-1} c_2 - a_2^{-1} b_2 (a_1^{-1} b_1 - a_2^{-1} b_2)^{-1} \cdot (a_2^{-1} c_2 - a_1^{-1} c_1) \quad (4')$$

   and then (4)−(4′) equals
$$-a_1^{-1} c_1 + a_2^{-1} c_2 + (-a_1^{-1} b_1 + a_2^{-1} b_2)(a_1^{-1} b_1 - a_2^{-1} b_2)^{-1} \cdot (a_2^{-1} c_2 - a_1^{-1} c_1) = 0,$$

   hence $\overline{x} = \overline{\overline{x}}$.                                                       Q.E.D.

   This shows that $g_1$ and $g_2$ have exactly one point of intersection.

2) $a_1 \neq 0$ but $a_2 = 0$, $b_2 = 0$.
   Equation (2) then becomes
$$b_2 y + c_2 = 0 \quad \text{or} \quad y = -b_2^{-1} c_2.$$

   Substituting this in (1) gives
$$a_1 x - b_1 b_2^{-1} c_2 + c_1 = 0 \quad \text{or} \quad x = -a_1^{-1} c_1 + a_1^{-1} b_1 b_2^{-1} c_2.$$

3) $a_1 = a_2 = 0$, $b_1 \neq 0$, $b_2 \neq 0$.
   This implies $a_1 = \lambda a_2$, $b_1 = \lambda b_2$, contrary to hypothesis.

### Does the parallel axiom hold?

Given line $g$ and point $P_1$ not on $g$, we seek $g_1$ parallel to $g$ through $P$.

Let $g_1$ have the equation

$$a_1 x + b_1 y + c_1 = 0.$$

Then [since $P_1 = (x_1, y_1)$ is on $g_1$]

$$a_1 x_1 + b_1 y_1 + c_1 = 0,$$

and therefore

$$c_1 = -a_1 x_1 - b_1 y_1.$$

Also [if $g$ is $a_0 x + b_0 y + c_0 = 0$] we must have

$$a_1 = \lambda a_0, \quad b_1 = \lambda b_0.$$

Thus $g_1$ has the equation

$$\lambda a_0 x + \lambda b_0 y - (\lambda a_0 x_1 + \lambda b_0 y_1) = 0,$$

which is unique up to the factor $\lambda$. Thus the parallel axiom holds.

As usual, one can associate an improper intersection point with two parallel lines and obtain the projective incidence theorems by introducing ideal elements. This can be done at once by means of homogeneous coordinates, but we have intentionally avoided this path.[6]

### Connecting two points.

For given $(x_1, y_1)$, $(x_2, y_2)$ we seek $\alpha$, $\beta$, $\gamma$ such that

$$\alpha x_1 + \beta y_1 + \gamma = 0$$
$$\alpha x_2 + \beta y_2 + \gamma = 0.$$

These equations determine $\alpha$, $\beta$, $\gamma$ up to a constant factor. The calculation is similar to that for the intersection of two lines.

**In space** there are analogous definitions:

Point: $P = (x, y, z)$.

Plane $\mathfrak{E}$: $Ax + By + Cz + D = 0$.

Planes $\mathfrak{E}_1$, $\mathfrak{E}_2$ are equal when $A_1 = \lambda A_2, \ldots, D_1 = \lambda D_2$.

A line is defined to be the set of common points of two distinct nonparallel planes.

We must now establish the *order* of points on a line and show that the Axioms II hold.

---

[6]Translator's note. Homogeneous coordinates are avoided, no doubt, because they do not work for the nonassociative systems considered in Section 2.7.

(a) Suppose that the line is parallel to an axis, say it is $by + c = 0$, so $y = -b^{-1}c = $ constant.

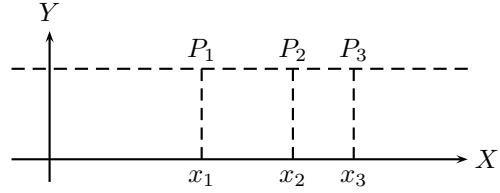

The $x_\nu$ are taken from a linearly ordered number system. Suppose that $x_i \neq x_j$ for $i \neq j$ and that $x_1 < x_2$. Then the compatible possibilities for ordering $x_1$, $x_2$, $x_3$ give the following schema.

$$x_1 < x_2 \begin{cases} x_2 < x_3 & \to & x_1 < x_2 < x_3 \\ x_2 > x_3 & \begin{cases} x_3 < x_1 & \to & x_3 < x_1 < x_2 \\ x_3 > x_1 & \to & x_1 < x_3 < x_2 \end{cases} \end{cases}$$

If $x_1 > x_2$ we similarly get

$$x_1 > x_2 > x_3$$
$$\text{or} \quad x_3 > x_1 > x_2$$
$$\text{or} \quad x_1 > x_3 > x_2.$$

**Definition:** $P_2$ lies between $P_1$ and $P_3$, in other words $(P_1 P_2 P_3)$, if

$$x_1 < x_2 < x_3 \quad \text{or} \quad x_1 > x_2 > x_3.$$

This establishes the possibility of ordering, and one easily sees that ordering is transitive.

(b) Suppose that the line is arbitrary: $ax + by + c = 0$ with $a, b \neq 0$. Then for two different points we have $x_i \neq x_j$ and $y_i \neq y_j$.

**Definition:** $(P_1 P_2 P_3)$ holds when $x_1 < x_2 < x_3$ or $x_1 > x_2 > x_3$.

We then have: $x_i < x_k < x_l$ or $x_i > x_k > x_l$ implies $y_i < y_k < y_l$ or $y_i > y_k > y_l$.

Proof. The equation $ax_i + by_i + c = 0$ gives $y_i = -b^{-1}(c + ax_i)$. Now it follows from

$$x_i < x_k < x_l,$$

by the monotonicity theorem for multiplication that

$$ax_i < ax_k < ax_l \quad \text{for} \quad a > 0$$
$$ax_i > ax_k > ax_l \quad \text{for} \quad a < 0,$$

hence by the monotonicity theorem for addition that

$$ax_i + c < ax_k + c < ax_l + c$$
$$\text{or} \quad ax_i + c > ax_k + c > ax_l + c,$$

and by the monotonicity theorem for multiplication again that

$$b^{-1}(ax_i + c) < b^{-1}(ax_k + c) < b^{-1}(ax_l + c)$$
$$\text{or} \quad b^{-1}(ax_i + c) > b^{-1}(ax_k + c) > b^{-1}(ax_l + c).$$

That is, $x_i < x_k < x_l$ or $x_i > x_k > x_l$ implies $y_i < y_k < y_l$ or $y_i > y_k > y_l$.                                                                    Q.E.D.

This shows that the linear order axioms hold. In order to investigate the axiom of Pasch (II.4), we must investigate the possibility of dividing the points of the plane into two classes by a line $g$. For precisely the points of $g$ we have $ax + by + c = 0$.

**Definition:**
Class 1 is defined by $g(x, y) = ax + by + c > 0$.
Class 2 is defined by $g(x, y) = ax + by + c < 0$.

Now we must prove:

1) When $(P_1 P_2 P_3)$ and $g(P_3) > 0$, $g(P_1) > 0$ then also $g(P_2)$.
   Proof. Suppose the connecting line $h$ of $P_1$ and $P_2$ is given by

   $$\alpha x + \beta y + \gamma = 0 \quad \text{with } \beta \neq 0,$$

   so that
   $$\alpha x_i + \beta y_i + \gamma = 0 \quad \text{for } i = 1, 2, 3.$$

   It follows that
   $$\alpha(x_1 - x_3) + \beta(y_1 - y_3) = 0$$
   $$\alpha(x_2 - x_1) + \beta(y_2 - y_1) = 0.$$

   The first of these equations times $(x_1 - x_3)^{-1}$ and the second times $(x_2 - x_1)^{-1}$ are
   $$\alpha + \beta(y_1 - y_3)(x_1 - x_3)^{-1} = 0$$
   $$\alpha + \beta(y_2 - y_1)(x_2 - x_1)^{-1} = 0.$$

   By subtraction, the latter give
   $$\beta(y_1 - y_3)(x_1 - x_3)^{-1} = \beta(y_2 - y_1)(x_2 - x_1)^{-1}.$$

   Since $\beta \neq 0$, it follows that
   $$(x_1 - x_3)^{-1}(x_2 - x_1) = (y_1 - y_3)^{-1}(y_2 - y_1) = \rho.$$

   Now if $(P_1 P_2 P_3)$ holds then
   $$x_1 < x_2 < x_3 \quad \text{or} \quad x_1 > x_2 > x_3$$

   and
   $$y_1 < y_2 < y_3 \quad \text{or} \quad y_1 > y_2 > y_3.$$

Hence in either case, $\rho < 0$. Also, $|\rho| < 1$ because

$$x_1 < x_2 < x_3 \quad \text{or} \quad x_1 > x_2 > x_3$$

implies

$$0 < x_2 - x_1 < x_3 - x_1 \quad \text{or} \quad 0 > x_2 - x_1 > x_3 - x_1,$$

hence

$$(x_3 - x_1)^{-1}(x_2 - x_1) < 1.$$

Now

$$x_2 = (x_1 - x_3)\rho + x_1$$
$$y_2 = (y_1 - y_3)\rho + y_1$$

hence

$$\begin{aligned}
g(P_2) &\equiv ax_2 + by_2 + c \\
&= a[(x_1 - x_3)\rho + x_1] + c\rho + b[(y_1 - y_3)\rho + y_1] - c\rho + c \\
&= (ax_1 + by_1 + c)\rho - (ax_3 + by_3 + c)\rho + ax_1 + by_1 + c \\
&= g(P_1)\rho - g(P_3)\rho + g(P_1) \\
&= g(P_1)(1 + \rho) - g(P_3)\rho \\
&> 0 \quad \text{since } g(P_1), g(P_2), 1 + \rho > 0, \rho < 0.
\end{aligned}$$

as required.                                                                                 Q.E.D.

Thus when $P_1$ and $P_2$ lie on the same side of $g$, so do all points of the segment $P_1P_3$.

2) *When $g(P_3) < 0$ and $g(P_1) > 0$ there is a $P_2$ on $g$ with $(P_1P_2P_3)$.*
Proof. Let

$$g \equiv ax + by + c = 0$$
$$h \equiv \alpha x + \beta y + \gamma = 0 \quad \text{with} \quad \alpha x_i + \beta y_i + \gamma = 0 \quad \text{for} \quad i = 1, 3.$$

In 1) it was proved that if $(P_1P_2P_3)$, where $P_3$ is any point between $P_1$ and $P_3$ on the line $P_1P_3$, then

$$(x_1 - x_3)^{-1}(x_2 - x_1) = (y_1 - y_3)^{-1}(y_2 - y_1) = \rho,$$

where $\rho < 0$ and $|\rho| < 1$, and also that

$$x_2 = (x_1 - x_3)\rho + x_1$$
$$y_2 = (y_1 - y_3)\rho + y_1.$$

One now seeks to determine $\rho$ so that $P_2$ at the same time lies on $g$, that is, so that

$$ax_2 + by_2 + c = 0.$$

This gives the following condition for $\rho$:

$$a(x_1 - x_3)\rho + ax_1 + b(y_1 - y_3)\rho + by_1 + c = 0,$$

or

$$g(P_2) = g(P_1)(1+\rho) - g(P_3)\rho = 0$$
$$[g(P_1) - g(P_3)]\rho = -g(P_1)$$
$$\rho = -[g(P_1) - g(P_3)]^{-1} \cdot g(P_1).$$

Hence $\rho < 0$ and also $|\rho| < 1$ because $0 < g(P_1) < g(P_1) - g(P_3)$.
Thus $P_2$ lies on $h$ and $g$. \hfill Q.E.D.

The results 1) and 2) show that a segment whose endpoints lie in the same class contains only points of this class, and points of different classes determine a segment containing a point of $g$. This enables us to prove the *theorem of Pasch*.

**Hypothesis:** $g$ meets $AB$ at $P$ with $(APB)$ and $CB$ at $Q$ with $(CQB)$.

Claim: *Either $g \parallel AC$ or $g$ has an intersection point $S$ with $AC$ and $(ASC)$ is false.*

Proof. When $g$ is parallel to $AC$ there is no intersection point, and hence the side $AC$ of the triangle certainly is not met by $g$ at an interior point.

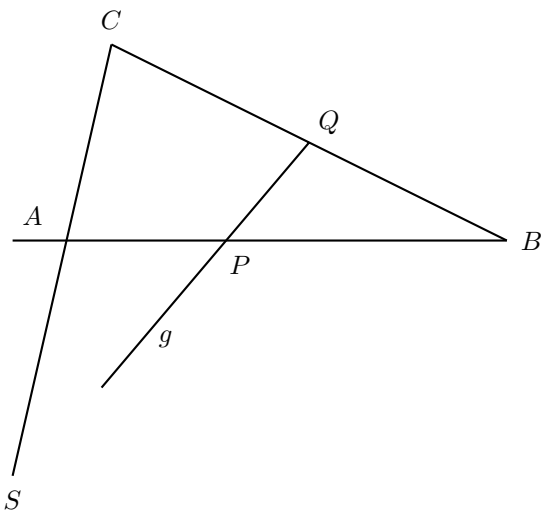

If an intersection point $S$ exists, then $(ASC)$ is false, otherwise $A$ and $C$ would lie on different sides of $g$, hence $B$ and $C$ would lie on the same side of $g$ (by the hypothesis that $A$ and $B$ lie on different sides), contrary to the hypothesis that $CB$ contains a point of $g$. \hfill Q.E.D.

**Summary:**

Over each ordered skew field it is possible to construct a coordinate geometry in which the planar (and spatial) incidence and order axioms hold. The theorem of Desargues also holds, as we shall show in 4 a).

4. **Two results about Desargues' theorem in the plane.**

   a) The following is a preparatory remark about linear dependence.
   Let
   $$\lambda_1 g_1 + \lambda_2 g_2 \equiv (\lambda_1 a_1 + \lambda_2 a_2)x + \cdots + (\lambda_1 c_1 + \lambda_2 c_2) = 0$$
   be a line through $P_0$, the intersection of $g_1$ and $g_2$. Given an arbitrary $\lambda_1 = \lambda$, $\lambda_2 = -\lambda g_1(P_1) \cdot g_2(P_2)^{-1}$ is determined by the requirement that the line go through $P_1$, if $P_1$ does not lie on $g_2$.
   Each line through $P_0$ can therefore be represented in the form

   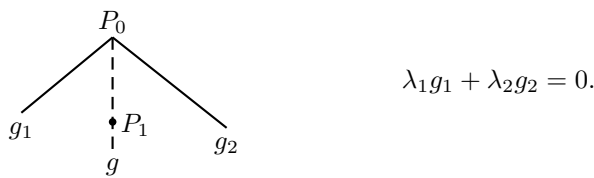

   $$\lambda_1 g_1 + \lambda_2 g_2 = 0.$$

   This includes the case where $P_0$ is an improper point.

   **Analytic proof of Desargues' theorem:**

   The theorem reads: *If two triangles lie in the plane so that the lines through corresponding vertices go through a point, then the intersections of corresponding sides lie on a line, and conversely.*

   It suffices to prove the second part of the assertion.

   The hypothesis is that

   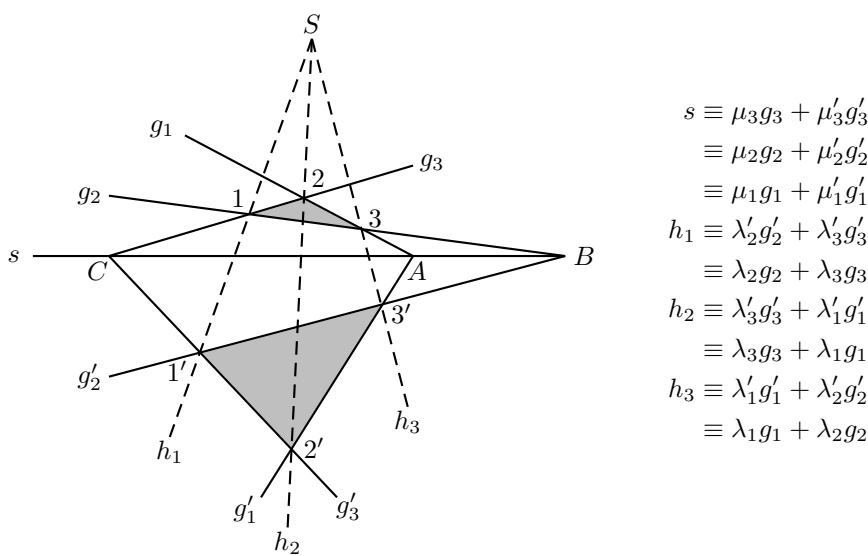

   $$\begin{aligned} s &\equiv \mu_3 g_3 + \mu_3' g_3' \\ &\equiv \mu_2 g_2 + \mu_2' g_2' \\ &\equiv \mu_1 g_1 + \mu_1' g_1' \\ h_1 &\equiv \lambda_2' g_2' + \lambda_3' g_3' \\ &\equiv \lambda_2 g_2 + \lambda_3 g_3 \\ h_2 &\equiv \lambda_3' g_3' + \lambda_1' g_1' \\ &\equiv \lambda_3 g_3 + \lambda_1 g_1 \\ h_3 &\equiv \lambda_1' g_1' + \lambda_2' g_2' \\ &\equiv \lambda_1 g_1 + \lambda_2 g_2 \end{aligned}$$

   And the conclusion is that $h_1$, $h_2$, $h_3$ go through a point $S$, that is, there are $c_1$, $c_2$, $c_3$ such that
   $$c_1 h_1 + c_2 h_2 + c_3 h_3 = 0 \quad \text{with not all } c_i = 0.$$

**Proof.**

By hypothesis,
$$\mu'_1 g'_1 - \mu'_2 g'_2 \equiv \mu_1 g_1 - \mu_2 g_2.$$

The line represented by the left hand side goes through the point $3'$, and the above identity says that it also goes through the point $3$, so the line is identical with $h_3$. Consequently one has

$$\mu'_2 g'_2 - \mu'_3 g'_3 \equiv \mu_2 g_2 - \mu_3 g_3 \equiv h_1$$
$$\mu'_3 g'_3 - \mu'_1 g'_1 \equiv \mu_3 g_3 - \mu_1 g_1 \equiv h_2$$
$$\mu'_1 g'_1 - \mu'_2 g'_2 \equiv \mu_1 g_1 - \mu_2 g_2 \equiv h_3$$

whence it follows by addition that

$$h_1 + h_2 + h_3 = 0$$

and hence the $c_i = 1$. \hfill Q.E.D.

This proof uses only the computation rules of a skew field. Thus the commutative law of multiplication is unnecessary. The converse statement, that perspective triangles also lie axially, now follows by applying the assertion just proved to the triangles

$$\Delta C11' \quad \text{and} \quad \Delta B33',$$

which are in perspective from $B$ by hypothesis. Since $h_1$ goes through $S$, as was just proved, the two triangles lie axially.

  b) The theorem of Desargues is not provable in a plane geometry satisfying the planar Hilbert axioms (I.1–3, II.1–4, III.1–4, IV, V) *except* III.5.

The proof of this assertion was given by Moulton (*Trans. Amer. Math. Soc.* 1902, vol. II, pp. 193–195), using a suitable model geometry.

Its points are defined as usual. Its lines (see the sketch) are defined as follows.

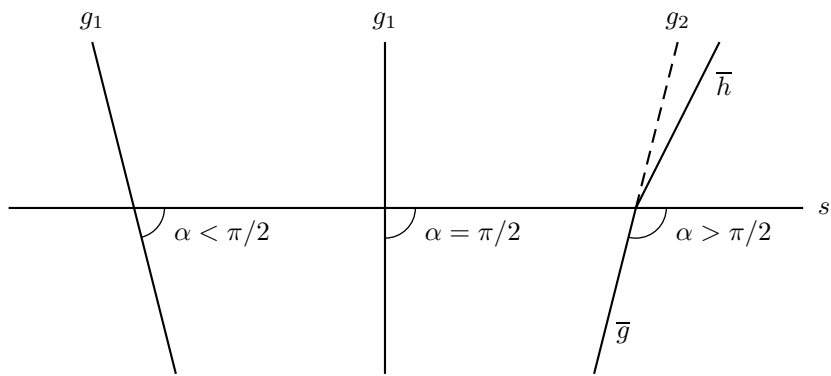

A fixed line $s$ is chosen, and an ordinary line $g_1$ that meets $s$ from below at an angle $\leq \frac{\pi}{2}$ is also a line of the artificial geometry. An ordinary line $g_2$ that meets $s$ from below at an angle $\alpha > \frac{\pi}{2}$ is bent above $s$ and continued at half the slope. The corresponding artificial line is thus a broken path consisting of two half lines $\overline{g}$ and $\overline{h}$, a so-called "bent line".

We must now establish the valid axioms.

The **connection of two points** $P$, $Q$ is seen to be unique by sketching the possible cases, of which only the last is nontrivial. One drops a perpendicular from $P$ to $s$ and makes $BP = PC$. Then if $CQ$ meets $s$ at $H$, the desired path is the bent line $PHQ$.

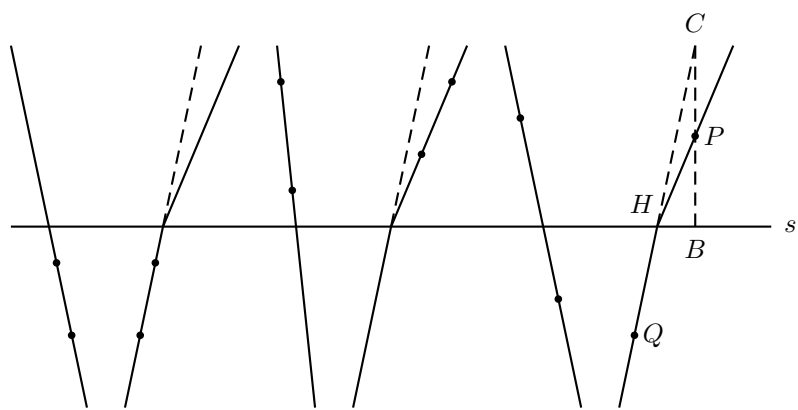

**Intersection of two lines:** Either $g \parallel h$ or $g \not\parallel h$.

$\alpha$) $g \parallel h$. For $\alpha \leq \frac{\pi}{2}$ the parallel postulate holds in the usual form. It also holds for $\alpha > \frac{\pi}{2}$, because when the lower half of a bent line $g$ is parallel to the lower half of $h$ this is also true for the upper half.

$\beta$) $g \not\parallel h$, that is, there is a point of intersection. The possible cases are sketched below.

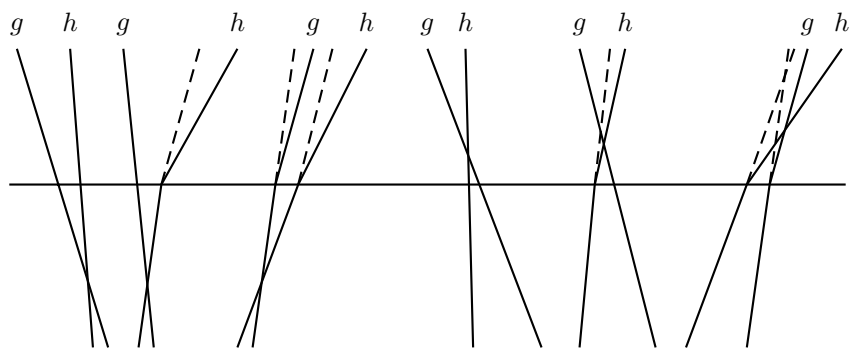

**Desargues' theorem** is false in this model geometry!

One sees this from the figure, in which $\Delta 123$ lies axially with $\Delta 1'2'3'$. (The axis of perspectivity is the line at infinity, since corresponding sides are taken to be parallel.) The triangles are *not* in perspective, however, since $h_1$ does not go through $S$.

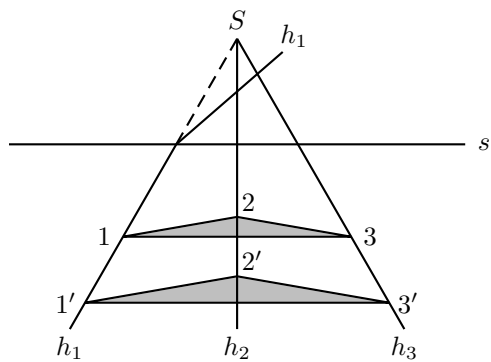

**Length measure** is the same as usual for ordinary lines, and for bent lines it is the sum of the ordinary lengths of the parts above and below the bend.

**Angle measurement** occurs only in the upper half plane. In order to determine the angle between $g$, $g'$, for example, we draw parallels $h$, $h'$ to $g$, $g'$ through an arbitrary fixed point $O$ in the upper half plane. Then let $\angle(g, g') \cong \angle(h, h')$.

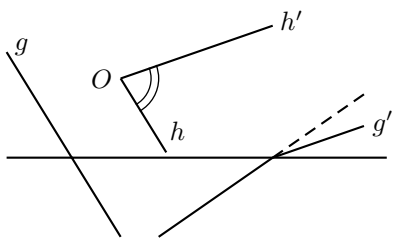

Verification of the remaining axioms is then evident, and we do not go further into them.

Axiom II.5 does not hold, as may be seen from the following figure, in which

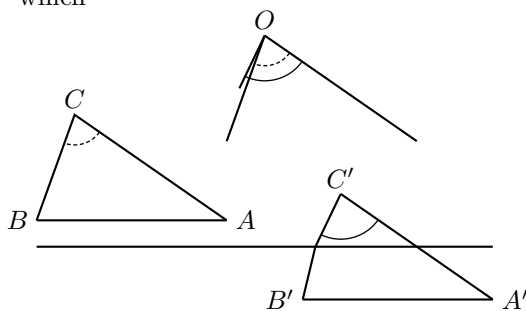

$$BA \cong B'A'$$
$$CA \cong C'A'$$
$$\angle CAB \cong \angle C'A'B'$$
but $\angle BCA \not\cong \angle B'C'A'$.

One sees in particular that Desargues cannot be proved from the planar axioms of incidence and order. Herein lies a fundamental difference between planar and spatial geometry: it is well known that Desargues' theorem in space is an immediate consequence of the spatial incidence axioms.

## 2.6  The theorem of Pappus

1. **The theorem of Pappus**[7] says:

   *If the vertices $1, 2, 3, 4, 5, 6$ of a hexagon lie alternately on lines $g$ and $h$, then the intersections of opposite sides, namely the points $P, Q, R$, which may be denoted*
   $$12 \cap 45, \quad 23 \cap 56, \quad 34 \cap 61,$$
   *lie on a line.*

   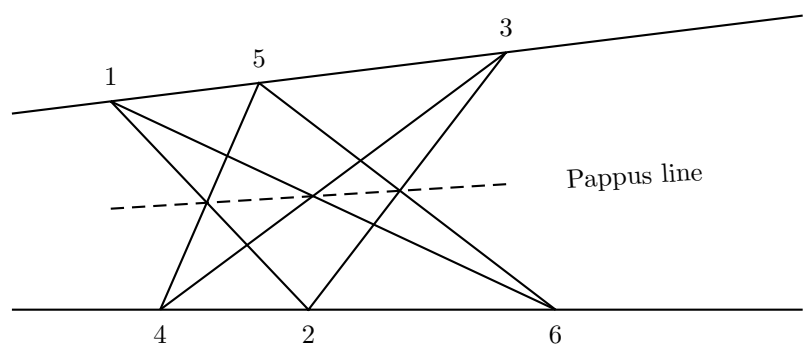

   []8

   The theorem is a special case of Pascal's theorem for conic sections, in which the conic degenerates to a pair of lines. In what follows we are concerned only with this special case. When the Pappus line is the line at infinity, one speaks of the affine Pappus theorem.

2. **Pappus' theorem does not follow from Desargues' theorem.**

   The proof is by construction of a model geometry, in which the incidence theorems, the parallel postulate, and Desargues' theorem hold. The geometry is the planar cartesian coordinate geometry over the Hilbert number system with parameters $s$ and $t$, described on p. 38. Thus
   $$ts = 2st \neq st.$$

---

[7]Translator's note. Like Hilbert, Moufang calls this the theorem of Pascal, but she points out that the theorem is only a special case of the usual Pascal's theorem. We revert to the historically accurate name for the special case, the theorem of Pappus, and use the term "Pascal's theorem" only for the general theorem about a hexagon with vertices on a conic section.

[8]Translator's note. Here we omit Moufang's sentence "This notation for the intersection of two lines will often be convenient, and we shall also use it without the bar over the connected points, since the distinction between segments and lines is no longer relevant." since we have already dispensed with the bar notation for lines. See page 8.

We make the following construction.

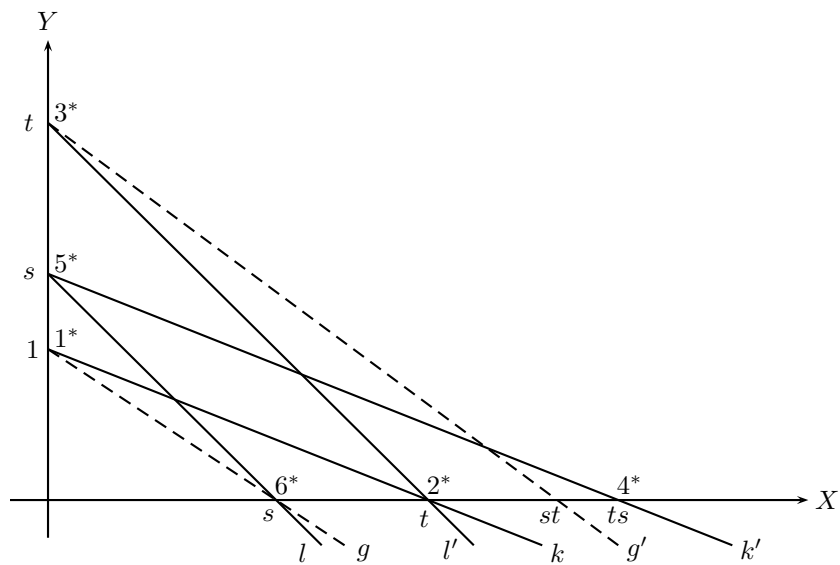

On the two coordinate axes we consider the ordinates $s$ and $t$, and on the $y$-axis also the ordinate 1. We now draw the lines

$$g \equiv sy + x - s = 0$$
$$g' \equiv sy + x - st = 0$$

which are parallel, with $g'$ meeting the $x$-axis at the point $st$;

$$k \equiv ty + x - t = 0$$
$$k' \equiv ty + x - ts = 0$$

which are parallel, with $k'$ meeting the $x$-axis at the point $ts$;

$$l \equiv x + y - s = 0$$
$$l' \equiv x + y - t = 0$$

which are parallel, with $l'$ meeting the $x$-axis at the point $t$.

According to Pappus' theorem, the connection of points $3^*$ and $4^*$ must be parallel to the connection of points $1^*$ and $6^*$, but this is false because $ts \neq st$. Thus Pappus' theorem is false in our geometry.

3. **Desargues' theorem follows from Pappus' theorem.**

This was shown by Hessenberg (*Math. Ann.*, vol. 61) by means of the incidence theorems and the parallel postulate. We make the following construction.

Let triangles $BCA$ and $B'C'A'$ be in perspective with center $O$. []⁹

---
⁹Translator's note. Here we omit Moufang's sentence "From now on we write $BCA$ persp. $B'C'A'$" because in fact she never uses this notation.

If $BA \parallel B'A'$ and $AC \parallel A'C'$ then the theorem says that $BC \parallel B'C'$.

We draw the construction line $m$ through $A$ parallel to $OB'$ and let

$$m \cap A'C' = M.$$

Also let $n = MB'$, $n \cap BA = N$, and $L = m \cap CC'$.

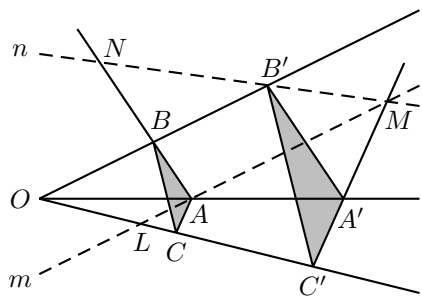

We now consider the following hexagons in turn:

α) In hexagon[10] $NAMA'B'O$, $NA \parallel B'A'$ and $OB' \parallel AM$, hence also $ON \parallel A'M = A'C'$ by Pappus' theorem.

β) In hexagon $NOBCAL$, $OB \parallel LA$ and $ON \parallel A'C' \parallel AC$, hence also $NL \parallel BC$ by Pappus' theorem.

γ) In hexagon $NOB'C'ML$, [$OB' \parallel LM$ and $ON \parallel AC \parallel MC'$ hence also] $NL \parallel B'C'$ [by Pappus' theorem.]

It follows from β) and γ) that

$$BC \parallel B'C'$$

and hence the theorem of Desargues is proved. Q.E.D.

4. **Pappus' theorem suffices as a basis for projective geometry.**

Suppose we are given two lines $g$ and $g'$ and a harmonic quadruple of points on each. A projective correspondence between $g$ and $g'$ is a one-to-one correspondence between their points. The fundamental theorem of projective geometry says that it is uniquely determined by three pairs of corresponding points, $AA'$, $BB'$, $CC'$, and that the correspondence may be set up by the operations of connection and intersection alone.

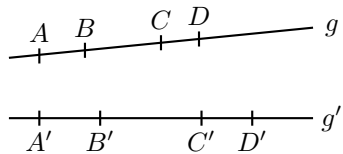

---

[10]Translator's note. Moufang writes the hexagons in α), β), γ) as $NB'MOAA'$, $NBAOLC$, $NB'MOLC'$ respectively, for some reason preferring to list the vertices out of order. I have rewritten their vertices in the correct cyclic order for the use of Pappus' theorem.

The provability of the fundamental theorem from Pappus' theorem will be shown in detail in Section 3.4.

Here we show the converse, namely, that *Pappus' theorem follows from the fundamental theorem.*

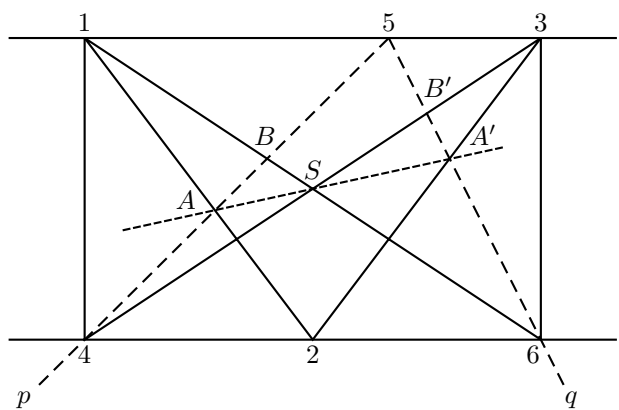

One maps the line $p$ by a perspectivity onto the pencil 1 and, likewise, $q$ onto the pencil 3. The product of two perspectivities is a projectivity. On the other hand, the pencil 1 is perspectively related to the point series 4, 2, 6 and also to the pencil 3, hence the pencils 1 and 3 are projectively related to each other and therefore $p$ is projectively related to $q$. But the latter projectivity is a perspectivity, since the intersection of $p$ and $q$ is fixed by the map, so the connections of corresponding points, namely $4B'$, $B6$, $AA'$ go through a single point $S$.

This means that
$$S = 16 \cap 34$$

lies on $AA'$. But this is the theorem of Pappus.                                Q.E.D.

## 2.7 Alternative fields

In contrast to ordinary fields, which are commutative and associative, and to skew fields, which are associative and noncommutative, alternative fields are noncommutative and nonassociative. That is, we have in general

$$a(bc) \neq (ab)c$$
$$ab \neq ba.$$

We require however the weak associative rules

$$a(ab) = (aa)b$$
$$(ba)a = b(aa)$$
$$a(ba) = (ab)a.$$

1. Alternative fields are realized by the hypercomplex algebra of Cayley numbers,[11] which today is the only known alternative field.[12] An alternative field satisfies the following axioms:

    1) **Addition:**
       $a + b = c$ is always defined and unique.
       Associative law: $(a + b) + c = a + (b + c)$.
       Commutative law: $a + b = b + a$.
       There is a zero element 0 such that $a + 0 = a = 0 + a$.

    2) **Multiplication:**
       $ab = c$ is always uniquely defined.

    3) **Distributive laws:**
       $$a(b + c) = ab + ac$$
       $$(b + c)a = ba + bc.$$

    4) **Weak associative laws:**
       $$a(ab) = (aa)b$$
       $$(ba)a = b(aa)$$
       $$a(ba) = (ab)a$$

       There is a unit element 1 such that $a \cdot 1 = 1 \cdot a = a$.

    5) Each element $a \neq 0$ has a *reciprocal* $a^{-1}$ such that
       $$a \cdot a^{-1} = 1 = a^{-1} \cdot a.$$

    4') There is also the following *special associative law*, but it can be proved from Axioms 1) to 4). (See R. Moufang: Alternativkörper und der Satz vom vollständigen Vierseit, *Hamburg Abh.* 1933.)
       $$a(a^{-1}b) = (aa^{-1})b = b, \quad (ba)a^{-1} = b(aa^{-1}) = b.$$

       Max Zorn (*Hamburg Abh.* vol. 8, 1931, p. 123 ff.) showed that these relations follow from the remaining axioms.[13]

       Cayley realized this alternative field as a hypercomplex algebra. The Cayley numbers use eight units over the real numbers. Their multiplication is given by the following table.

---

[11]Translator's note. I have retained the term *Cayley numbers* for what are now more often called the *octonions*, since it should not lead to confusion. However, Cayley gets more than his share of credit for the octonions, which were first discovered by John Graves in December 1843 and rediscovered by Cayley in 1845. Cayley was first to mention them in print, mainly because Hamilton promised to communicate Graves' discovery but didn't get around to it.

[12]Translator's note. More precisely, they *are* the only alternative field over the real numbers. There are several theorems about the uniqueness of the Cayley numbers, for which a good reference is *Numbers* by Ebbinghaus et al., Springer-Verlag 1991.

[13]Translator's note. In fact 4') is *equivalent* to 4) in the presence of the other axioms. Moufang gives the converse of Zorn's result in her paper cited above. She attributes the proof of $a(ba) = (ab)a$ to Reidemeister.

|  | $e_0 = 1$ | $e_1$ | $e_2$ | $e_3$ | $e_4$ | $e_5$ | $e_6$ | $e_7$ |
|---|---|---|---|---|---|---|---|---|
| $e_0 = 1$ | 1 | $e_1$ | $e_2$ | $e_3$ | $e_4$ | $e_5$ | $e_6$ | $e_7$ |
| $e_1$ | $e_1$ | $-1$ | $e_3$ | $-e_2$ | $e_5$ | $-e_4$ | $-e_7$ | $e_6$ |
| $e_2$ | $e_2$ | $-e_3$ | $-1$ | $e_1$ | $e_6$ | $e_7$ | $-e_4$ | $-e_5$ |
| $e_3$ | $e_3$ | $e_2$ | $-e_1$ | $-1$ | $e_7$ | $-e_6$ | $e_5$ | $-e_4$ |
| $e_4$ | $e_4$ | $-e_5$ | $-e_6$ | $-e_7$ | $-1$ | $e_1$ | $e_2$ | $e_3$ |
| $e_5$ | $e_5$ | $e_4$ | $-e_7$ | $e_6$ | $-e_1$ | $-1$ | $-e_3$ | $e_2$ |
| $e_6$ | $e_6$ | $e_7$ | $e_4$ | $-e_5$ | $-e_2$ | $e_3$ | $-1$ | $-e_1$ |
| $e_7$ | $e_7$ | $-e_6$ | $e_5$ | $e_4$ | $-e_3$ | $-e_2$ | $e_1$ | $-1$ |

2. There is a *simpler presentation* of this system due to *Dickson*,[14] which takes the quaternions as ground field. The general number of the alternative field is defined to be
$$a = q + Qe,$$
where $q$ and $Q$ are quaternions over the real numbers and $e$ is an arbitrary symbol. We now define:

**Addition:** $a + b = (q + Qe) + (r + Re) = (q + r) + (Q + R)e$.

**Multiplication:** $a \cdot b = (qr - \overline{R}Q) + (Rq + Q\overline{r})e$,

where $\overline{q} = \alpha_0 e_0 - \alpha_1 e_1 - \alpha_2 e_2 - \alpha_3 e_3$

is the conjugate of $q = \alpha_0 e_0 + \alpha_1 e_1 + \alpha_2 e_2 + \alpha_3 e_3$.

When one sets $a = \alpha_0 e_0 + \alpha_1 e_1 + \alpha_2 e_2 + \cdots + \alpha_7 e_7$,

then the eight units in this presentation are

$$e_0 = 1, \quad e_1, \quad e_2, \quad e_3,$$
$$e_0 e = e_4, \quad e_1 e = e_5, \quad e_2 e = e_6, \quad e_3 e = e_7.$$

3. The **computation rules** of the alternative field may now be derived from the Dickson presentation.

   **Addition** is commutative and associative, since this is true for the quaternions.

   Zero element: $a = 0$ when $q = Q = 0$.

   $$a + \tilde{a} = 0 \quad \longrightarrow \quad a = -q - Qe, \quad \text{with} \quad \tilde{a} = -a.$$

   **The distributive law**
   $$a(b_1 + b_2) = ab_1 + ab_2$$
   will now be proved. If
   $$b_1 + b_2 = b$$
   $$= (r_1 + R_1 e) + (r_2 + R_2 e)$$
   $$= (R_1 + r_2) + (R_1 + R_2)e$$

---

[14]Translator's note. *Linear Algebras*, Cambridge University Press 1914, p. 15. It is worth mentioning that Dickson's construction also produces the complex numbers from the real numbers, and the quaternions from the complex numbers.

then
$$a(b_1 + b_2) = q(r_1 + r_2) - (\overline{R}_1 + \overline{R}_2)Q + \{(R_1 + R_2)q + Q(\overline{r}_1 + \overline{r}_2)\}e$$
$$= qr_1 - \overline{R}_1 Q + (R_1 q + Q\overline{r}_1)e + qr_2 - \overline{R}_2 Q + (R_2 q + Q\overline{r}_2)e$$
$$= ab_1 + ab_2 \qquad \text{Q.E.D.}$$

The other distributive law is proved similarly.

**Multiplication.**

First we exhibit the existence of inverses. The identity element is defined by $Q = 0$, $q = 1$.

Here we need to insert a remark on the *norm theorem*. The norm of a Cayley number, like that of a complex number or quaternion, is the product of a Cayley number $a = q + Qe$ by its conjugate $\overline{a} = \overline{q} - Qe$.

Thus it follows from the definition of multiplication that

$$a\overline{a} = N(a) = q\overline{q} + Q\overline{Q}.$$

The norm theorem says that *the norm of a product equals the product of the norms of the factors*.

For *complex numbers* this is easily verified by computation:

$$a = \alpha_0 + \alpha_1 e_1 \qquad N(a) = \alpha_0^2 + \alpha_1^2$$
$$b = \beta_0 + \beta_1 e_1 \qquad N(b) = \beta_0^2 + \beta_1^2$$
$$ab = (\alpha_0 + \alpha_1 e_1)(\beta_0 + \beta_1 e_1)$$
$$= \alpha_0 \beta_0 - \alpha_1 \beta_1 + (\alpha_0 \beta_1 + \alpha_1 \beta_0)e_1$$
$$N(a) \cdot N(b) = (\alpha_0^2 + \alpha_1^2)(\beta_0^2 + \beta_1^2)$$
$$= \alpha_0^2 \beta_0^2 + \alpha_0^2 \beta_1^2 + \alpha_1^2 \beta_0^2 + \alpha_1^2 \beta_1^2$$
$$N(ab) = (\alpha_0 \beta_0 - \alpha_1 \beta_1)^2 + (\alpha_0 \beta_1 + \alpha_1 \beta_0)^2$$
$$= \alpha_0^2 \beta_0^2 + \alpha_0^2 \beta_1^2 + \alpha_1^2 \beta_0^2 + \alpha_1^2 \beta_1^2,$$

hence $N(ab) = N(a) \cdot N(b)$. \hfill Q.E.D.

One also calls this the Euler-Lagrange identity.[15]

The norm theorem also holds for quaternions. In fact:

$$ab = (\alpha_0 + \alpha_1 e_1 + \alpha_2 e_2 + \alpha_3 e_3)(\beta_0 + \beta_1 e_1 + \beta_2 e_2 + \beta_3 e_3)$$
$$= (\alpha_0 \beta_0 - \alpha_1 \beta_1 - \alpha_2 \beta_2 - \alpha_3 \beta_3)$$
$$+ (\alpha_0 \beta_1 + \alpha_1 \beta_0 + \alpha_2 \beta_3 - \alpha_3 \beta_2)e_1$$
$$+ (\alpha_0 \beta_2 + \alpha_2 \beta_0 + \alpha_3 \beta_1 - \alpha_1 \beta_3)e_2$$
$$+ (\alpha_0 \beta_3 + \alpha_3 \beta_0 + \alpha_1 \beta_2 - \alpha_2 \beta_1)e_3,$$

---

[15] Translator's note. This term, if used at all today, is probably applied to the four square identity that follows. The four square identity was discovered by Euler in 1748, and used by Lagrange in 1770 to prove that every natural number is the sum of four squares.

The two square identity, corresponding to the norm theorem for complex numbers, appears to have been known to Diophantus. He mentions a special case of it in Book III, Problem 19 of his *Arithmetica*.

$$N(ab) = (\alpha_0^2 + \alpha_1^2 + \alpha_2^2 + \alpha_3^2)(\beta_0^2 + \beta_1^2 + \beta_2^2 + \beta_3^2)$$
$$= (\alpha_0\beta_0 - \alpha_1\beta_1 - \alpha_2\beta_2 - \alpha_3\beta_3)^2$$
$$+ (\alpha_0\beta_1 + \alpha_1\beta_0 + \alpha_2\beta_3 - \alpha_3\beta_2)^2$$
$$+ (\alpha_0\beta_2 + \alpha_2\beta_0 + \alpha_3\beta_1 - \alpha_1\beta_3)^2$$
$$+ (\alpha_0\beta_3 + \alpha_3\beta_0 + \alpha_1\beta_2 - \alpha_2\beta_1)^2,$$

hence $N(ab) = N(a) \cdot N(b)$. Q.E.D.

The norm theorem also holds for the Cayley numbers. In order to show this, we first prove the following:

**Lemma:** *For any quaternions $Q$ and $R$*

$$\overline{Q \cdot R} = \overline{R} \cdot \overline{Q},$$

*where the bar denotes the conjugate.*

Proof. We have

$$N(Q) = Q \cdot \overline{Q}$$
$$N(R) = R \cdot \overline{R}$$
$$QR \cdot \overline{QR} = N(QR) = N(Q) \cdot N(R) = Q\overline{Q} \cdot R\overline{R}.$$

Left multiplying this by $Q^{-1}$ gives

$$R \cdot \overline{QR} = \overline{Q} \cdot R\overline{R}$$
$$= R\overline{R} \cdot \overline{Q} \quad \text{since } R\overline{R} \text{ is real.}$$

Then left multiplication by $R^{-1}$ gives

$$\overline{QR} = \overline{R} \cdot \overline{Q} \qquad \text{Q.E.D.}$$

Using the Dickson presentation, it may now be shown that *each nonzero Cayley number has an inverse.*

$$a = q + Qe$$
$$\overline{a} = \overline{q} - Qe$$
$$a\overline{a} = q\overline{q} + \overline{Q}Q + (-Qq + Qq)e$$
$$= \sum_{i=0}^{7} \alpha_i^2 = N(q) + N(Q) = N(a) = \text{real},$$

and therefore

$$\frac{\overline{a}}{N(a)} = \frac{\overline{q}}{N(a)} - \frac{Q}{N(a)} = a^{-1}.$$

Now we come to the proof of the norm theorem for Cayley numbers!

We have $ab = (q + Qe)(r + Re) = qr - \overline{R}Q + (Rq + Q\overline{r})$.

Therefore, since $N(a) = \sum_{i=0}^{7} \alpha_i^2 = q\bar{q} + \overline{Q}Q$,

$$N(ab) = (qr - \overline{R}Q)\overline{(qr - \overline{R}Q)} + (Rq + Q\bar{r})\overline{(Rq + Q\bar{r})}$$
$$= (qr - \overline{R}Q)(\bar{r}\,\bar{q} - \overline{Q}R) + (Rq + Q\bar{r})(\bar{q}\,\overline{R} + r\overline{Q})$$
by the lemma
$$= qr\bar{r}\,\bar{q} - qr\overline{Q}R - \overline{R}Q\bar{r}\,\bar{q} + \overline{R}Q\overline{Q}R + Rq\bar{q}\,\overline{R} + Rqr\overline{Q} + Q\bar{r}\,\bar{q}\,\overline{R} + Q\bar{r}r\overline{Q},$$
whereas
$$N(a)N(b) = (q\bar{q} + \overline{Q}Q)(r\bar{r} + \overline{R}R)$$
$$= q\bar{q}r\bar{r} + q\bar{q}\,\overline{R}R + \overline{Q}Qr\bar{r} + \overline{Q}QR\overline{R}.$$

When one observes that

$$q\bar{q}r = rq\bar{q} \quad \text{(because } q\bar{q} \text{ is real)}$$

and

$$\bar{q}q = q\bar{q},$$

then one sees that

$$N(ab) - N(a)N(b) = -qr\overline{Q}r - \overline{R}Q\bar{r}\,\bar{q} + Rqr\overline{Q} + Q\bar{r}\,\bar{q}\,\overline{R}$$
$$= -(qr\overline{Q}r + \overline{R}Q\bar{r}\,\bar{q}) + (Rqr\overline{Q} + Q\bar{r}\,\bar{q}\,\overline{R})$$

Multiplying $(Rqr\overline{Q} + Q\bar{r}\,\bar{q}\,\overline{R})$ on the left by $\overline{R}$ and on the right by $R$, then dividing by $\overline{R}R$, gives[16]

$$N(ab) - N(a)N(b) = -(qr\overline{Q}r + \overline{R}Q\bar{r}\,\bar{q}) + (qr\overline{Q}r + \overline{R}Q\bar{r}\,\bar{q}) = 0 \quad \text{Q.E.D.}$$

The norm theorem for Cayley numbers implies the following algebraic identity:[17]

$$\sum_0^7 \alpha_i^2 \cdot \sum_0^7 \beta_i^2 = (\alpha_0\beta_0 - \sum_1^7 \alpha_i\beta_i)^2$$
$$+ (\alpha_0\beta_1 + \alpha_1\beta_0 + \alpha_2\beta_3 - \alpha_3\beta_2 + \alpha_4\beta_5 - \alpha_5\beta_4 - \alpha_6\beta_7 + \alpha_7\beta_6)^2$$
$$+ (\alpha_0\beta_2 + \alpha_2\beta_0 - \alpha_1\beta_3 + \alpha_3\beta_1 + \alpha_4\beta_6 - \alpha_6\beta_4 + \alpha_5\beta_7 - \alpha_7\beta_5)^2$$
$$+ (\alpha_0\beta_3 + \alpha_3\beta_0 + \alpha_1\beta_2 - \alpha_2\beta_1 + \alpha_4\beta_7 - \alpha_7\beta_4 - \alpha_5\beta_6 + \alpha_6\beta_5)^2$$
$$+ (\alpha_0\beta_4 + \alpha_4\beta_0 - \alpha_1\beta_5 + \alpha_5\beta_1 - \alpha_2\beta_6 + \alpha_6\beta_2 - \alpha_3\beta_7 + \alpha_7\beta_3)^2$$
$$+ (\alpha_0\beta_5 + \alpha_5\beta_0 + \alpha_1\beta_4 - \alpha_4\beta_1 - \alpha_2\beta_7 + \alpha_7\beta_2 + \alpha_3\beta_6 - \alpha_6\beta_3)^2$$
$$+ (\alpha_0\beta_6 + \alpha_6\beta_0 + \alpha_1\beta_7 - \alpha_7\beta_1 + \alpha_2\beta_4 - \alpha_4\beta_2 - \alpha_3\beta_5 + \alpha_5\beta_3)^2$$
$$+ (\alpha_0\beta_7 + \alpha_7\beta_0 - \alpha_1\beta_6 + \alpha_6\beta_1 + \alpha_2\beta_5 - \alpha_5\beta_2 + \alpha_3\beta_4 - \alpha_4\beta_3)^2.$$

---

[16]Translator's note. These operations produce no change, because $(Rqr\overline{Q} + Q\bar{r}\,\bar{q}\,\overline{R})$ is an element plus its conjugate, hence real.

[17]Translator's note. The eight square identity was discovered by Degen, *Mém. l'Acad. Imp. Sci. St. Petersbourg*, VIII (1822), pp. 207–219, and rediscovered by Graves in 1843. This led to Graves' discovery of the octonions. Thus each hypercomplex number system— complex numbers, quaternions, and octonions (or Cayley numbers)—was foreshadowed by the corresponding norm theorem, in the form of an identity involving sums of squares.

Thus the product of two sums of eight squares is again a sum of eight squares, generalizing the Euler-Lagrange identity. Such identities for sums of $n$ squares exist only for $n = 1, 2, 4, 8$.[18]

4. Of the *weak associativity laws* for an alternative field we prove only two here; the proof of the others is analogous.

   1) Claim: $a^{-1}(ab) = (a^{-1}a)b = b$.
      Proof. Given
      $$a = q + Qe, \quad b = r + Re,$$
      we have
      $$ab = qr - \overline{R}Q + (Rq + Q\overline{r})e = r_1 + R_1 e.$$
      Hence
      $$a^{-1}(ab)$$
      $$= \frac{\overline{q}}{q\overline{q} + Q\overline{Q}} r_1 + \overline{R}_1 \frac{Q}{q\overline{q} + Q\overline{Q}} + \left( R_1 \frac{\overline{q}}{q\overline{q} + Q\overline{Q}} - \frac{Q}{q\overline{q} + Q\overline{Q}} \overline{r}_1 \right) e$$
      $$= \frac{\overline{q}}{N(a)}(qr - \overline{r}Q) + \overline{(Rq + Q\overline{r})} \frac{q}{N(a)}$$
      $$\qquad + \left\{ (Rq + Q\overline{r}) \frac{\overline{q}}{N(a)} - \frac{Q}{N(a)} \overline{(qr - \overline{R}Q)} \right\} e$$
      $$= \frac{N(q)r}{N(a)} - \frac{\overline{q}\overline{R}Q}{N(a)} + \frac{\overline{q}\overline{R}Q}{N(a)} + \frac{r\overline{Q}Q}{N(a)}$$
      $$\qquad + \left\{ \frac{RN(q)}{N(a)} + \frac{Q\overline{r}\,\overline{q}}{N(a)} - \frac{Q}{N(a)}\overline{r}\,\overline{q} + \frac{Q\overline{Q}R}{N(a)} \right\} e$$
      $$= \frac{N(q)r + rN(Q)}{N(q) + N(Q)} + \left\{ \frac{RN(q) + N(Q)R}{N(q) + N(Q)} \right\} e$$
      $$= r + Re = b \qquad\qquad \text{Q.E.D.}$$

   2) Claim: $a(ab) = (aa)b$
      Proof. We have
      $$a(ab)$$
      $$= (q + Qe)(r_1 + R_1 e) = qr_1 - \overline{R}_1 Q + (R_1 q + Q\overline{r}_1)e$$
      $$= q(qr - \overline{R}Q) - (\overline{q}\overline{R} + r\overline{Q})Q + \left\{ (Rq + Q\overline{r})q + Q(\overline{r}\,\overline{q} - \overline{Q}R) \right\} e$$
      $$= \underline{q(qr)} - \underline{\underline{q\overline{R}Q - \overline{q}\overline{R}Q}} - rN(Q) + \left\{ Rqq + \underline{\underline{Q\overline{r}q + Q\overline{r}\,\overline{q}}} - \underline{N(Q)R} \right\} e,$$

      Also, $aa = (q + Qe)(q + Qe) = qq - \overline{Q}Q + (Qq + Q\overline{q})e = q_1 + Q_1 e,$

---
[18]Translator' note. First proved by Hurwitz, *Göttinger Nachrichten*, 1898, pp. 309–316.

hence

$(aa)b$

$$= (q_1 + Q_1 e)(r + Re) = q_1 r - \overline{R} Q_1 + (R q_1 - Q_1 \overline{r})e$$
$$= (qq - N(Q))r - \overline{R}(Qq + Q\overline{q}) + \{R(qq - N(Q)) + (Qq + Q\overline{q})\overline{r}\} e$$
$$= \underline{(qq)r} - \underline{N(Q)r} - \underline{\overline{R}Qq} - \underline{\underline{\overline{R}Q\overline{q}}} + \left\{\underline{R(qq)} - \underline{RN(Q)} + \underline{\underline{Qq\overline{r} + Q\overline{q}\,\overline{r}}}\right\} e$$

When one recalls that quaternions are associative, and that they commute with real numbers, the underlined terms in this expression for $(aa)b$ can be seen to equal the underlined terms in the expression for $a(ab)$. The doubly underlined terms are equal because they each have a factor $q + \overline{q}$ which, being real, can also be commuted. Q.E.D.

To prove that the Cayley numbers are *nonassociative* in general, it suffices to find a triple $a$, $b$, $c$ for which $(ab)c \neq a(bc)$.

Claim:

*The triple $Q$, $e$, $q$ satisfies $Q(eq) \neq (QE)q$.*

Proof. We have

$$eq = (0 + 1 \cdot e)(q + 0 \cdot e) = 0 \cdot q - 0 \cdot 1 + (0 \cdot 0 + 1 \cdot \overline{q})e = \overline{q}e$$
$$Q(eq) = Q(\overline{q}e) = (Q + 0e)(0 + \overline{q}e) = Q \cdot 0 - q \cdot 0 + (\overline{q}Q + 0)e = \overline{q}Qe$$
$$(Qe)q = (0 + Qe)(q + 0e) = 0 \cdot q - 0 \cdot Q + (0 + Q\overline{q})e = Q\overline{q}e.$$

Since $\overline{q}Q \neq Q\overline{q}$ in general, for example, for $Q = e_1$, $\overline{q} = e_2$, the claim is proved. Thus alternative fields are nonassociative in general.      Q.E.D.

5. **General structural properties of alternative fields.**

   $\alpha$) **The associative laws** may be written in the following form.
   We denote $(xy)z - x(yz)$ by $[x, y, z]$, so $[x, y, z] \neq 0$ in general. However, we have the weak associative laws:
   
   $$[x, x, y] = 0, \quad [y, x, x] = 0, \quad [x, y, x] = 0.$$
   
   With this notation it is easy to prove:
   
   **Theorem:**
   
   *The symbol $[x, y, z]$ changes sign ("alternates") with the exchange of two elements.*
   
   Proof. We have
   
   1) $[a + b, a + b, c] = 0 = [a, a, c] + [b, a, c] + [b, b, c] + [a, b, c]$,
      hence $[b, a, c] = -[a, b, c]$.
   2) $[a, b + c, b + c] = 0 = [a, b, b] + [a, b, c] + [a, c, c] + [a, c, b]$
      hence $[a, c, b] = -[a, b, c]$.                                 Q.E.D.
   3) In particular, if we set $c = b$ then we get $[b, a, b] = -[a, b, b]$.
      This proves the third weak associative rule,
      
      $$a(ba) = (ab)a,$$
      
      mentioned at the beginning of this section.

From 1) and 2) one sees that, when $[x, y, z] = 0$, then the symbol is in fact zero for all six permutations of $x$, $y$, $z$.

$\beta$) **The rule** $(ab)^{-1} = b^{-1}a^{-1}$ holds in an alternative field.

By applying the identities

$$(u^{-1}u)v = v \quad \text{and} \quad (vu)u^{-1}$$

we get

$$b^{-1} = \left\{(ab)^{-1}(ab)\right\} b^{-1} = (ab)^{-1}a,$$

hence

$$b^{-1}a^{-1} = ((ab)^{-1}a)a^{-1} = (ab)^{-1} \qquad \text{Q.E.D.}$$

$\gamma$) **The four identity** holds in any distributive ring:

$$[ab, c, d] - [a, bc, d] + [a, b, cd] = a[b, c, d] + [a, b, c]d.$$

This follows from the distributive laws, as one sees by expanding the brackets.

$\delta$) **A commutative alternative field** is an ordinary field, and hence the general associative law holds. Namely

$$\begin{aligned}
0 = \ & (ab - ba)c - c(ab - ba) \\
& + (bc - cb)a - a(bc - cb) \\
& + (ca - ac)b - b(ca - ac) \\
= \ & [a, b, c] + [b, c, a] + [c, a, b] - [b, a, c] - [c, b, a] - [a, c, b] \\
= \ & 6[a, b, c]
\end{aligned}$$

Hence $[a, b, c] = 0$ if $ab = ba$, $ca = ac$, $bc = cb$. $\qquad$ Q.E.D.

$\varepsilon$) **Ordered alternative fields**

The axioms of linear order are evidently not satisfied by the Cayley numbers. Since the latter are the only known realization of an alternative field, it remains an open problem whether there is a linearly ordered alternative field.[19] However, we have the theorem:

*An archimedean ordered alternative field is an ordinary field.*

Proof. An alternative field contains an identity element. We define

$$\begin{aligned}
2 = 1 + 1, \quad & \ldots, \quad (m + 1) = m + 1, \ldots \\
-1 = 0 - 1, \quad & \ldots, \quad -(m + 1) = -m - 1, \ldots,
\end{aligned}$$

which are all different because the alternative field is linearly ordered. From these "integers" we build the "rational numbers" of the alternative field:

$$\alpha = mn^{-1}.$$

Then if $a$ is an arbitrary member of the alternative field we have

$$\alpha a = a\alpha.$$

---

[19] Translator's note. This question was answered in the negative by the theorem of Bruck and Kleinfeld, *Proc. Amer. Math. Soc.* 2, (1951), according to which every alternative field is an octonion algebra.

First we have $ma = am$ because

$$ma = (1 + 1 + \cdots + 1)a = a + a + \cdots + a = a(1 + 1 + \cdots + 1).$$

Likewise, $(ma)b = m(ab)$ because

$$(a + a + \cdots + a)b = ab + ab + \cdots + ab = (1 + 1 + \cdots + 1)ab = m(ab).$$

By $\alpha$), $[m, a, b] = 0$ implies

$$[a, m, b] = 0 = [a, b, m].$$

Therefore

$$\begin{aligned}
a(mn^{-1}) &= (am)n^{-1} = (ma)n^{-1} = m(an^{-1}) = \left((an^{-1})^{-1}m^{-1}\right)^{-1} \\
&= \left((na^{-1})m^{-1}\right)^{-1} = \left((a^{-1}n)m^{-1}\right)^{-1} = \left(a^{-1}(nm^{-1})\right)^{-1} \\
&= (nm^{-1})^{-1}a = (mn^{-1})a.
\end{aligned}$$

When $\alpha$ is a "rational number" of the alternative field one also has

$$(\alpha a)b = \alpha(ab), \quad \text{that is, } [\alpha, a, b] = 0.$$

Indeed, it follows from

$$(ma)b = m(ab)$$

first of all that

$$\{(ma)b\}^{-1} = \{m(ab)\}^{-1},$$

and this may be written

$$b^{-1}(a^{-1}b^{-1}) = (b^{-1}a^{-1})m^{-1}.$$

Then if $b^{-1} = v$ and $a^{-1} = u$ we have

$$v(um^{-1}) = (vu)m^{-1} \quad \text{that is, } [v, u, m^{-1}] = 0,$$

where $u$ and $v$ are arbitrary and $m$ is an integer.

The four identity now gives

$$[mn^{-1}, a, b] - [m, n^{-1}a, b] + [m, n^{-1}, ab] = m[n^{-1}, a, b] + [m, n^{-1}, a]b,$$

hence $[mn^{-1}, a, b] = [\alpha, a, b] = 0.$ \hfill Q.E.D.

This shows that the subdomain of numbers $\alpha = mn^{-1}$ satisfies the computation rules of the rational numbers, and hence is isomorphic to them.

Now suppose that the alternative field is Archimedeanly ordered. Let $a$ and $b$ be positive with $a < b$. Then there is a rational $\alpha$ with $a < \alpha < b$. It follows from this that the commutative law of multiplication holds.

If not, there would be at least one pair $a$, $b$ with $ab < \alpha < ba$. Multiplying this inequality on the right by $a$ gives

$$(ab)a < \alpha a,$$

and multiplication on the left similarly gives

$$a\alpha < a(ba).$$

Since $a\alpha = \alpha a$ and $a(ba) = a(ba)$, this is a contradiction.

From the commutative law, it follows by $\delta$) that the associative law is also valid. Thus it is proved that an Archimedeanly ordered alternative field is an ordinary field. The question of a geometric equivalent to an alternative field is treated in Section 3.3.

# Chapter 3

# Foundations of projective geometry in the plane

## 3.1 Special cases of Desargues' theorem

1. In future we denote the general form of Desargues' theorem on perspective triangles in the plane by $D_0$, and special cases of it by $D_1$, $D_2$, etc.[1]

   Also, Desargues' theorem will always be understood in its double form: that perspective triangles lie axially, and conversely. "Perspective" means the lines through corresponding vertices go through a point. "Axial" means that the intersections of corresponding sides lie on a line.[2]

   $D_1$: A triangle $\Delta 123$ has one vertex (3) on a side of the other triangle $\Delta 1'2'3'$, or else an intersection (3) of corresponding sides lies on the line through corresponding vertices (of triangles $\Delta 11'B$ and $\Delta 22'C$). We denote these two types of special Desargues' theorem by $D_1$.

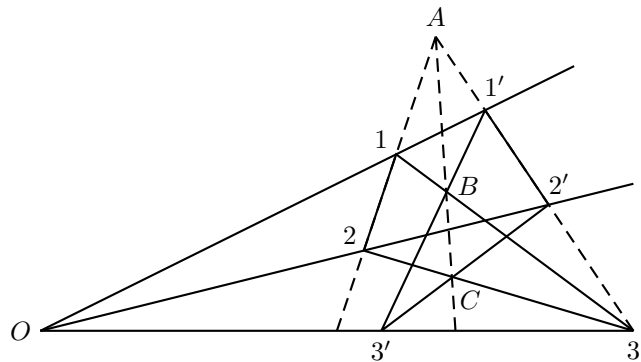

---

[1]Translator's note. In her 1930s papers, Moufang denoted these theorems by $D_{11}$, $D_{10}$, $D_9$, etc. At that time, the subscript denoted the number of free parameters in the configuration; here it evidently denotes the number of constraints.

[2]Translator's note. More explicit terms are normally used in English: "in perspective from a point" for "perspective", and "in perspective from a line" for "axial". However, Moufang's terminology seems clear enough, so I have not altered it.

Our figure contains yet another special Desargues' theorem, which Hessenberg[3] gave the name "little Desargues' theorem": it concerns triangles $\Delta 12O$ and $\Delta BC3'$ whose center of perspective (3) lies on the axis $A2$.

We call this the little Desarguish theorem to distinguish it from the other special cases $D_1$.

$D_2$: Each of the triangles $\Delta 123$ and $\Delta 1'2'3'$ has a vertex on a side of the other, or two intersections of corresponding sides lie on two lines connecting corresponding vertices. ($\Delta 11'B$ and $\Delta 22'C$), or one triangle has two vertices on sides of the other ($\Delta 11'A$ and $\Delta 33'C$).

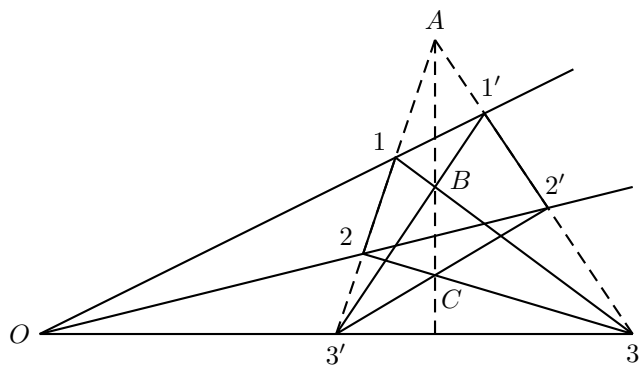

We denote these three special cases of Desargues' theorem by $D_2$.

It is worth mentioning here that $D_2$ (that is, any of its three cases) represents the geometric equivalent of an alternative field, as we shall show in the next section.

$D_3$: All three vertices 1, 2, 3 lie on sides of the triangle $\Delta 1'2'3'$.

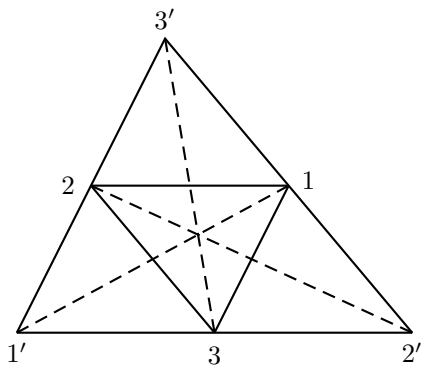

We content ourselves here with this rather unsystematic account of special Desargues' theorems. We shall not be concerned with theorem $D_3$, but $D_1$ and $D_2$ are fundamental for the investigations below.

---

[3]Translator's note. Hessenberg introduced this term in his paper Begründung der elliptischen Geometrie in *Math. Annalen* 61 (1906), p. 178.

2α) *Theorem $D_2$ is equivalent to the assertion that the position of the fourth harmonic point $D$ for three given points $A, B, C$ is independent of the quadrilateral used to construct it.*

First, $s_1$ through $C$ is arbitrary, because if $s_2 \neq s_1$ goes though $C$ then the quadrilateral $ABQQ'$ always leads to the point $D$, since $D_2$ in $\Delta PQB$ and $\Delta P'Q'A$ implies that the points $O$, $M$, $N$ are collinear. Conversely,[4] if the position of $D$ is independent of the construction elements $A, B, O$ then $O, M, N$ are collinear, that is, theorem $D_2$ holds.

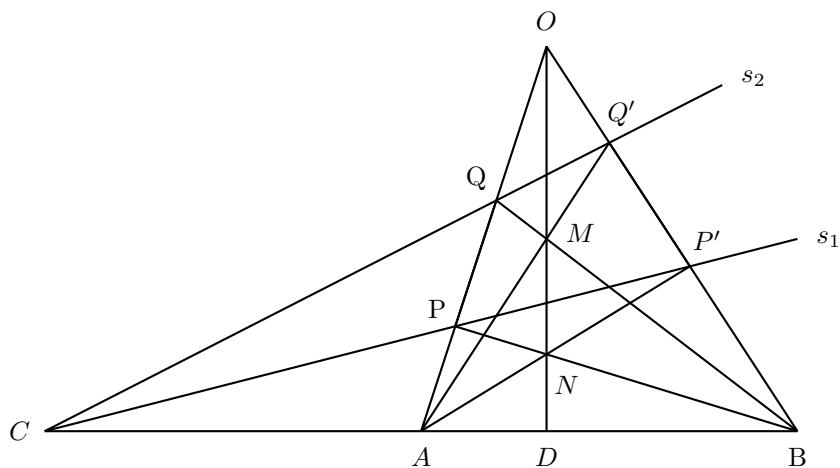

It remains to show, with the help of theorem $D_2$, that moving $O$ to $O_1$ results in the same point $D$. Here $s$ remains fixed.

We repeat the construction with the quadrilateral $ABP_1Q_1$ and claim that $O_1, M_1, D$ *are collinear.*

---

[4]Translator's note. This "converse" is in fact incorrect. It repeats a mistake made on p. 762 of Moufang's paper in *Math. Annalen* 106, (1932)—assuming that the diagonal points of a complete quadrilateral cannot lie in a line. They can (e.g. in the 7-point plane), and only a partial converse holds: *if no complete quadrilateral in the projective plane has collinear diagonal points, then the theorem of the complete quadrialateral (CQ) inplies the little Desargues theorem* (see Pickert *Projective Ebenen*, Springer-Verlag, 1955). The mistake was apparently first pointed out by Marshall Hall, on p. 269 of his paper on projective planes in *Trans. Amer. Math. Soc.* 54 (1943).

Because of the mistake, Moufang's research took a wrong turn in 1932–3. Believing CQ to be equivalent to $D_2$, she rewrote the main result of her 1930 dissertation—that $D_2$ is equivalent to coordinatization by an alternative field—as a "proof" that CQ is equivalent to coordinatization by an alternative field. However, in her "proof" that CQ ⇒ alternative field (p. 766 in *Math. Annalen* 106) she actually invokes little Desargues at each point, thinking it is a consequence of CQ. Thus in fact she has proved that little Desargues ⇒ alternative field.

In her followup paper, *Abh. Math. Sem. Hamburg* 9 (1933), Moufang continues to assume that CQ ⇔ $D_2$, but in fact she really proves that alternative field ⇒ $D_2$. This proof is easily modified to show that alternative field ⇒ little Desargues, and no doubt Moufang would have done so had she noticed her mistake about CQ. This is why Moufang gets credit for the theorem: little Desargues ⇔ alternative field.

In the present manuscript, Moufang returns to stating the theorem in its original form, $D_2$ ⇔ alternative field, but *it is still not known whether this theorem is correct*. The first to find a valid replacement for $D_2$—namely little Desargues—seems to have been Hall in the 1943 paper cited above. Moufang gives no sign of having seen Hall's paper.

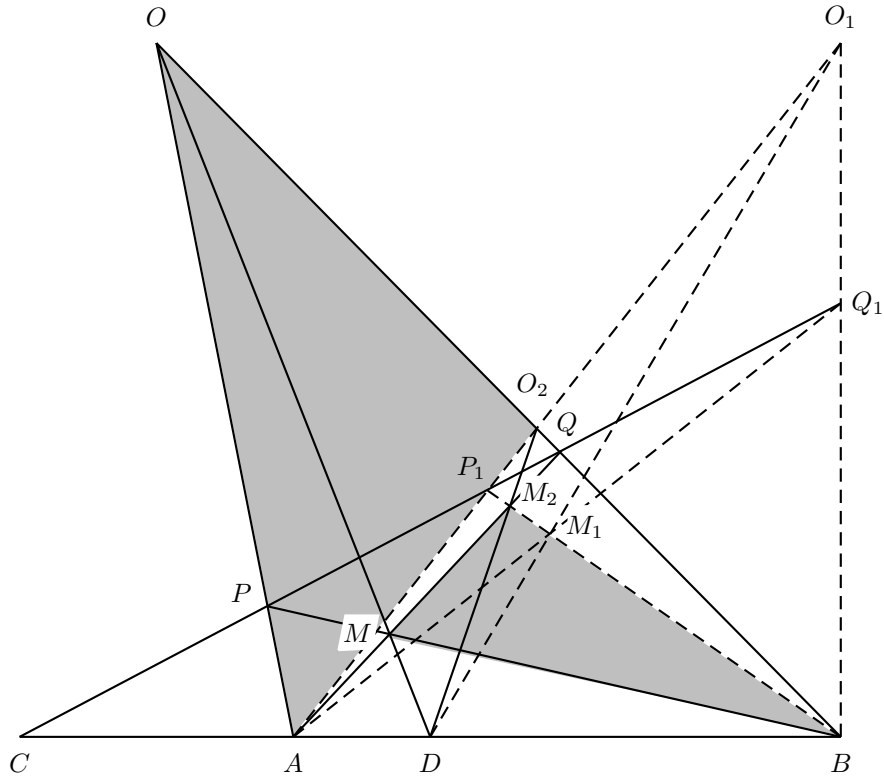

Proof. Consider $O_2$ on $OB$ and the [shaded] triangles $\Delta BMM_2$ and $\Delta AOO_2$. They are axial, since

$$BM \cap AO = P, \quad BM_2 \cap AO_2 = P_1, \quad MM_2 \cap OO_2 = Q$$

are collinear.

This is a $D_2$ configuration. Consequently, the triangles in question are in perspective [from $D$] and $O_2, M_2, D$ are collinear. Now we formally replace $O$ by $O_2$ and $O_2$ by $O_1$.

It follows by the same argument that $O_1, M_1, D$ are collinear.    Q.E.D.

Assuming theorem $D_2$, this proves that the position of the fourth harmonic point is independent of the choice of $O$ and also of the choice of line $s$ through $C$. The independence of the fourth harmonic point is called the *theorem of the complete quadrilateral*, and this shows that it is equivalent to $D_2$.

2β) Now we prove that $D_2$ *implies all properties of harmonic quadrilaterals, namely exchangeability of pairs and invariance under projection.*

For the sake of brevity we write $\{ABCD\}$ to denote that $D$ is the fourth harmonic point of $A, B, C$.

(a) To begin with, $\{ABCD\}$ implies $\{BACD\}$ and $\{ABDC\}$.

The first is trivial because of the symmetric roles of $A$ and $B$ in the quadrilateral construction. The second is proved as follows.

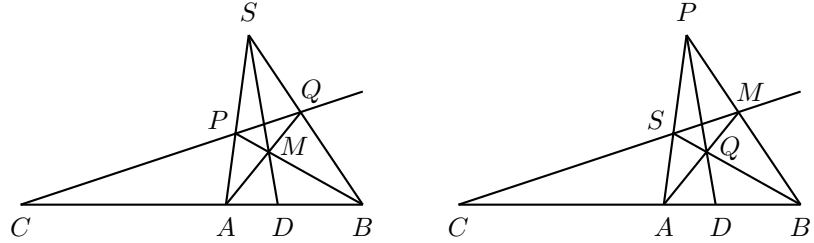

[By relabelling the first figure] there is a generating quadrilateral over $AB$ with adjacent vertex $D$ and adjacent side going through $C$, namely the quadrilateral $ABPQ$ in the second figure, with adjacent vertices $S, M, D$.

(b) **Exchangeability of pairs.** $\{ABCD\}$ *implies* $\{CDAB\}$.

I have to show that *there is a generating quadrilateral over $CD$ with adjacent vertex $B$ and adjacent side through $A$*. To do this I must prove, for example, that $CRMD$ is such a quadrilateral, where $R = CO \cap PB$.

The points $A, B, C, D$ and $P, Q, S, M$ correspond to those in the first figure above [but with $S$ now replaced by $O$]. Let $CM \cap AO = F$. Then $\{APOF\}$, because the quadrilateral $APBQ$ stands over $F$.

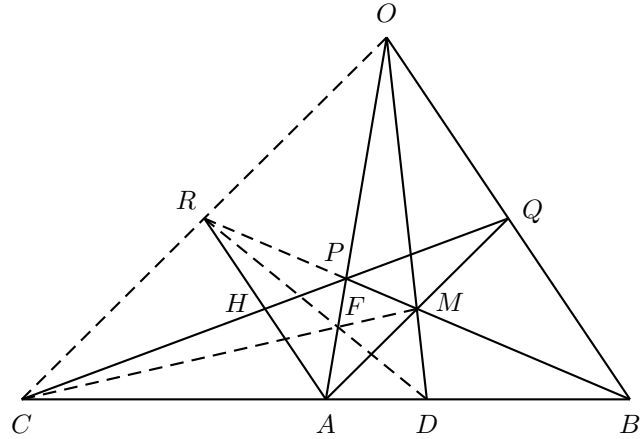

In quadrilateral $APRC$, let $H = RA \cap CP$. Then $H, F, B$ are collinear, since one can also generate $\{APOF\}$ with the quadrilateral $APCR$.

The relation $\{ABCD\}$, generated by the quadrilateral $ABHP$ instead of $ABPQ$, shows that $R, F, D$ are collinear. The quadrilateral $CDRM$ therefore has adjacent vertices $B, O, F$; and $O, F, A$ are collinear. That is, we also have $\{CDBA\}$. [Hence $\{CDAB\}$ by part (a).]                                                                Q.E.D.

We therefore have

$$\{ABCD\} \longrightarrow \{CDAB\}, \{DCAB\}, \{BACD\}, \{BADC\},$$
$$\{ABDC\}, \{CDBA\}, \{DCBA\}.$$

(c) **Invariance of the harmonic quadruple under projection.**
  Suppose that $\{ABCD\}$ is projected from $O$ onto $g'$.
  Claim: *The relation $\{A'B'C'D'\}$ holds.*

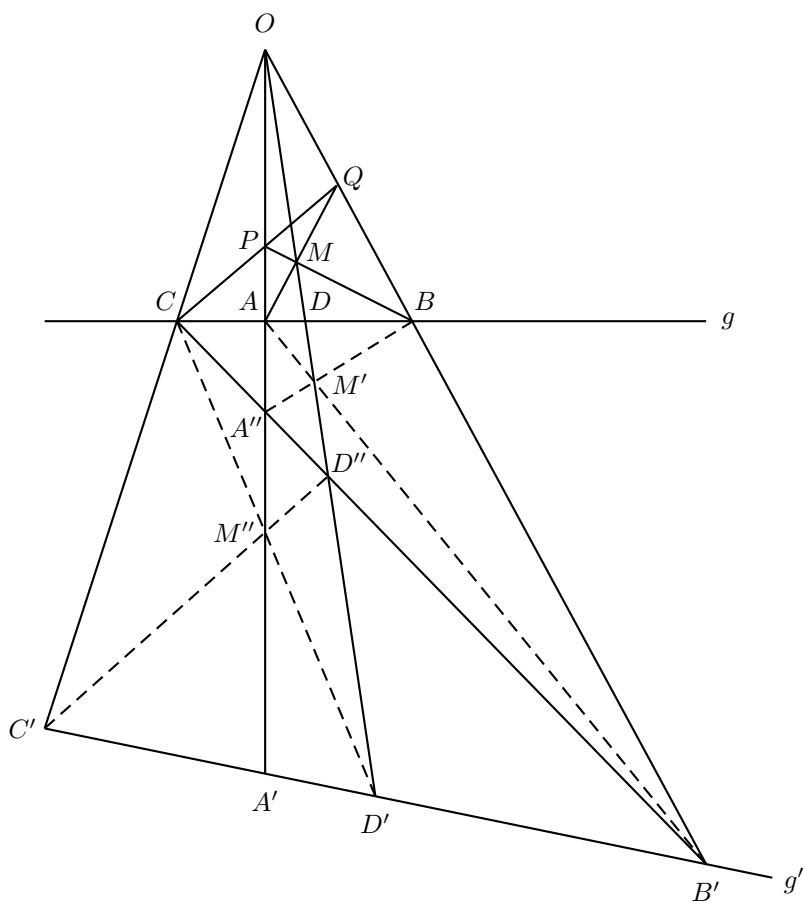

Proof. First of all, $O$ and $M$ are collinear with $M' = A''B \cap AB'$, since $\{ABCD\}$ may be generated by the quadrilateral $ABPQ$ as well as by the quadrilateral $ABA''B'$. Thus $\{A''B'CD''\}$, since it is generated by the quadrilateral $A''B'BA$.

It follows that $\{CD''A''B'\}$ by (b). This quadruple is generated by the quadrilateral $C'D'CD''$, and hence $O$ and $A$ are collinear with $M'' = C'D'' \cap CD'$. But then $\{C'D'B'A'\}$ is also a harmonic quadruple, hence so too is $\{A'B'C'D'\}$. Q.E.D.

3. **Theorem $D_2$ implies theorem $D_1$.**[5]

One has to prove that the perspective triangles $\Delta 123$ and $\Delta 1'2'3'$ in the accompanying figure lie axially, that is, that $A, B, C$ are collinear.

To prove this one introduces the points

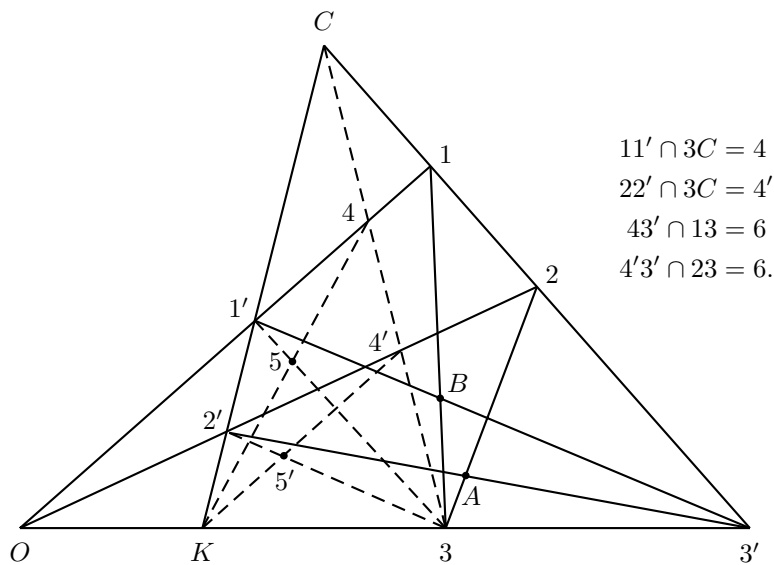

$$11' \cap 3C = 4$$
$$22' \cap 3C = 4'$$
$$43' \cap 13 = 6$$
$$4'3' \cap 23 = 6.$$

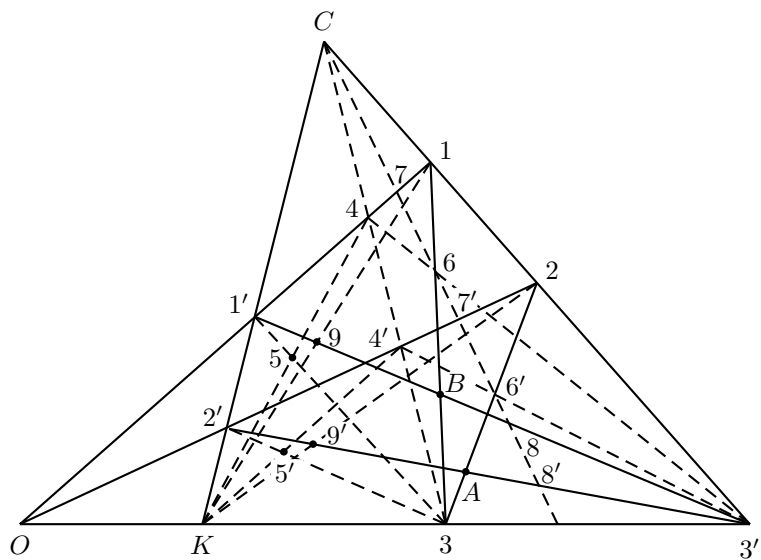

Theorem $D_2$ for triangles $\Delta 123$ and $\Delta 44'3'$ says that $C, 6, 6'$ are collinear.

---

[5]Translator's note. Moufang's proof of this theorem also makes the implicit assumption that the diagonal points of a complete quadrilateral are collinear, and I believe it is still open whether $D_2$ implies $D_1$ in an arbitrary projective plane. For a proof under the explicit assumption that the diagonal points are not collinear, see Pickert *Projective Ebenen* (Springer-Verlag 1955) p. 193 or Heyting *Axiomatic Projective Geometry* (North-Holland 1963) p. 55.

We now construct the points

$$7 \equiv 11' \cap 66'$$
$$7' \equiv 22' \cap 66'$$
$$8 \equiv 1'3' \cap 66'$$
$$8' \equiv 2'3' \cap 66'$$
$$9 \equiv K1 \cap 1'3'$$
$$9' \equiv K2 \cap 2'3'.$$

The points $C, 8, 8'$ are collinear, since $C, 6, 6'$ are collinear and by construction 8 is collinear with 6 and $6'$, and $8'$ likewise.

Then it follows by Theorem $D_2$ in $\Delta 1'3'4$ and $\Delta K13$ that $O, 9, 6$ are also collinear.

By Theorem $D_2$ in $\Delta 2'4'3'$ and $\Delta K32$ we have $O, 9', 6'$ collinear. And Theorem $D_2$ in $\Delta 12K$ and $\Delta 1'2'3'$ shows that $C, 9, 9'$ are collinear.

The generating quadrilateral $4133'$ gives $\{41O7\}$, and the quadrilateral $4'233'$ likewise gives $\{4'2O7'\}$.

One now projects $\{41O7\}$ from 6 onto $1'3'$, obtaining the harmonic quadruple $\{3'B98\}$; and $\{4'2O7\}$ from $6'$ onto $2'3'$, obtaining the harmonic quadruple $\{3'A9'8'\}$.

The harmonic quadruples $\{3'B98\}$ and $\{3'A9'8'\}$ with the common point $3'$ are therefore in perspective. That is, since $C, 9, 9'$ are collinear, and also $C, 8, 8'$, so too are $C, B, A$. Q.E.D.

4. The converse of the theorem just proved from $D_2$, namely that *axially lying triangles $\Delta 123$ and $\Delta 1'2'3'$, with $3'$ on 12, are also in perspective*, follows in a similar way by reversing the argument.

One can also reach the conclusion by an indirect argument. Suppose that the three lines $11'$, $22'$, $33'$ have intersections

$$O_3 = 11' \cap 22', \quad O_1 = 22' \cap 33', \quad O_2 = 33' \cap 11' \quad \text{with } O_3 \neq O_2.$$

Connecting $O_2$ with $2'$ then determines a unique point $\overline{2} \neq 2$ on $13'$, and connecting $\overline{2}$ with 3 determines a unique point $\overline{A} \neq A$ on $2'3'$. The triangles $\Delta 1'2'3'$ and $\Delta 1\overline{2}3$ then satisfy the hypotheses of the previous section and hence are axial, that is, $C, B, \overline{A}$ are collinear. But $C, B, A$ are also collinear by hypothesis, hence $A = 2'3' \cap BC$ and $\overline{A} = 2'3' \cap BC$, so $A = \overline{A}$. Retracing our steps, we then find $\overline{2} = 2$ and hence $O_3 = O$. Q.E.D.

## 3.2 The Hilbert segment calculus

To understand the meaning of $D_2$, we construct a plane geometry in which the plane axioms of incidence and order, the parallel postulate, and $D_2$ hold. In this geometry we use $D_2$ to introduce a segment calculus independent of congruence axioms and continuity, that is, we construct a "number system" geometrically.

Sum and product are defined as geometric operations.[6] It is shown that each computation rule is equivalent to an incidence theorem, and the rules provable in our geometric system turn out to be the axioms of an alternative field.

We remark that the parallel postulate is used only to introduce ideal points in the plane as simply as possible. Actually, the incidence axioms are needed only in their projective form (that is, two distinct points determine a unique line, and two distinct lines determine a unique point).

1. **"Addition"**

    1) In what follows we distinguish two lines, intersecting at $O$, as the $x$- and $y$-axes, and a unit point on each.

    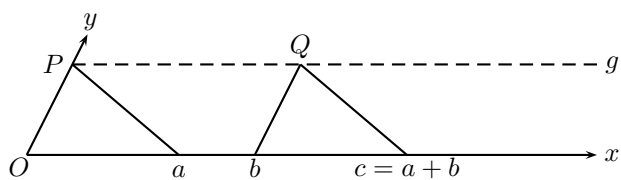

    The "sum" $a + b$ of segments [ending at] $a$ and $b$ is defined as the segment with endpoint $c$ determined as follows. Let $g$ be any parallel to the $x$-axis. Through $b$ draw a parallel to the $y$-axis, meeting $g$ at $Q$. Then through $Q$ draw a parallel to $Pa$, meeting the $x$-axis at $c$. Thus $c$ at first seems to depend on the special choice of $g$. However, this is not the case, because if $g'$ is another parallel to the $x$-axis we have the following picture.

    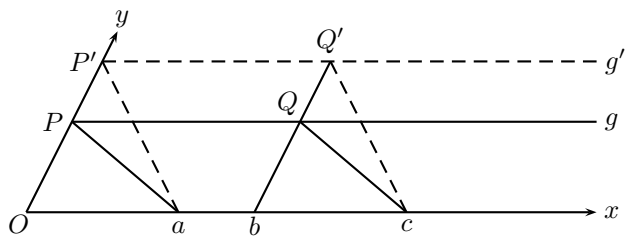

    The triangles $\Delta P'Pa$ and $\Delta Q'Qc$ are in perspective by construction, hence by the little Desargues theorem they are also axial, so the sum $c$ is uniquely determined.

    2) $a + b = b + a$

    To prove this we consider the following sketch.

    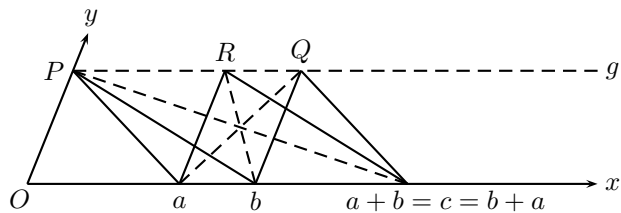

---

[6]Translator's note. Moufang uses the symbols $\hat{+}$ and $\circ$ for these operations, but I have opted for the usual symbols, since there seems to be no danger of confusion.

Hypothesis: $Pa \parallel Qc$,   $aR \parallel bQ$.

Claim: $Pb \parallel Rc$.

Proof. $\Delta PRa$ is axial with $\Delta cbQ$, hence by $D_2$ they are in perspective. Hence $\Delta PbQ$ is in perspective with $\Delta cRa$, so by $D_2$ they are axial. That is, $PB \parallel Rc$.                                     Q.E.D.

3) $(a+b)+c = a+(b+c)$

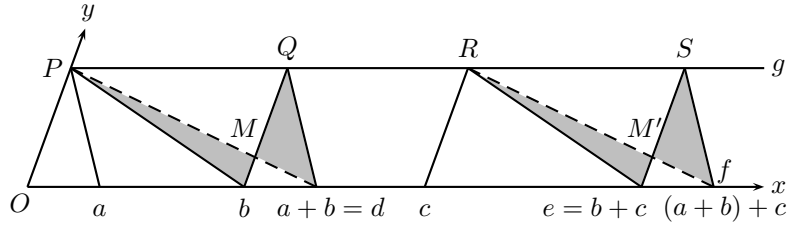

Hypothesis: $Pb \parallel Re$,   $Pd \parallel Rf$,   $Qb \parallel Rc \parallel Se \parallel y$-axis.

Claim: $Pa \parallel Sf$, or $Qd \parallel Sf$.

Proof. $\Delta PMb$ is axial with $\Delta RM'e$, since all corresponding sides are parallel. The little Desargues theorem implies that $MM' \parallel x$-axis, and hence $\Delta QMd$ is in perspective with $\Delta SM'f$ [where $f$ equals $(a+b)+c$]. Then it follows from little Desargues again that $Qd \parallel Sf$ [and hence $f$ also equals $a+(b+c)$].                     Q.E.D.

Thus the associative law holds for segment addition.

4) The **zero element** is the segment $a = 0$ whose endpoints are both at the origin.

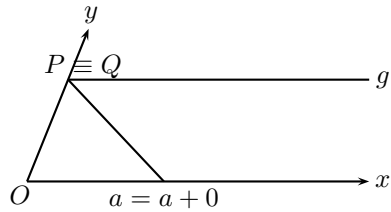

We have $a + 0 = a = 0 + a$.

The **additive inverse**[7] $-a$ of $a$ satisfies $(-a) + a = 0 = a + (-a)$. For its construction see the sketch!

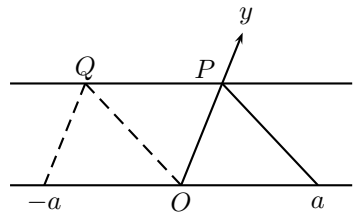

Thus all the computation rules for addition are valid.

---
[7]Translator's note. Moufang uses the symbol $\tilde{a}$ for the additive inverse of $a$, but again I have opted for the usual notation.

*The roles of the two axes are interchangeable.* To switch from the $x$-axis to the $y$-axis we must change our previous notation, writing $x_a$ for $a$ and $x_e$ for the unit point $e$ on the $x$-axis, etc. For each $a$, $y_a$ is defined by $y_e x_e \parallel y_a x_a$.

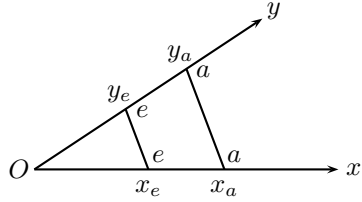

**Interchangeability of the axes with respect to sum** results as follows. Construct $x_c$ by means of the auxiliary line $g$, and $y_c$ by means of the auxiliary line $h$. Then

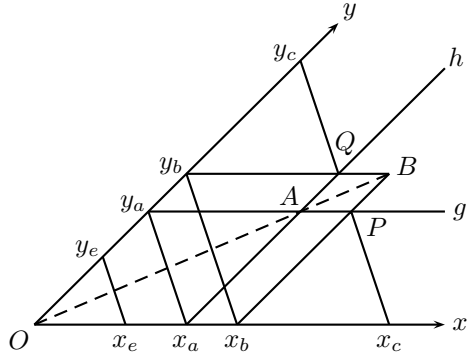

$Oy_a \parallel x_a A \parallel x_b P$
$x_e y_e \parallel x_a y_a \parallel P x_c \parallel x_b y_b$
$Oy_a \parallel x_a Q \parallel x_b P$
$x_a y_a \parallel Q y_c.$

Let $B = Qy_b \cap Px_b$. Then the triangles $\Delta x_a y_a A$, $\Delta x_b y_b B$, which lie axially by construction, are in perspective by $D_2$. That is, $O, A, B$ are collinear. Applying $D_2$ to the perspective triangles $\Delta x_a O y_a$, $\Delta QBP$ gives $x_a y_a \parallel PQ$, hence $x_c, P, Q, y_c$ lie on a parallel to $x_e y_e$. Q.E.D.

2. **"Multiplication"**

   1) **Definition** of $x_a x_b = x_c$.

      We are given $x_e y_e \parallel x_a y_a$. Through $y_b$ one draws the parallel to $y_e x_a$; it cuts the $x$-axis at $x_c = x_a x_b$. We abbreviate this by $ab = c$.

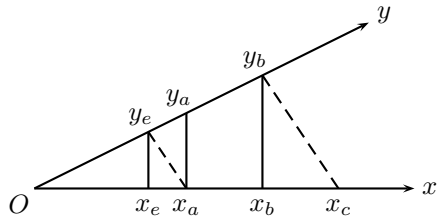

The sketch shows the construction, with $y_e x_a \parallel y_b x_c$.

2) **Interchangeability of the axes with respect to product.**
By construction we have

$$x_e y_e \parallel x_a y_a \parallel x_b y_b \parallel x_c y_c$$
$$x_a y_e \parallel x_c y_b.$$

We have to prove that $x_e y_a \parallel x_b y_c$.

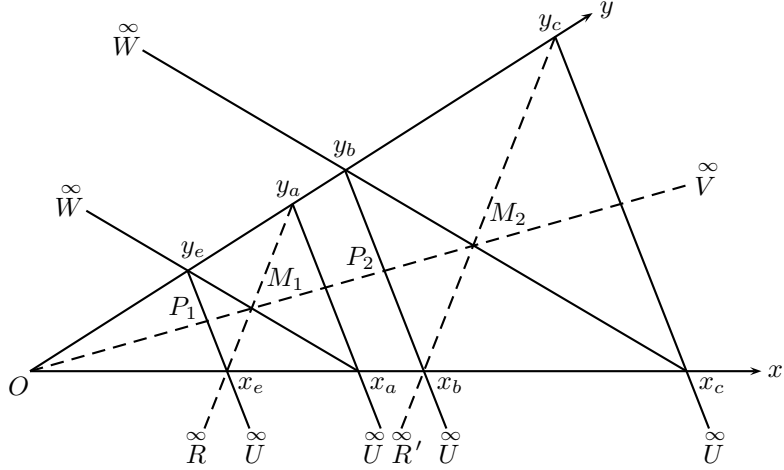

Proof. We have two quadrilaterals with $O$ as common adjacent vertex, namely $y_e x_e y_a x_a$ and $y_b x_b y_c x_c$. The point at infinity $\overset{\infty}{U}$ of $x_e y_e$ is also a common adjacent vertex.

Let $M_1 = y_e x_a \cap y_a x_e$, $M_2 = y_c x_b \cap y_b x_c$.

Then $O, M_1, M_2$ are collinear by the quadrilateral theorem.

$\{y_e x_e P_1 \overset{\infty}{U}\}$ and $\{y_b x_b P_2 \overset{\infty}{U}\}$ are harmonic quadruples. I project them from $M_1$ and $M_2$ respectively onto the line at infinity and obtain—when $\overset{\infty}{W}$ denotes the point at infinity on $x_a y_e$, $\overset{\infty}{R}$ the point at infinity on $x_e y_a$, and $\overset{\infty}{R}'$ the point at infinity on $x_b y_c$—the harmonic quadruples $\{\overset{\infty}{W}\overset{\infty}{R}\overset{\infty}{V}\overset{\infty}{U}\}$ and $\{\overset{\infty}{W}\overset{\infty}{R}'\overset{\infty}{V}\overset{\infty}{U}\}$.

It follows that $\overset{\infty}{R} = \overset{\infty}{R}'$, that is, $x_e y_a \parallel x_b y_c$.

Thus the axes are interchangeable. Q.E.D.

3) **Existence of identity and inverse elements.**
The identity element is $x_e$. The proof that

$$x_a x_e = x_a = x_e x_a$$

is clear from the pictures.

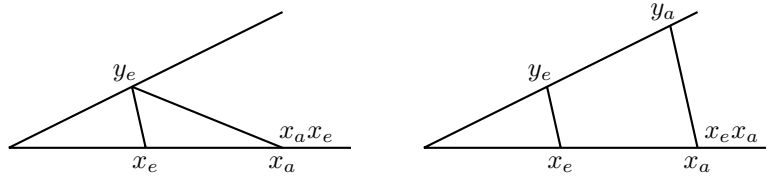

The inverse[8] $x_{a^{-1}}$ of $x_a$ must satisfy

$$x_a x_{a^{-1}} = x_e = x_{a^{-1}} x_a.$$

We determine $x_{a^{-1}}$ and $y_{a^{-1}}$ by the following construction. Through $x_e$ draw the parallel to $x_a y_e$. It cuts the $y$-axis at $y_{a^{-1}}$, whence $x_{a^{-1}}$ is determined by $y_{a^{-1}} x_{a^{-1}} \parallel y_e x_e$.

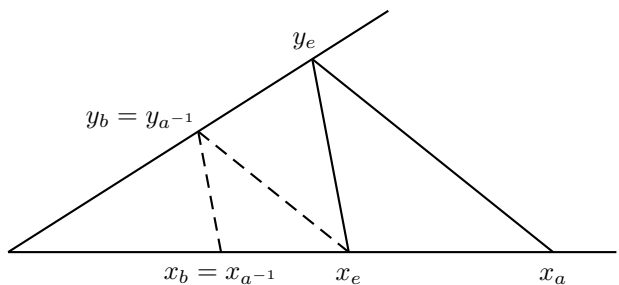

Then indeed $x_a x_{a^{-1}} = x_e$.

Also $y_{a^{-1}} y_a = y_e$, hence also $x_{a^{-1}} x_a = x_e$.

3. **Computation rules for multiplication.**

1) $a(ab) = (aa)b$

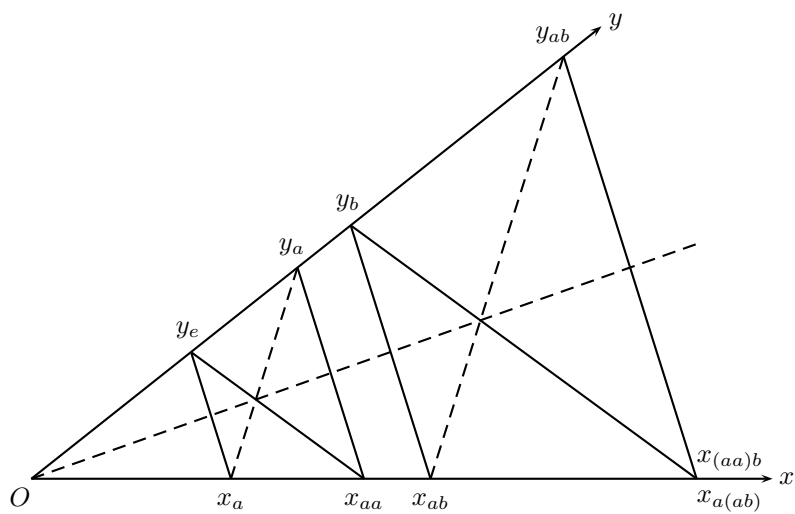

The figure is identical—apart from notation—with the figure for 2), so the proof is analogous to the one there.

---

[8]Translator's note. Moufang denotes the multiplicative inverse of $x_a$ by $x_{\tilde{a}}$. Since $x_a$ replaces the previous notation $a$, it seems to me more logical to denote the multiplicative inverse by $x_{a^{-1}}$. Indeed, Moufang herself switches to the latter notation before long.

2) $b(aa) = (ba)a$

By construction we have (see sketch)

$$x_e y_e \parallel x_a y_a,$$
$$y_e x_b \parallel y_a x_c, \quad x_e y_a \parallel x_a y_f, \quad y_e x_b \parallel y_a x_c \parallel y_f x_d,$$

[where $c = ba$, $d = b(aa)$, $f = aa$] and we have to show $y_e x_c \parallel y_a x_d$.

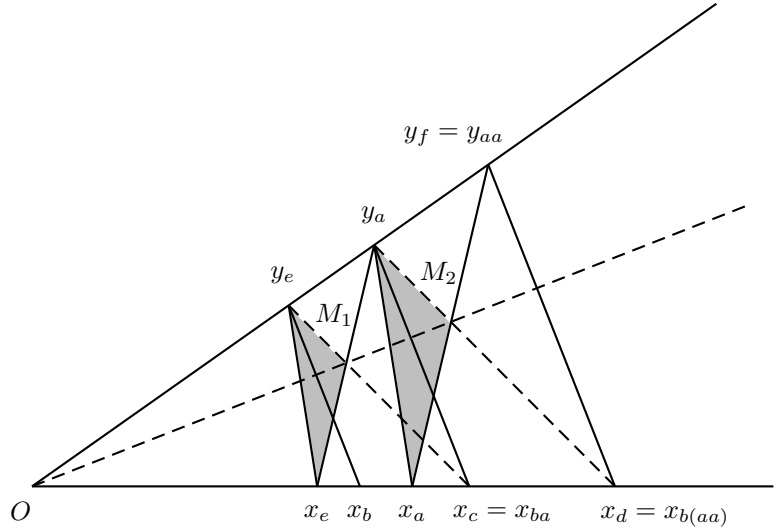

Proof. Let $M_1 = x_e y_a \cap x_c y_e$, $M_2 = x_a y_f \cap x_d y_a$.

$\triangle y_e x_e M_1$ and $\triangle y_a x_a M_2$ lie axially, because all corresponding sides are parallel. Since $y_a$ lies on $x_a M_1$, it follows by $D_1$ that the triangles are in perspective, and hence $O, M_1, M_2$ are collinear.

Thus $\triangle y_a M_1 x_c$ is in perspective with $\triangle y_f M_2 x_d$ by the above. Since $y_a$ lies on $M_2 x_d$, $D_1$ applies again, and it shows that the triangles are axial. That is, $y_a x_c \parallel y_f x_d$. Q.E.D.

3) $a^{-1}(ab) = (a^{-1}a)b = b$
$(ba)a^{-1} = b(aa^{-1}) = b$

The first of these is trivial by the interchangeability of the axes.

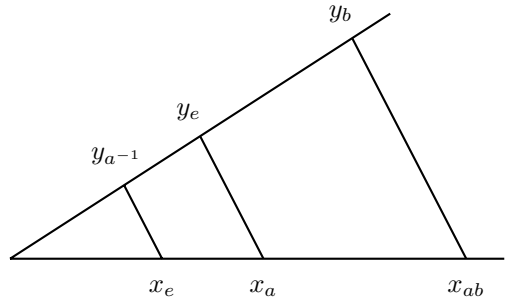

To prove the second equation one has to show that the figure in which

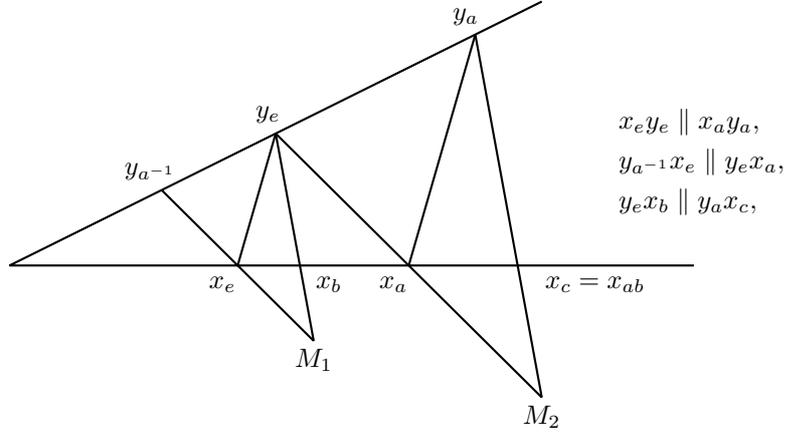

$$x_e y_e \parallel x_a y_a,$$
$$y_{a^{-1}} x_e \parallel y_e x_a,$$
$$y_e x_b \parallel y_a x_c,$$

where $c = ab$, also has

$$y_{a^{-1}} x_b \parallel y_e x_c.$$

Proof. Let $m_1 = y_{a^{-1}} x_e \cap y_e x_b$, $M_2 = y_e x_a \cap y_a x_c$.

Then $D_1$, applied to the the axial triangles $\Delta y_e x_e M_1$ and $\Delta y_a x_a M_2$, shows that $O, M_1, M_2$ are collinear. The assertion then follows by applying $D_1$ to the perspective triangles $\Delta y_{a^{-1}} x_b M_1$ and $\Delta y_e x_c M_2$.
Q.E.D.

4) **Does $(ab)c = a(bc)$ hold?**

One cannot prove this from $D_2$, but it does follow from two applications of the general Desargues theorem ($D_0$).

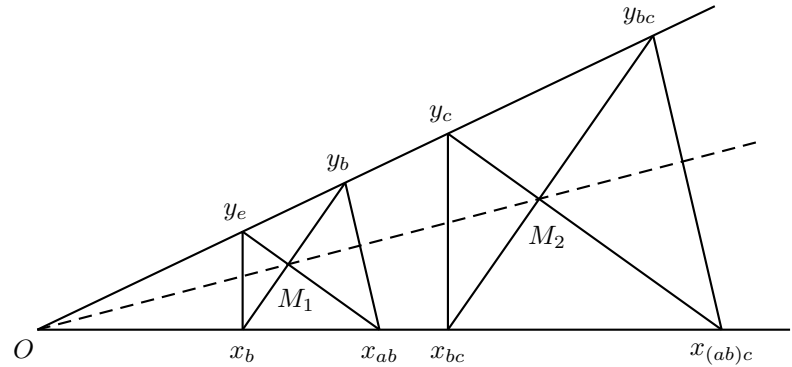

Since $O, M_1, M_2$ are collinear by $D_0$, it follows that

$$y_b x_{ab} \parallel y_{bc} x_{(ab)c}.$$

The attempt to prove the general associative law from $D_2$ fails, as it must, for reasons we explain later.

5) **Does $ab = ba$ hold?**

This law immediately follows from the theorem of Pappus, as one sees from the figure.

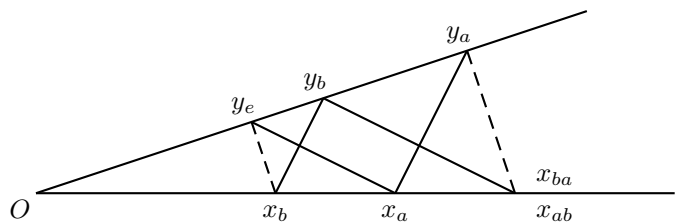

We proved earlier that Desargues' theorem follows from Pappus' theorem (Hessenberg), but not conversely (using Hilbert's number system). Since $D_2$ is a special case of Desargues' theorem, it also cannot follow from Pappus' theorem. Later we shall even prove, as already mentioned, that $D_2$ does not imply $D_0$.

6) **The distributive laws.**

$\alpha$) $a(b+c) = ab + ac$
By hypothesis:[9] $x_e y_a \parallel x_b y_{ab} \parallel x_c y_{ac} \parallel Q y_d \parallel P x_f$.

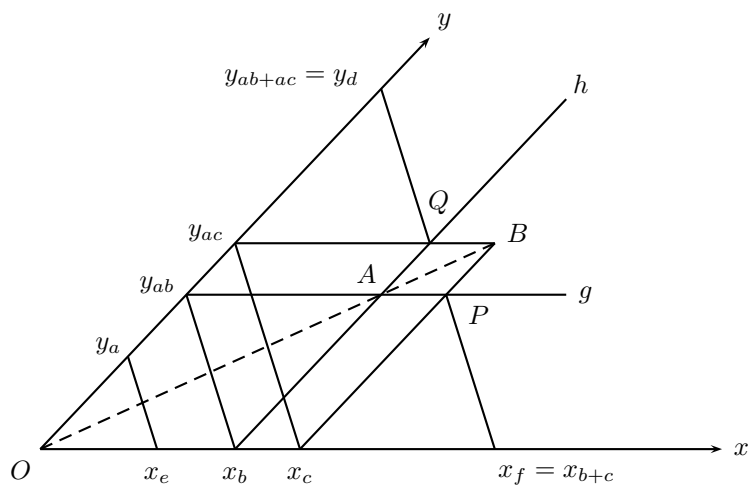

It is to be proved that $x_f, P, Q, y_d$ are collinear. The figure is the same as the figure used earlier to prove interchangeability of the axes with respect to addition, and the proof is word for word the same as before.

---

[9]Translator's note: Moufang's figure for this proof is somewhat confusing, since she wants to compare it with a geometrically identical figure used earlier. She does so by putting the old labels (in brackets) next to the new ones. To my mind this is a distraction, so I have taken the liberty of omitting the old labels from the figure, and have also omitted a couple of Moufang's sentences that refer to it.

β) $(b+c)a = ba + bc$.

I construct the segments on the $y$-axis, where, for the sake of brevity, only the subscripts are written. We have

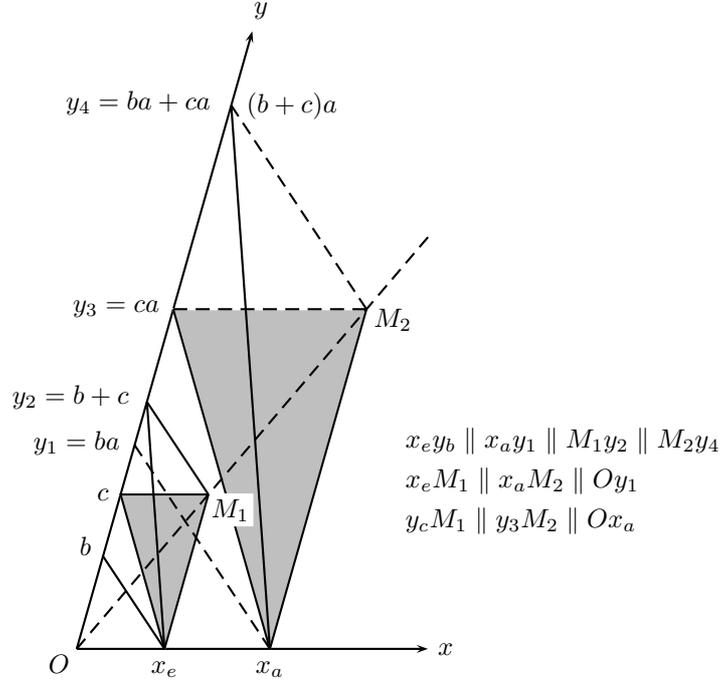

$$x_e y_b \parallel x_a y_1 \parallel M_1 y_2 \parallel M_2 y_4$$
$$x_e M_1 \parallel x_a M_2 \parallel O y_1$$
$$y_c M_1 \parallel y_3 M_2 \parallel O x_a$$

and we have to prove that $x_e y_2 \parallel x_a y_4$.

Proof. By construction, $\Delta x_e y_c M_1$ and $\Delta x_a y_3 M_2$ are axial, so by theorem $D_1$ they are also in perspective, that is, $O, M_1, M_2$ are collinear. Then the assertion follows from the same theorem in the triangles $\Delta x_e y_2 M_1$ and $\Delta x_a y_4 M_2$.    Q.E.D.

All the computation rules for the sum and product operations of an alternative field are thereby proved, using only $D_2$ and its consequences.

**Overview and classification**

Suppose we assume projective incidence properties, that is, the ordinary incidence axioms plus the parallel postulate.
Then the additional assumptions

$$\begin{array}{rcl} D_2 & \longrightarrow & \text{all computation rules of an alternative field} \\ D_0 & \longrightarrow & \text{all computation rules of a skew field} \\ \text{Pappus} & \longrightarrow & \text{all computation rules of a field.} \end{array}$$

Hilbert based the segment calculus on the general Desargues' theorem. The result that the weaker Desargues' theorem $D_2$ yields the computation rules of an alternative field is due to R. Moufang. *Math. Ann.* 106 (1932). p. 755ff.

4. **Introducing plane coordinates and coordinate geometry.**

   A point $P$ is represented by a number pair $(x, y)$.

   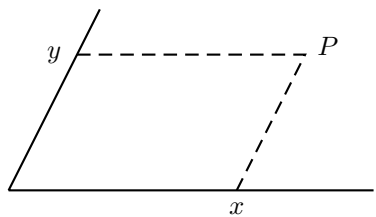

   In order to present the equation of a line, we investigate several cases.

   1) A line $g$ through $O$.

   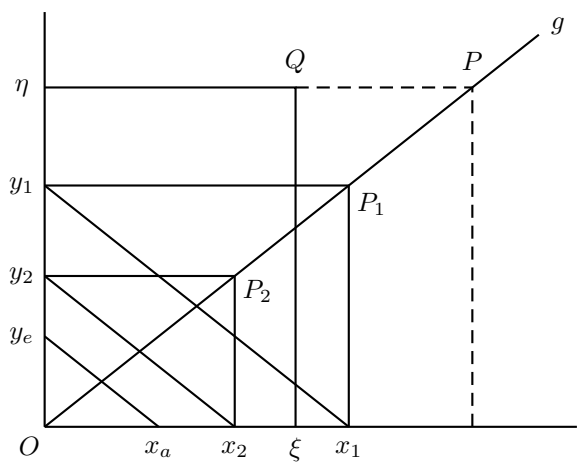

   Let $P_1$ and $P_2$ be arbitrary points on $g$. Then, in triangles $\Delta y_1 x_1 P_1$ and $\Delta y_2 x_2 P_2$, $D_2$ obviously gives $x_1 y_1 \parallel x_2 y_2$.

   Suppose that the parallel to $x_1 y_1$ through $y_e$ cuts the $x$-axis at $x_a$. Then $x_1 = a y_1$ and $x_2 = a y_2$, and hence for all points on $g$ we have $x = ay$ or $y = bx$.

   For $a = 0$ one obtains all the points of the $y$-axis; for $b = 0$ the points of the $x$-axis.

   Conversely, the coordinates $(\xi, \eta)$ of a point $Q$ not on $g$ do not satisfy the equation $x = ay$, since $\xi \neq a\eta$ by the uniqueness of the product operation.

   2) A line $g$ not through $O$.

   Let $P$ be an arbitrary point on $g$ with coordinates $(x, y)$. One draws parallels to $g$ through the unit point and the point $x$ on the $x$-axis. They cut the $y$-axis in the points $a$ and $ax$ respectively. Then it follows from the definition of addition that

   $$y_0 = y + ax$$

   for all points on $g$. Here $a$ and $y_0 = c$ are fixed.

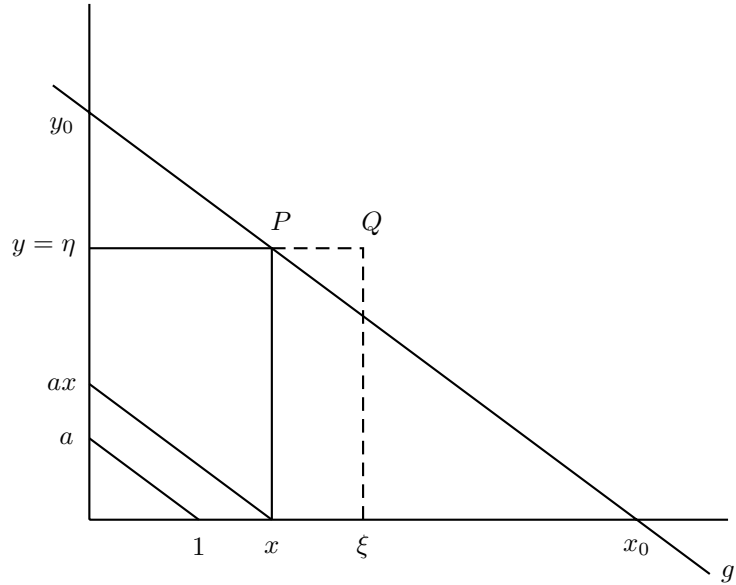

Conversely, the coordinates $(\xi, \eta)$ of a point $Q$ not on $g$ do not satisfy this equation, otherwise for $\eta = y$, $\xi \neq x$ we should have

$$y_0 = ax + y = a\xi + \eta, \quad \text{hence} \quad a(x - \xi) = 0,$$

which is false because $a \neq 0$ and $x \neq \xi$.

Similarly, one can obtain the equation of $g$ in the form

$$x + by = \text{const}.$$

Thus we have *two types of lines*:

$$x = c \quad \text{and} \quad y + ax = c \text{ with } a \neq 0,$$

or

$$y = c \quad \text{and} \quad x + by = c.$$

The lines $x = c$ are not included in the form $y + ax = c$.

## 3.3   Coordinate geometry over an alternative field ("harmonic geometry")

In Section 2.5 it was proved that the incidence theorems and Desargues' theorem hold in any coordinate geometry over a skew field. In this section it will be shown that the incidence theorems and the special Desargues theorem $D_2$ hold in any coordinate geometry over an alternative field. (Such a geometry is called *harmonic* because $D_2$ allows harmonic quadruples to be defined and their main properties proved.) Thus we are about to follow the path opposite to that in the previous section, by starting with the analytic definition of points and lines.

1. **Definition of points and lines**

   $P = (x, y)$; $P_1 = P_2$ means that $x_1 = x_2$ and $y_1 = y_2$.

   A line is represented by either one of two essentially different equations:
   $$x - b = 0 \quad \text{or} \quad y + ax - b = 0.$$

   Two lines are called equal when all their coefficients coincide.

   It will be shown below that the computation rules of an alternative field suffice to prove the projective incidence theorems for our geometry (in the affine form), and to prove the special Desargues theorem $D_2$.

2. **Incidence theorems**

    1) **Intersection of two lines**

        Two distinct lines either have exactly one point of intersection or else are parallel. Lines $g_1$ and $g_2$ are called parallel when both have the form $x - b = 0$ or both have the form $y + ax - b = 0$ with same $a$. Thus two lines of different types are never parallel.
        $$\begin{aligned} g_1 : & \quad x - b_1 = 0 \\ g_2 : & \quad x - b_2 = 0 \quad [\text{with } b_1 \neq b_2] \end{aligned}$$
        have no point of intersection, that is, they are always parallel.
        $$\begin{aligned} g_1 : & \quad x - b_1 = 0 \\ g_2 : & \quad y + a_2 x - b_2 = 0 \end{aligned}$$
        have exactly one intersection: $x = b_1$, $y = b_2 - a_2 x = b_2 - a_2 b_1$.
        $$\begin{aligned} g_1 : & \quad y + a_1 x - b_1 = 0 \\ g_2 : & \quad y + a_2 x - b_2 = 0 \end{aligned}$$
        For $a_1 = a_2$ we have $g_1 \parallel g_2$. For $a_1 \neq a_2$, subtraction gives
        $$(a_1 - a_2)x - b_1 + b_2 = 0.$$
        And $a_1 - a_2 \neq 0$, so $(a_1 - a_2)^{-1}$ exists. Multiplying on the left by it [and using $a^{-1}(ab) = b$ with $a = a_1 - a_2$ and $b = x$] gives
        $$x - (a_1 - a_2)^{-1}(b_1 - b_2) = 0,$$
        hence $x = (a_1 - a_2)^{-1}(b_1 - b_2)$, $y = b_1 - a_1 \{(a_1 - a_2)^{-1}(b_1 - b_2)\}$, [using the equation for $g_1$]. Thus we have found a point $P$ on $g_1$. $P$ also lies on $g_2$, as calculation [of $y + a_2 x - b_2$] shows:
        $$\begin{aligned} & b_1 - a_1\left\{(a_1-a_2)^{-1}(b_1-b_2)\right\} + a_2\left\{(a_1-a_2)^{-1}(b_1-b_2)\right\} - b_2 \\ &= b_1 - b_2 + (a_1 - a_2)\left\{(a_1 - a_2)^{-1}(b_1 - b_2)\right\} \\ &= b_1 - b_2 - (b_1 - b_2) = 0 \end{aligned}$$

        
        This proves the existence of an intersection point for $g_1$ and $g_2$. Now we must prove its uniqueness.

If $(x_1, y_1)$ and $(x_2, y_2)$ are two intersection points then we have

$$\left.\begin{array}{l} y_1 + a_1 x_1 - b_1 = 0 \\ y_2 + a_1 x_2 - b_1 = 0 \end{array}\right\} \quad y_1 - y_2 + a_1(x_1 - x_2) = 0 \qquad (1)$$

$$\left.\begin{array}{l} y_1 + a_2 x_1 - b_2 = 0 \\ y_2 + a_2 x_2 - b_2 = 0 \end{array}\right\} \quad y_1 - y_2 + a_2(x_1 - x_2) = 0 \qquad (2)$$

Then (1)–(2) gives $(a_1 - a_2)(x_1 - x_2) = 0$. [Since $a_1 \neq a_2$] this implies $x_1 = x_2$, and hence $y_1 = y_2$ from (1).

Thus in all cases where $g_1 \not\parallel g_2$ there is exactly one intersection point.
Q.E.D.

2) **Connecting two points by a line**

Suppose $P_1 \equiv (x_1, y_1)$ is different from $P_2 \equiv (x_2, y_2)$.

Thus if $x_1 = x_2$, $y_1 \neq y_2$.

In this case the connecting line is $x - b = 0$, with $b = x_1 = x_2$. The formula $y + ax - b = 0$ is impossible here, because if

$$y_1 + ax_1 - b = 0$$
$$y_2 + ax_2 - b = 0$$

then subtraction gives $y_1 - y_2 + a(x_1 - x_2) = 0$, with $x_1 - x_2 = 0$, and this contradicts $y_1 \neq y_2$.

When $x_1 \neq x_2$ I take the formula $y + ax - b = 0$. Then the two points give the two equations

$$y_1 + ax_1 - b = 0$$
$$y_2 + ax_2 - b = 0$$

and subtraction gives $y_1 - y_2 + a(x_1 - x_2) = 0$. Since $x_1 - x_2 \neq 0$, $(x_1 - x_2)^{-1}$ exists, and multiplying on the right by it gives:

$$(y_1 - y_2)(x_1 - x_2)^{-1} - a = 0, \quad \text{hence} \quad a = (y_1 - y_2)(x_2 - x_1)^{-1}.$$

[Substituting this value of $a$ in the first equation gives]

$$b = y_1 + \left\{(y_1 - y_2)(x_2 - x_1)^{-1}\right\} x_1.$$

These values of $a$ and $b$ also satisfy the second equation, because

$$y_2 + \left\{(y_1 - y_2)(x_2 - x_1)^{-1}\right\} x_2 = y_1 + \left\{(y_1 - y_2)(x_2 - x_1)^{-1}\right\} x_1$$

Thus it is shown that there is always exactly one line through two different points. Q.E.D.

3) **The parallel postulate**

Let $g$ be $y + ax - b = 0$ and suppose that $P_1 \equiv (x_1, y_1)$ is a point not on $g$. We have to show that there is exactly one line $g'$ through $P_1$ with $g \parallel g'$.

Let $g'$ be $y + ax - b' = 0$. Since $P_1$ lies on $g'$, $b' = y_1 + ax_1$ is uniquely determined. Q.E.D.

If $g$ has the form $x - b = 0$, $g'$ is uniquely determined as $x - x_1 = 0$.

The line connecting a proper point $(x_1, y_1)$ to an improper point—given either by the collection of parallels $y + ax - b_i = 0$ to a line $y + ax - b = 0$, or by the collection of lines $x - b_j = 0$ parallel to a line $x - b = 0$—is thereby uniquely determined.

3. **Condition for three points to be collinear.**

The line $g$ through $P_1$ and $P_2$ is uniquely determined by 2) above. $P_3$ on $g$ therefore satisfies

$$y_3 + \{(y_1 - y_2)(x_2 - x_1)^{-1}\} x_3 - y_1 - \{(y_1 - y_2)(x_2 - x_1)^{-1}\} x_1 = 0 \quad (1)$$

First notice that this relation is symmetric in $P_1$ and $P_2$, because

$$y_1 - y_2 = \{(y_2 - y_1)(x_1 - x_2)^{-1}\}(x_1 - x_2)$$

implies that

$$y_1 + \{(y_1 - y_2)(x_2 - x_1)^{-1}\} x_1 = y_2 + \{(y_2 - y_1)(x_1 - x_2)^{-1}\} x_2 \quad (2)$$

Now we have to show that $P_3$ on $g$ implies $P_2$ on $P_1 P_3$, that is, that equation (1) also holds when the indices 2 and 3 are exchanged:

$$y_2 + \{(y_1 - y_3)(x_3 - x_1)^{-1}\} x_2 - y_1 - \{(y_1 - y_3)(x_3 - x_1)^{-1}\} x_1 = 0 \quad (3)$$

To do this, we deduce from (1) that

$$x_3 = \{(y_1 - y_2)(x_2 - x_1)^{-1}\}^{-1} (y_1 - y_3) + x_1$$
$$= \{(x_2 - x_1)(y_1 - y_2)^{-1}\} (y_1 - y_3) + x_1$$
$$x_3 - x_1 = \{(x_2 - x_1)(y_1 - y_2)^{-1}\} (y_1 - y_3)$$
$$(x_3 - x_1)^{-1} = (y_1 - y_3)^{-1} \{(y_1 - y_2)(x_2 - x_1)^{-1}\}.$$

Substituting in the left hand side of (3) gives

$$\{(y_1 - y_2)(x_2 - x_1)^{-1}\} x_2 + y_2 - y_1 - \{(y_1 - y_2)(x_2 - x_1)^{-1}\} x_1$$
$$= \{(y_1 - y_2)(x_2 - x_1)^{-1}\} (x_2 - x_1) + y_2 - y_1,$$

which indeed equals zero.

*The condition for collinearity of three points $P_1, P_2, P_3$ is therefore*

$$(y_1 - y_3)(x_3 - x_1)^{-1} = (y_1 - y_2)(x_2 - x_1)^{-1}$$

*and likewise for each permutation of the indices.*

There is an analogous condition for three lines

$$y + a_i x - b_i = 0, \quad i = 1, 2, 3$$

to pass through a point (likewise for each permutation of the indices):

$$(a_2 - a_1)^{-1}(b_2 - b_1) = (a_3 - a_1)^{-1}(b_3 - b_1)$$

It also follows that *if $P_3$ lies on $P_1P_2$ and $P_4$ lies on $P_1P_2$ then $P_3$ lies on $P_3P_4$*, because

$$(y_2 - y_1)(x_2 - x_1)^{-1} = (y_3 - y_1)(x_3 - x_1)^{-1}$$

and

$$(y_2 - y_1)(x_2 - x_1)^{-1} = (y_4 - y_1)(x_4 - x_1)^{-1}$$

then

$$(y_3 - y_1)(x_3 - x_1)^{-1} = (y_4 - y_1)(x_4 - x_1)^{-1} \qquad \text{Q.E.D.}$$

4. **Proof of the special Desargues theorem $D_2$.**

The proof breaks into five steps.

1) Let $g$ be any line through $O$, $P$ and $P'$ any points on $g$, $h$ any line through $P$, $Q$ the intersection of $h$ with the $x$-axis, $PQ \parallel P'Q'$, $k$ through $P$ parallel to the $x$-axis, $l$ through $Q$ parallel to $g$, $R$ the intersection of $k$ and $l$, $k'$ through $P'$ parallel to the $x$-axis, $l'$ through $Q'$ parallel to $g$, and $R'$ the intersection of $k'$ and $l'$.

   Claim: $O, R, R'$ are collinear.

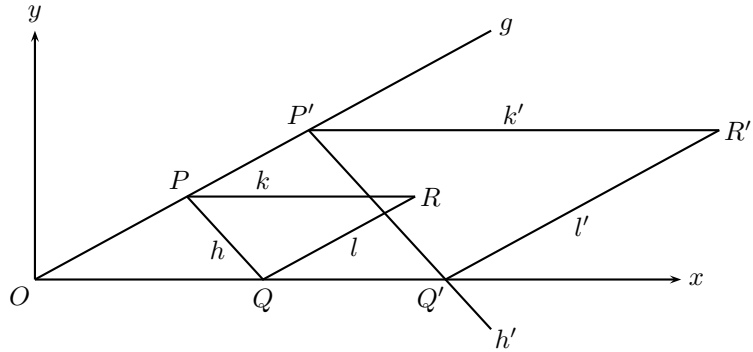

Proof. Assuming

$$\begin{aligned} g: &\quad y + ax = 0 \\ h: &\quad y + \alpha x - \beta = 0, \end{aligned}$$

it follows that $P = \big((a-\alpha)^{-1}\beta, -a\{(a-\alpha)^{-1}\beta\}\big)$, $Q = (\alpha^{-1}\beta, 0)$.
Similarly

$$\begin{aligned} k: &\quad y + a\{(a-\alpha)^{-1}\beta\} = 0 \\ l: &\quad y + ax - a(\alpha^{-1}\beta) = 0, \end{aligned}$$

imply $R = \big(\alpha^{-1}\beta + (a-\alpha)^{-1}\beta, -a\{(a-\alpha)^{-1}\beta\}\big) = (\xi, \eta)$ say.

Then

$$\xi\eta^{-1} = -\left[\alpha^{-1}\beta + (a-\alpha)^{-1}\beta\right]\left[\{(a-\alpha)^{-1}\beta\}^{-1} a^{-1}\right]$$
$$= -a^{-1} - (\alpha^{-1}\beta)\left[\{\beta^{-1}(a-\alpha)\} a^{-1}\right]$$
$$= -a^{-1} - (\alpha^{-1}\beta)\left[\beta^{-1} - (\beta^{-1}\alpha)a^{-1}\right]$$
$$= -a^{-1} - \alpha^{-1} + a^{-1}$$
$$= \alpha^{-1}, \quad \text{which is independent of } \beta.$$

Hence for $R' \equiv (\xi', \eta')$ we have $\eta'\xi'^{-1} = \eta\xi^{-1}$, and therefore the three points $O, R, R'$ are collinear by the criterion derived above,

$$(y_1 - y_2)(x_3 - x_1)^{-1} = (y_1 - y_2)(x_2 - x_1)^{-1},$$

in the special case where $(x_1, y_1) = (0,0)$, $(x_2, y_2) = (\xi, \eta)$, and $(x_3, y_3) = (\xi', \eta')$. \hfill Q.E.D.

2) In order to free the figure from its special position in the coordinate system, we make the coordinate transformation:

$$x = x' + c,$$
$$y = y' + d.$$

This is a collineation, that is, points go to points and lines to lines, and incidence of points and lines is preserved.

In fact

$$x + a = 0 \quad \text{goes to} \quad x' + a' = 0$$
$$y + ax + b = 0 \quad \text{goes to} \quad y' + ax' + b' = 0.$$

Consequently, the special Desargues theorem also holds when the center of perspective is an arbitrary point in the plane.

The same calculation may be carried out when $x$ and $y$ are exchanged, that is, when equations of lines are taken in the forms $y + b = 0$ and $x + ay + b = 0$. This proves those instances of the Desargues theorem $D_2$ in which the lines connecting vertices are parallel to the $x$-axis or the $y$-axis.

3) Now we free $QQ'$ from its special position in the coordinate system. This is done via the transformation

$$x' = x$$
$$y' = cx + y \quad \text{or} \quad y = -cx' + y'.$$

This is also a collineation, because

$$x + a = 0 \quad \text{goes to} \quad x' + a = 0,$$
$$y + ax + b = 0 \quad \text{goes to} \quad y' - cx' + ax + b = 0$$
$$\qquad\qquad\qquad \text{or} \quad y' + a'x' + b' = 0.$$

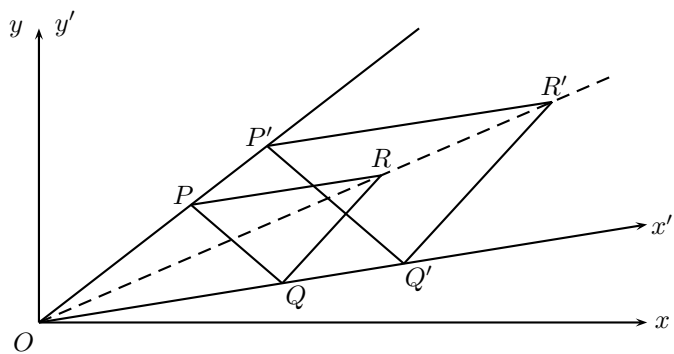

This collineation leaves the $y$-axis fixed, because for $x = 0$ we have $x' = 0$ and $y = y'$. However, the $x$-axis goes to the line $cx + y = 0$, which is a line through the origin. Then further application of the transformation 2) proves the Desargues theorem $D_2$ in the case where the axis of perspectivity is the line at infinity.

4) Now we transform the axis of perspectivity from infinity to the $y$-axis.

$$x' = x^{-1}$$
$$y' = yx^{-1}$$

is a collineation, because

$$x + a = 0 \quad \text{goes to} \quad x'^{-1} + a = 0 \quad \text{or} \quad x' + a^{-1} = 0$$
$$y + ax + b = 0 \quad \text{goes to} \quad y'x'^{-1} + ax'^{-1} + b = 0$$
$$\text{or} \quad y' + a + bx' = 0.$$

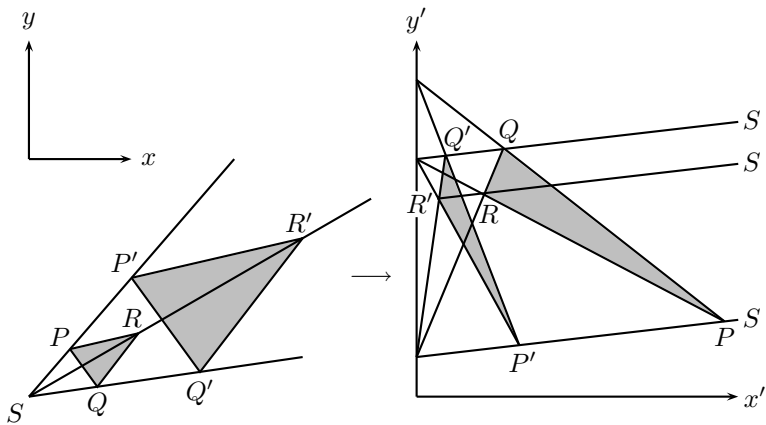

Parallel lines go to lines that meet on the $y$-axis. Indeed:

$$\left. \begin{array}{l} y + ax + b_1 = 0 \\ y + ax + b_2 = 0 \end{array} \right\} \quad \text{go to} \quad \left\{ \begin{array}{l} y' + b_1 x' + a = 0 \\ y' + b_2 x' + a = 0, \end{array} \right.$$

which intersect at $x' = 0$, $y' = -1$. That is, the intersection lies on the $y'$-axis.

    5) Finally, one can free the axis of perspectivity from its special position by the collineation

$$y' = y$$
$$x' = cy + x,$$

as in 3).

Thus $D_2$ holds in general, and the algebraic equivalent of $D_2$ is a plane coordinate geometry over an alternative field. The general Desargues theorem $D_0$, however, does not hold. If this theorem were valid, the configuration on page 80 would hold everywhere in the plane. If one constructs a plane coordinate geometry over the Cayley numbers in the manner given here, and takes $a = e_1$, $b = e_2$, $c = e$ (where $e_1, e_2, e$ are special Cayley numbers) then, since $(e_1 e_2)e \neq e_1(e_2 e)$, the configuration on page 80 is false, so the general Desargues theorem cannot hold, even though $D_2$ does.

5. The results above establish the following logical relationships between Pappus' theorem and the Desargues theorems $D_0$ and $D_2$.

$$\text{Pappus' theorem} \quad \overset{\longrightarrow}{\not\leftarrow} \quad D_0 \quad \overset{\longrightarrow}{\not\leftarrow} \quad D_2,$$

where $\longrightarrow$ means "implies" and $\not\longrightarrow$ means "does not imply".

It is important to note that the statement Desargues $\not\longrightarrow$ Pappus is valid assuming the plane incidence *and* order axioms, because the Hilbert number system is ordered. However, the statement $D_2 \not\longrightarrow D_0$ is proved assuming only the incidence axioms, because the Cayley number system is not ordered. It is an open problem whether $D_0$ is provable from $D_2$ assuming the plane incidence *and* order axioms. A solution of this problem would answer the question whether ordered alternative fields exist.

6. These investigations show the central importance of the segment calculus. It enables a partial algebraicization of a geometry in which certain configuration theorems hold. Investigating the structure of such a geometry by purely synthetic methods is generally very difficult. More precisely: if $S_1$ and $S_2$ are two configuration theorems and $S_2$ follows from $S_1$ on the basis of plane incidence (and order) then the proof can be carried out synthetically, given sufficient skill. But if $S_2$ does not follow from $S_1$, this can hardly be proved by purely synthetic methods. One therefore resorts to algebraic methods, whereby one constructs a model coordinate geometry in which $S_1$ is true and $S_2$ is false.

The segment calculus achieves the passage from synthetic geometry to coordinate geometry. The computation rules for the coordinates depend on the nature of the basic configuration theorem $S_1$. This gives a systematic approach to the very difficult and obdurate problem of classifying plane configuration theorems.

## 3.4   Foundations of projective geometry

1. *Pappus' theorem is equivalent to the fundamental theorem of projective geometry.*

   In Section 2.6 (page 55) we proved that Pappus' theorem follows from the fundamental theorem of projective geometry. We now prove the converse.

   From Section 3.2 (page 73) we know that Pappus' theorem gives a segment calculus over a field, whose elements serve as coordinates in an analytic geometry. This allows us to introduce projective coordinates in the plane and to define projective maps. Here it suffices to do this for single figures. A mapping is then projective if and only if it preserves the cross ratio ($CR$). That is, a quadruple $A, B, C, P$ can be mapped projectively to $A', B', C', P'$ if and only if

   $$CR = \frac{x_3 - x_1}{x_3 - x_2} : \frac{x_4 - x_1}{x_4 - x_2} = \frac{x_3' - x_1'}{x_3' - x_2'} : \frac{x_4' - x_1'}{x_4' - x_2'},$$

   where the $x_i$ are the coordinates of $A, B, C, P$ and the $x_i'$ are those of $A', B', C', P'$.

   When $A, B, C$ and $A', B', C'$ are fixed, the $P'$ corresponding to a given $P$ is uniquely determined. If $x_4' = x'$ and $x_4 = x$, while the others remain fixed then the above equation may be solved for $x'$ as

   $$x' = \frac{\alpha x + \beta}{\gamma x + \delta}, \quad \text{or} \quad \begin{cases} \rho \xi_1' = \alpha \xi_1 + \beta \xi_2 \\ \rho \xi_2' = \gamma \xi_1 + \delta \xi_2 \end{cases}$$

   in the homogeneous coordinates $\xi_1, \xi_2$. Thus a projective map is described by a linear transformation. Each linear transformation leaves the $CR$ invariant, as one may compute directly.

   Suppose $P$ is given and the numbers $x_i, x_i'$ come from a field. Then the linear transformation of a line onto itself,

   $$x' = \frac{\alpha x + \beta}{\gamma x + \delta},$$

   is uniquely determined by three corresponding point pairs. If

   $$\frac{\xi_1}{\xi_2} \text{ goes to } \frac{\xi_1'}{\xi_2'}, \quad \frac{\eta_1}{\eta_2} \text{ goes to } \frac{\eta_1'}{\eta_2'}, \quad \frac{\zeta_1}{\zeta_2} \text{ goes to } \frac{\zeta_1'}{\zeta_2'}$$

   then

   $$\begin{aligned} \rho_1 \xi_1' &= \alpha \xi_1 + \beta \xi_2 & \rho_1 \xi_2' &= \gamma \xi_1 + \delta \xi_2 \\ \rho_2 \eta_1' &= \alpha \eta_1 + \beta \eta_2 & \rho_2 \eta_2' &= \gamma \eta_1 + \delta \eta_2 \\ \rho_3 \zeta_1' &= \alpha \zeta_1 + \beta \zeta_2 & \rho_3 \zeta_2' &= \gamma \zeta_1 + \delta \zeta_2 \end{aligned}$$

   or

   $$x_i' = (\gamma x_i + \delta)^{-1}(\alpha x_i + \beta).$$

   Of these four coefficients, only three are essential, but without the help of the commutative law none of them can be removed. With the commutative law, one coefficient can be normalized to 1 and there remain three linear equations for the three essential coefficients:

$$\gamma x_i x'_i + \delta x'_i - \alpha x_i - \beta = 0.$$

These equations can be solved for $\alpha, \beta, \gamma, \delta$ in the usual way. Thus the projective map is uniquely determined by three corresponding pairs. It follows from well known facts of analytic geometry that it may be realized by joining points and intersecting lines.

One sees, further, that Pappus' theorem implies not only the fundamental theorem of projective geometry, but in fact *all* plane configuration theorems are consequences of Pappus' theorem. A configuration theorem is in general a statement of the following kind. Given finitely many points and lines with prescribed incidences, one constructs in a prescribed way finitely many new points. Then the configuration theorem says that certain incidences hold in the extended figure, independently of the choice of initial elements.

2. **Basing plane projective geometry on weaker hypotheses.**

   (a) Projective incidence+Pappus $\longrightarrow$ Field
   (b) Projective incidence+$D_0$ $\longrightarrow$ Skew field
       +Archimedean postulate $\longrightarrow$ Field
   (c) Projective incidence+$D_2$ $\longrightarrow$ Alternative field
       +Archimedean postulate $\longrightarrow$ Field

   Thus addition of the Archimedean postulate—as was shown earlier—makes each skew or alternative field into a field, which is necessary for the fundamental theorem of projective geometry (by the first item in this section).

3. **Meaning of the general Desargues theorem $D_0$.**

   It is well known that $D_0$ can be proved from the spatial incidence axioms. In other words, embeddability of the plane in space is sufficient for the validity of $D_0$. It is also necessary.

   Indeed, suppose we have a plane geometry with incidence and $D_0$. Then we may introduce a coordinate geometry over a skew field in which the usual definitions of linear figures extend to three or more dimensions. The incidence theorems there follow by solving linear equations, which is always possible in a skew field.

4. **Foundation of projective geometry in space.**

   Spatial incidence and the Archimedean postulate suffice to derive the fundamental theorem. In fact, spatial incidence implies Desargues' theorem, and hence a geometry over a skew field, and addition of the Archimedean postulate makes the skew field a field. We know from the first item in this section that the field properties suffice to prove the fundamental theorem.

5. **Order axioms are inessential in the above investigations.**

6. **The parallel postulate may be eliminated.**

   It suffices to indicate how to derive the projective incidence properties of the plane.

   (a) The work of Schur or Pasch makes it possible to introduce ideal elements into the spatial incidence axioms (Axioms I) and the order axioms, without the parallel postulate.

   (b) Ordinary incidence and order in the plane do not suffice to introduce ideal elements. One needs a configuration theorem such as Desargues' theorem. It is essential to have at least the special Desargues theorem $D_2$.

## 3.5 Proof of Pappus' theorem from the congruence and parallel axioms

1. We mention in advance that the parallel postulate makes it possible to draw a line $b$, through a point $A$ outside line $b$, that meets $a$ at a prescribed angle $\angle \alpha$.

   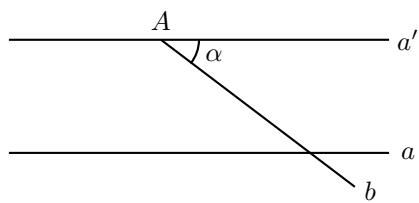

   To show this we draw through $A$ the unique parallel $a'$ to $a$, and transport the angle $\alpha$ to $a'$ at $A$. Its other leg is the desired line $b$, since alternate angles at parallels are equal.

2. Suppose that $A, B, C$ lie on $g$ and $A', B', C'$ lie on $g'$, and that $BC' \parallel B'C$ and $AC' \parallel A'C$. We have to prove that $AB' \parallel A'B$. The proof is based on the theorem about angles in a circle.

   $\alpha$) Draw $D'B$ through $B$ so that $\angle OCA' \cong \angle OD'B$. Since $CA' \parallel C'A$ we also have $\angle OCA' \cong \angle OAC'$, hence $\angle OD'B \cong \angle OAC'$. Then $A'D'BC$ is a cyclic quadrilateral, hence $\angle CBA' \cong \angle OD'C$ also. Both angles are subtended by the chord $A'C$ of the circle $\mathfrak{K}_1$

   $\beta$) Since $\angle OD'B \cong \angle OAC'$, and both are subtended by the chord $BC'$, $BD'C'A$ is a cyclic quadrilateral (in circle $\mathfrak{K}_2$). It follows that $\angle OAD' \cong \angle OC'B$ and, since $B'C \parallel C'B$, also $\angle OAD' \cong OB'C$.

   $\gamma$) Since $\angle OAD' \cong \angle OB'C$, $CD'B'A$ is a cyclic quadrilateral in circle $\mathfrak{K}_3$, where both angles are subtended by the chord $D'C$. Thus $\angle OAB' \cong \angle OD'C$ also.

   But $\alpha$) also implies $\angle OD'C \cong \angle OBA'$, from which it follows that $\angle OBA' \cong \angle OAB'$, that is, $B'A \parallel BA'$.                Q.E.D.

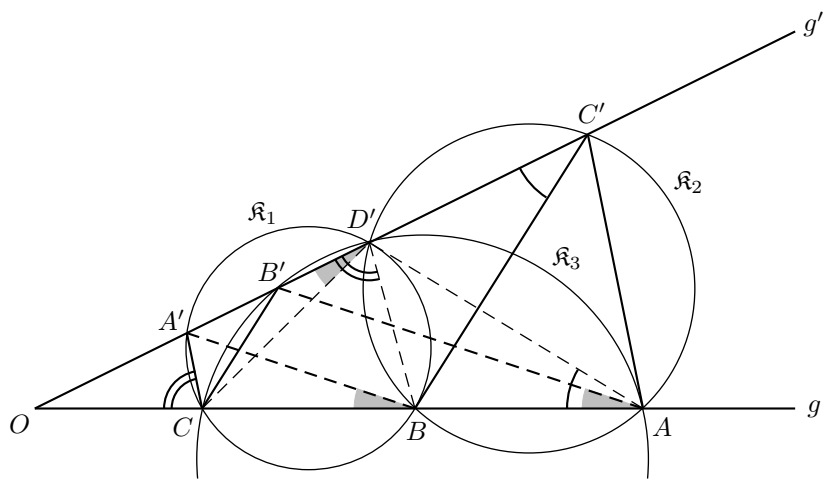

In *Math. Annalen*, vol. 64, Hjelmslev proved that projective geometry can be based on the plane axioms of incidence, order and congruence, *without* the parallel postulate and continuity. Thus there are three essentially different paths to projective geometry without the parallel postulate:

$$\begin{cases} \textbf{with continuity} \text{ and with } D_2 \\ \textbf{without continuity} \begin{cases} \text{with congruence} \\ \text{with Pappus' theorem} \end{cases} \end{cases}$$

We have yet to prove that the foundation can be built on continuity and the special Desargues theorem $D_2$, since our derivation used the parallel postulate. It remains to show that the ordinary plane axioms of incidence and order, together with $D_2$, suffice to introduce ideal elements in the plane (see Chapter 4).

## 3.6 Proof of Pappus' theorem in space

That is, from the axiom groups I, II and congruence axioms but without continuity.

1. Pappus' theorem may be proved with the help of spatial incidence and congruence. In fact, here we prove the general Pascal theorem on the hexagon in a curve of second order, whence the Pappus-Pascal theorem for a line pair follows.[10]

    As Dandelin 1825 showed, a one-sheeted hyperboloid contains two families of lines, $g_i$ and $h_i$, that are skew to the axis of the hyperboloid. All the lines in a family are skew to each other, but each line $g_i$ cuts each line $h_i$. Thus six lines that belong alternately to the two families form a spatial hexagon. For this "mystic hexagon" one can prove a theorem that gives Pascal's theorem immediately.

---

[10]Translator's note. Recall from the earlier note on page 52 that Moufang calls Pappus' theorem "Pascal's theorem" throughout. This is the first time she actually mentions Pappus, and the first time she takes up the general Pascal theorem.

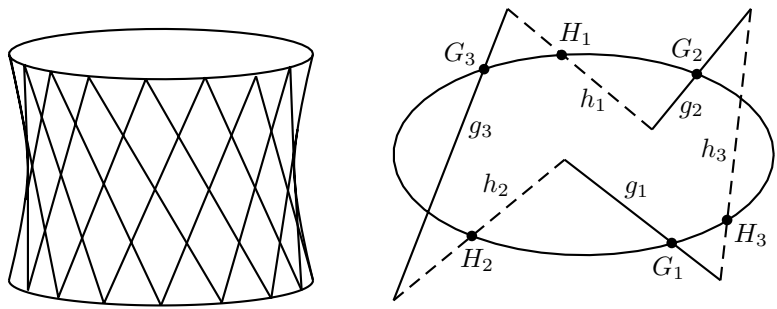

2. The intersection of a surface of degree 2 with a plane is a curve of degree 2, $\mathfrak{C}_2$. On $\mathfrak{C}_2$ we take six points $H_i$, $G_i$ ($i = 1, 2, 3$) as in the figure, for which we shall prove *Pascal's theorem*.

   It will be convenient to use the following abbreviations:

   $$\begin{aligned}[GH] &\quad \text{for the line connecting } G \text{ and } H \\ [gh] &\quad \text{for the plane connecting } g \text{ and } h \\ (gh) &\quad \text{for the intersection point of } g \text{ and } h\end{aligned}$$

   Claim:
   $$\left.\begin{aligned}([G_3H_1][G_1H_3]) &= Q \\ ([G_1H_2][G_2H_1]) &= R \\ ([G_2H_3][G_3H_2]) &= P\end{aligned}\right\} \quad \text{are collinear,}$$

   Or in words: *the intersections of opposite sides of hexagon on a $\mathfrak{C}_2$ lie on a line.*

   Proof. Let $g_i$, $h_i$ be the lines in the alternate families of lines on the hyperboloid that go through the points $G_i$, $H_i$ respectively.

   The pairs $g_1$, $h_3$ and $h_1$, $g_3$ span planes that meet in a line (similarly for the index pairs 2,1 and 3,2).

   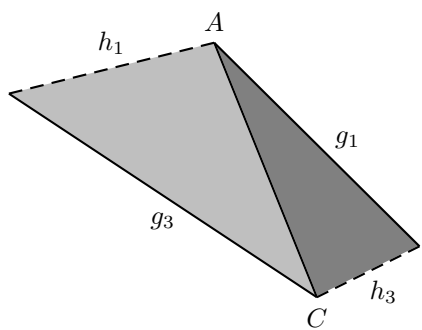

   Writing this down in detail we get

   $$\begin{aligned}([g_1h_3][g_3h_1]) &= [(g_1h_1)(g_3h_3)] = [AC] \\ ([g_2h_1][g_1h_2]) &= [(g_2h_2)(g_1h_1)] = [BA] \\ ([g_3h_2][g_2h_3]) &= [(g_3h_3)(g_2h_2)] = [CB].\end{aligned}$$

$\Delta ABC$ spans a plane $\mathfrak{E}$. Let the plane of the conic section be $\mathfrak{S}$. We have $\mathfrak{E} \neq \mathfrak{S}$, because $\mathfrak{E} = \mathfrak{S}$ implies, for example, that $(g_1 h_1)$ lies in $\mathfrak{S}$; but $g_1$ has only the point $G_1$ in common with $\mathfrak{S}$, so we should have $(g_1 h_1) = G_1$ and analogously $(h_1 g_1) = H_1$, which is a contradiction because $G_1 \neq H_1$.

Thus the intersection of $\mathfrak{E}$ and $\mathfrak{S}$ is a line $s$. But $P, Q, R$ lie on $s$, as we shall prove for $P$. (For $Q$ and $R$ the proof is analogous.)

It is clear that $P$ lies in $\mathfrak{S}$; we have to show that $P$ also lies in $\mathfrak{E}$. First, $P = ([H_2 G_3][H_3 G_2])$ lies on $[(g_2 h_2)(g_3 h_3)]$, because $[(g_2 h_2)(g_3 h_3)]$ is the intersection of the planes $[g_3 h_2]$ and $[g_2 h_3]$, and the latter planes are cut by the third plane $\mathfrak{S}$ in the lines $[G_3 H_2]$ and $[G_2 H_3]$, hence the latter lines meet on $[(g_2 h_2)(g_3 h_3)] = [BC]$.

But all points of $BC$ lie in $\mathfrak{E}$, by definition of $\mathfrak{E}$. Hence $P$ lies in $\mathfrak{E}$, and similarly so do $R$ and $Q$. That is, $P, Q, R$ lie simultaneously in $\mathfrak{E}$ and $\mathfrak{S}$, and hence in a line. $\hfill$ Q.E.D.

One therefore has the theorem: *the intersection lines of opposite planes of the mystic hexagram are the sides of a triangle whose vertices are the intersection points of opposite sides of the hexagram.*

This proves Pascal's theorem for the hexagon $G_1 H_1 G_2 H_2 G_3 H_3$.

To pass from a proper curve of degree 2 to a line pair it suffices to choose the cutting plane $\mathfrak{S}$ as a tangent plane of the hyperboloid.

3. **Existence of the one-sheeted hyperboloid on the basis of the congruence axioms** remains to be shown.

To do this we regard the one-sheeted hyperboloid as the collection of the two families of lines described at the beginning of this section, and look at the implications of the congruence axioms for the properties of motions and reflections. The most important properties are:

$\alpha$) Two intersecting lines are mapped into each other by reflection in a unique plane, the so-called meridian plane or symmetry plane. Conversely, two lines mapped into each other by reflection intersect in a proper or improper point, that is, they lie in a plane.

$\beta$) If $b$ results from $a$ by rotation about $d$, and if $a$ and $b$ do not lie in a plane perpendicular to $d$, then $a$ and $b$ are skew.

Proof. Rotation about $d$ maps all, and only, the points of $d$ into themselves. Since neither $a$ nor $b$ meets the axis $d$, no point of $a$ or $b$ is fixed. Moreover, no point $P$ of $a$ goes to a point $Q \neq P$ of $a$, otherwise $PQ = a$ would lie in a plane perpendicular to $d$, contrary to hypothesis. Thus $a$ and $b$ can have no point in common: as a point of $b$, such a point must come from a point of $a$, but as a point of $a$ it comes from another point of $a$, contrary to what we have just proved. Thus $a$ and $b$ are skew. $\hfill$ Q.E.D.

We also mention without proof:

$\gamma$) An odd number of reflections in planes with a common axis is replacable by reflection in one plane through this axis.

δ) An even number of reflections in planes through an axis may be replaced by a rotation about $d$.

The existence of a mystic hexagram may now be proved. Take the cutting plane $\mathfrak{S}$ tangential to the hyperboloid, so the intersection curve becomes a line pair $g$, $h$. We therefore begin with two intersecting lines $g$, $h$ with symmetry plane $\mathfrak{E}$. In $\mathfrak{E}$ one chooses a line $d$ not passing through $(gh)$ and not perpendicular to $[gh]$. Three points $G_i$ are given on $g$, and three points $H_i$ on $h$. Let $S$ be the reflection in $\mathfrak{E}$ that exchanges $g$ and $h$.

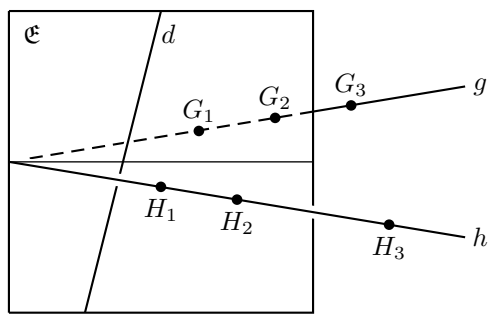

One constructs the planes

$$[dG_i] = \mathfrak{E}_{i1} \quad \text{and} \quad [dH_i] = \mathfrak{E}_{i2}.$$

They constitute a pencil of planes with axis $d$. I reflect the line $h$ to the line $h_i$ in the plane $\mathfrak{E}_{i2}$, abbreviating this by

$$h \xrightarrow{T_i} h_i$$

Similarly, let $g \xrightarrow{S_i} g_i$ denote reflection of $g$ in $\mathfrak{E}_{i1}$.

I claim: $g_i$ and $h_i$ are the sides of a mystic hexagram.

To begin with, $g_i$ goes through $G_i$, $h_i$ goes through $H_i$.

We must therefore show

1) Each $g_i$ meets each $h_k$.
2) Each $g_i$ is skew to each $g_k$ (and similarly for the $h_i$).

Proof of 1): $h_i \xrightarrow{T_i} h \xrightarrow{S} g \xrightarrow{S_k} g_k$.

Hence by γ) and α) there is an intersection point $(h_i g_k)$.

Proof of 2): Since the square of each reflection is the identity, they are involutions. That is, the same reflection that brings $h$ to $h_i$ also brings $h_i$ to $h$. Thus

$$h_i \xrightarrow{T_i} h \xrightarrow{T_k} h_k$$

and, by δ, the map from $h_i$ to $h_k$ is a rotation about $d$. Hence $h_i$ and $h_k$ are skew by β).

The existence of the mystic hexagram on the basis of congruence axioms (or more precisely, from theorems on reflections and rotations they imply) is thereby proved.                                                          Q.E.D.

**Summary. Projective geometry in space can be founded in two ways:**

| Incidence | Incidence |
|---|---|
| Order | Order |
| Continuity (Archimedean postulate) | Congruence (mystic hexagram) |
| Without parallel postulate | Without parallel postulate |

The introduction of ideal elements is achieved with the help of spatial incidence and order. The approach on the left is taken more often than necessary, because the usual coordinate geometry is Archimedean, whereas projective geometry need not be. For example, the coordinate geometry over the field of rational functions in $t$ or its extension by square roots to a field of algebraic functions of $t$, which is ordered but not Archimedeanly ordered, is a projective geometry, since the incidence axioms and Pappus' theorem hold.

It is still an open problem to find how much must be assumed about space—including incidence and order but omitting congruence—in order to prove Pappus' theorem.

# Chapter 4

# Introducing ideal elements in the plane using $D_2$

The introduction of ideal elements using only plane incidence and order has not yet been achieved. In space, however, the introduction of ideal elements on the basis of Hilbert's incidence and order axioms offers no difficulty, as Schur and Pasch have shown. The fundamental difference between the plane and space appears here: in the plane it is necessary to add configuration theorems in order to prove theorems that in space can be proved from the incidence axioms alone.

In what follows we use an idea of Dehn (carried out in R. Moufang, *Math. Ann.*, vol. 105 (1931), p. 759 ff.) for introducing ideal elements in the plane on the basis of the theorem $D_2$. It is well known that the general Desargues theorem is sufficient. The corresponding constructions with inaccessible points in the plane are given, for example, by Enriques (*Questions in Elementary Geometry*).

Here is how *Dehn's model geometry* is constructed.

1. We introduce *artificial points and artificial lines*, for which it can be shown that two points always have a common line and two lines always have exactly one common point. We divide the projective plane by three pairwise intersecting lines into four regions, numbered 1 to 4.

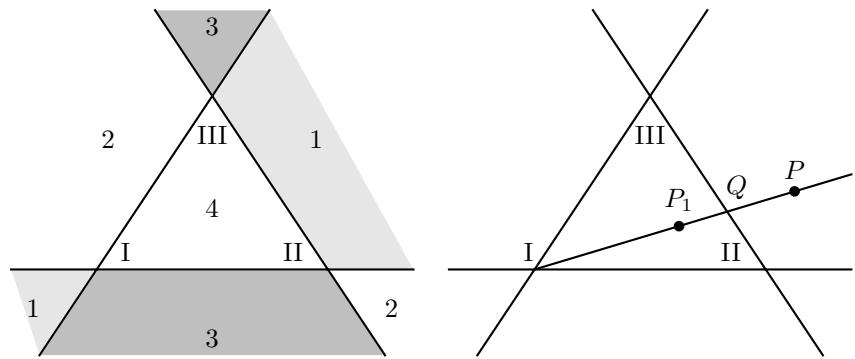

In the figure [on the left], the same number in two regions indicates that they correspond to one piece of the projective plane. Now we map the outer regions into the interior of the triangle I, II, III by collinear reflection so that I$QPP_1$ is a harmonic quadruple [in the figure on the right]. The triangle I, II, III is thereby covered four times, and each point of the triangle carries an index 1, 2, 3, or 4 according to which of the regions the corresponding point $P_i$ comes from. The latter points $P_i$ are the artificial points of this model geometry.

A line is here represented by a triangle $ABC$, whose construction follows from the definition of the points: the triangles I, II, III and $ABC$ are reciprocal, that is, I$A$, II$B$, III$C$ meet at a point.

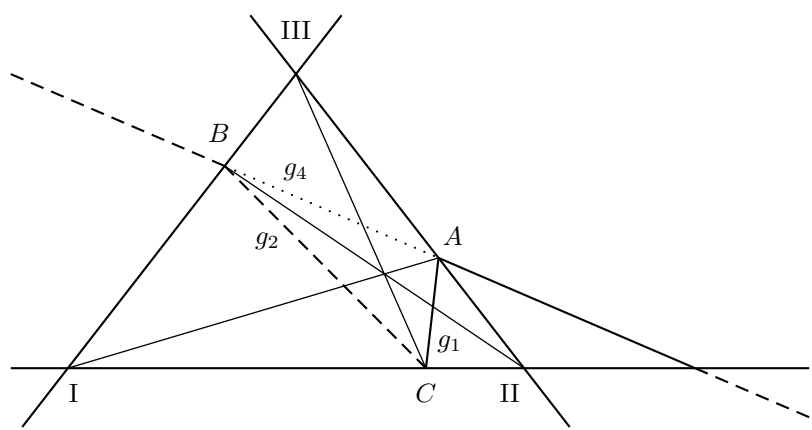

The sides of the triangle can be labelled in $\binom{4}{3} = 4$ ways, as in the following table:

| $AB$ | $BC$ | $CA$ |
|------|------|------|
| $g_4$ | $g_1$ | $g_2$ |
| $g_3$ | $g_1$ | $g_2$ |
| $g_2$ | $g_4$ | $g_3$ |
| $g_1$ | $g_3$ | $g_4$ |

The vertices of the triangle simultaneously carry all four indices, since they lie in the four regions simultaneously. Thus all four sheets come together at the vertices. The sides of the triangle are doubly covered: for example, the points I, II carry either the indices 3, 4 or the indices 1, 2, so sheets 3 and 4 meet along side I II and so too do sheets 1 and 2.

2. **Proof of the projective incidence axioms**

For this we use ordinary plane incidence and order, and theorem $D_2$ for accessible points only.

In Section 3.1 we used $D_2$ to derive the properties of harmonic quadruples, briefly: exchangeability of pairs, invariance under projection, harmonic quadruples in the complete quadrilateral. All three theorems use only accessible points.

1) **Connection of two points**

   $\alpha$) Points with the same index: the connecting line is shown in the figure as the inner triangle.
   Essentially, one has to construct the unique triangle $ABC$ through $P_iQ_i$ and reciprocal to I, II, III, and then one has to label the sides according to the table above.

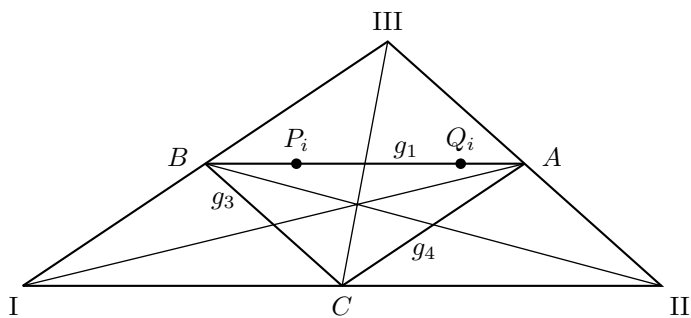

   $\beta$) Points with difference indices: for example, $P_4$, $P_1$.
   The connecting line is constructed as follows.

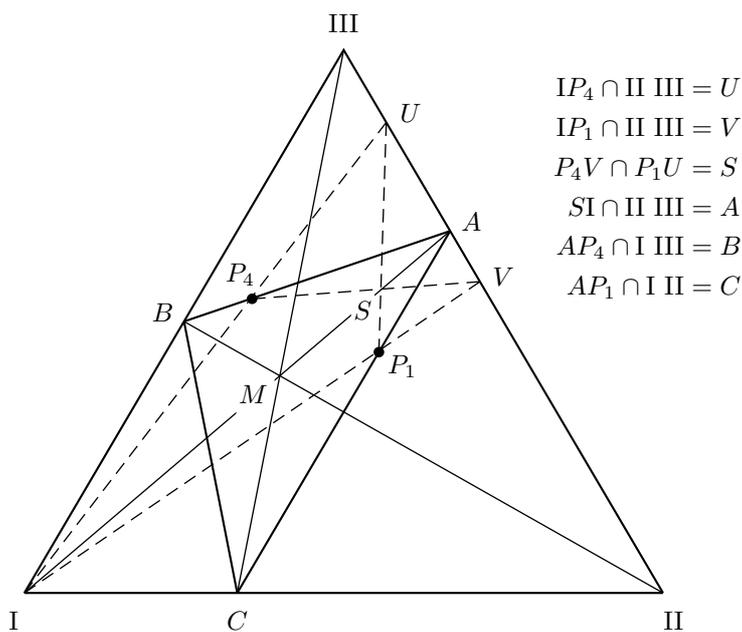

$$\text{I}P_4 \cap \text{II III} = U$$
$$\text{I}P_1 \cap \text{II III} = V$$
$$P_4V \cap P_1U = S$$
$$S\text{I} \cap \text{II III} = A$$
$$AP_4 \cap \text{I III} = B$$
$$AP_1 \cap \text{I II} = C$$

$\Delta ABC$ is then the artificial triangle, without labelling as yet. If $C\text{III} \cap \text{I}A = M$ we must show that $\Delta ABC$ is reciprocal to $\Delta \text{I II III}$, that is, that $B, M, \text{II}$ are collinear.

To prove this we consider the quadrilateral $\text{I}P_4AV$ with adjacent vertex $U$. Let
$$P_4A \cap US = G.$$

Then $\{GP_1US\}$ by definition of harmonic quadruples. Let

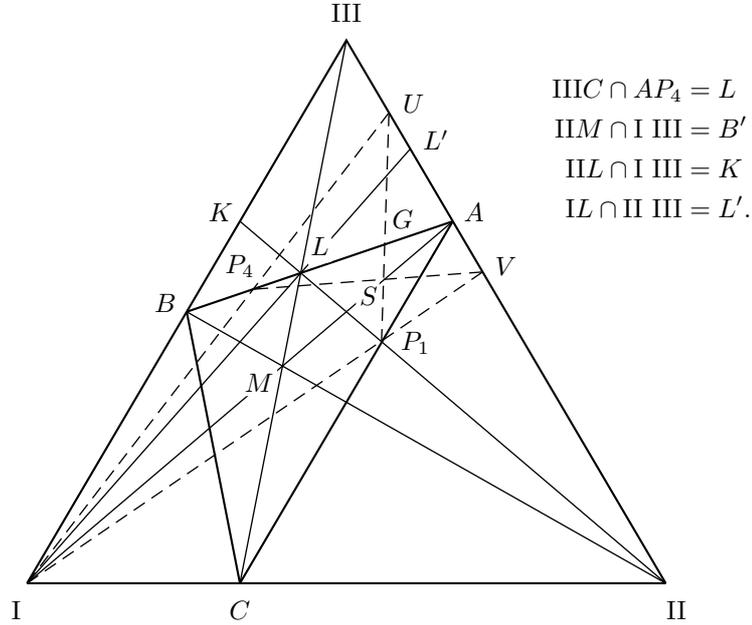

$$\text{III}C \cap AP_4 = L$$
$$\text{II}M \cap \text{I III} = B'$$
$$\text{II}L \cap \text{I III} = K$$
$$\text{I}L \cap \text{II III} = L'.$$

We must therefore show that $B = B'$.
We have
$$\{GP_1US\} \stackrel{A}{\overline{\wedge}} \{LC\text{III}M\} \stackrel{\text{II}}{\overline{\wedge}} \{K\text{I II}B'\}$$

This notation means, for example, that the harmonic quadruple $\{GP_1US\}$ goes to the harmonic quadruple $\{LC\text{III}M\}$ by projection from $A$. The sign $\overline{\wedge}$ is usual in projective geometry for a perspective or projective relationship.
We also have
$$\{LC\text{III}M) \stackrel{\text{I}}{\overline{\wedge}} \{L'\text{II III}A\} \stackrel{L}{\overline{\wedge}} \{IK\text{III}B\}.$$

Therefore, by the uniqueness of the fourth harmonic point to three given points, we have $B = B'$, as required.

This secures the existence of the reciprocal triangle. Now we have to prove its uniqueness. For this we must show that if $ABC$ is a reciprocal triangle to $\Delta$ I II III, and if $P_4$ is an arbitrary point on $BA$ and $P_1$ is an arbitrary point on $AC$, and if $U, V, S$ have the old meaning, then $S, A, I$ are always collinear.

Consider the quadrilateral $BA\text{III}$ with the harmonic quadruple $\{LC\text{III}M\}$.

$$\{LC\text{III}M\} \stackrel{A}{\overline{\wedge}} \{GP_1US'\}, \quad \text{whence} \quad S' = MA \cap UP_1.$$

We therefore have to show that $S = S'$. We consider the quadrilateral $IP_4AV$. Let

$$AI \cap P_4V = T, \quad UT \cap AP_4 = G', \quad UT \cap IV = P_1'.$$

Thus $\{G'P_1'UT\}$.

Now the points $G', P_1', U, T$ are harmonic in the quadrilateral $IP_4AV$, independently of whether $B_1A$ and $IV$ meet at an inaccessible point or not. Indeed, if one draws the auxiliary line $VG'$, the quadrilateral $P_4TAU$ has purely accessible vertices and hence a harmonic quadruple, which goes to $G', P_1', U, T$ by projection from I.

$$\{G'P_1'UT\} \stackrel{A}{\barwedge} \{GZUS'\} \quad \text{with} \quad Z = AP_1' \cap UP_1,$$

therefore $Z = P_1$. But then $P_1 = P_1'$ also, hence $T = S = S'$, as required.

Thus through two given points with different indices there is exactly one triangle $ABC$ reciprocal to I II III, obtainable by the construction above. Moreover, its labelling is uniquely determined by the table.                                                Q.E.D.

Thus it is shown that two distinct points have exactly one connecting line in our artificial geometry.

2) **Intersection of two lines**

In the accompanying figure,

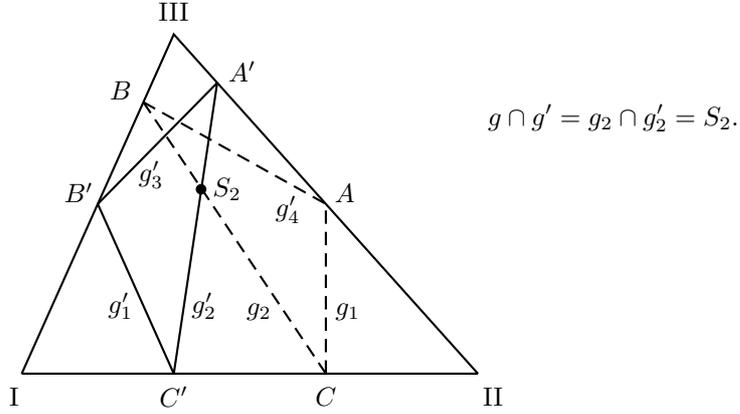

$g \cap g' = g_2 \cap g_2' = S_2.$

In finding the intersection of the lines $g$ and $g'$, the only points that come into consideration are the crossing points of like-labelled sides of the triangles $ABC$ and $A'B'C'$.

Two reciprocal triangles never cross at more than four places. Of these, only one is a crossing in the sense of our geometry, namely the intersection of two like-labelled triangle sides.

To prove this we first develop a few auxiliary results.

By Pasch's theorem, an ordering is preserved by projection:

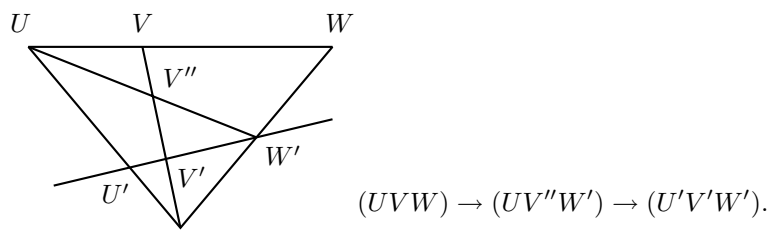

$$(UVW) \to (UV''W') \to (U'V'W').$$

We use this to show that $(BB'\text{III})$ and $(AA'\text{III})$ imply $(CC'\text{II})$. Let

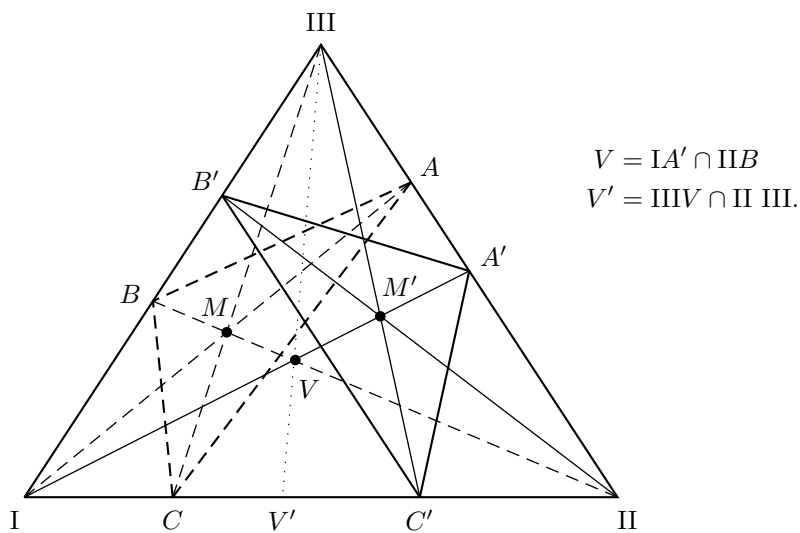

$$V = \text{I}A' \cap \text{II}B$$
$$V' = \text{III}V \cap \text{II III}.$$

Now
$$(BB'\text{III}) \stackrel{\text{II}}{\barwedge} (VM'A') \stackrel{\text{III}}{\barwedge} (V'C'\text{II})$$
and
$$(A'A\text{III}) \stackrel{\text{I}}{\barwedge} (VMB) \stackrel{\text{III}}{\barwedge} (V'C\text{I})$$

Together, these give $(CC'\text{II})$, as required.

It follows that $BC$ and $B'C'$ do not meet, by Pasch's axiom applied to $\Delta \text{I}BC$ and $B'C'$. The same holds for $A'C'$. Thus the number of crossings is bounded by 4. (One easily sees that $(AA'\text{III})$ and $(BB'\text{III})$ give the same result. This exhausts all essentially different possibilities.)

One now sees easily that there is exactly one crossing point of like-labelled lines, by following up the $\binom{4}{2} = 6$ possibilities in the table on page 101.  Q.E.D.

3. **Topological connectivity** of the artificial geometry in the triangle.

   1) The surface is not orientable, because if one follows the given [dashed] closed path through $P_1$, which passes all the vertices of I II III, the orientation is reversed. The strip bounded by $g$ and $h$ therefore has the connectivity of the Möbius band.

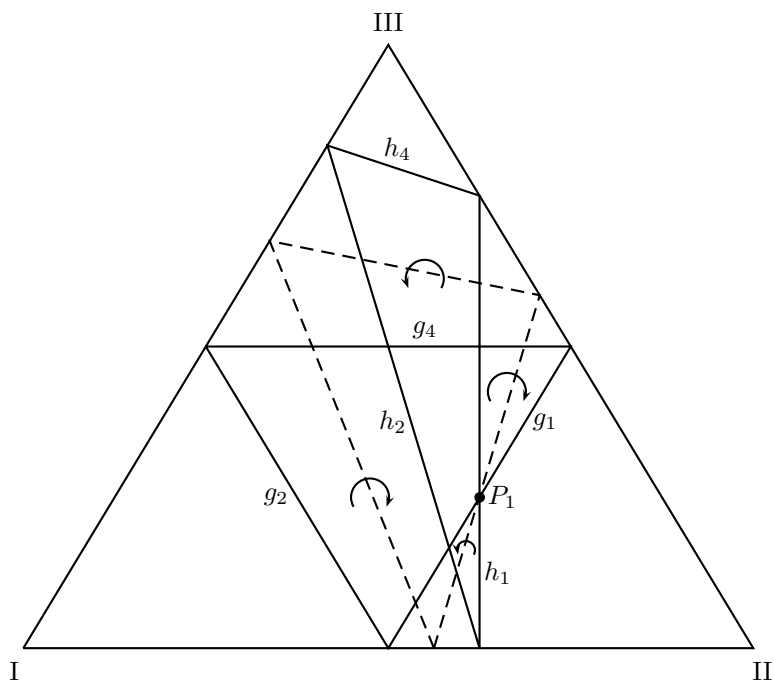

   The surface therefore contains an orientation-reversing path.

   2) The connectivity of the fourfold covering is represented in the accompanying sketch, where points to be identified are labelled the same.[1]

   The [Euler] characteristic is

   $$\chi = \text{faces} - \text{edges} + \text{vertices} = 4 - 6 + 3 = 1.$$

   Nonorientability and characteristic 1 characterize the projective plane. [In fact, nonorientability is redundant, since the projective plane is the only surface with Euler characteristic 1.]

---

[1]Translator's note. Moufang's sketch consists of four triangles, named 1, 2, 3, 4, with virtually illegible edge labelling in my copy of the manuscript. However, the edge labelling can be reconstructed from the diagram on page 100, and it is more instructive to do so. One sees from page 100 that the covering consists of four triangles 1, 2, 3, 4, each with the three vertices labelled I, II, III. The edges of the triangles are joined in pairs—for example, the edge between I and III on triangle 4 is joined to the edge between I and III on triangle 2, whereas the edge between I and III on triangle 1 is joined to the edge between I and III on triangle 3. This creates a surface with 4 triangular faces, 6 edges, and 3 vertices.

4. **Order**

    By cutting the covering of the triangle along $g$, each of the three sheets is cut once, decomposing them each into a quadrilateral and a triangle. Thus altogether there are three triangles, three quadrilaterals and an uncut triangle.

    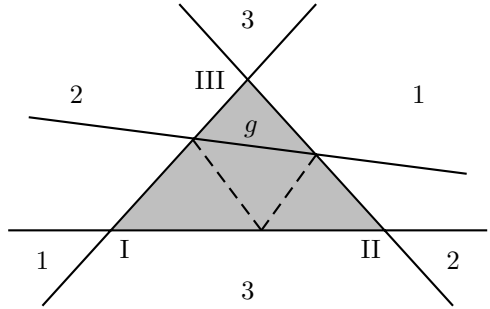

    The quadrilaterals are connected in such a way as to form a topological disk. Each line divides the disk into two parts so that two points on opposite sides of the cut cannot be connected without crossing the cut. But this is essentially the axiom of Pasch. Thus our model geometry satisfies not only projective incidence, but also projective order.

5. **The special Desargues theorem** $D_2$ holds in this artificial incidence geometry, as was proved by Smid (*Math. Ann.* 111 (1935)).

    Thus the significance of $D_2$ lies in the projective extendability of the plane and also in the algebraicization of synthetic plane geometry via a cartesian coordinate geometry over an alternative field.

    At the same time we have established that the fundamental theorem can be derived from $D_2$ and the Archimedean postulate, starting from incidence and order in the bounded plane.

# Appendix 1

## Methods for introducing ideal elements in the plane

1. **With the help of $D_2$**

   [As in the previous section.]

2. **With the help of $D_0$**

   This achieves the construction of the line connecting an accessible point $P$ with the inaccessible intersection of two lines $g$ and $h$ (see figure).

   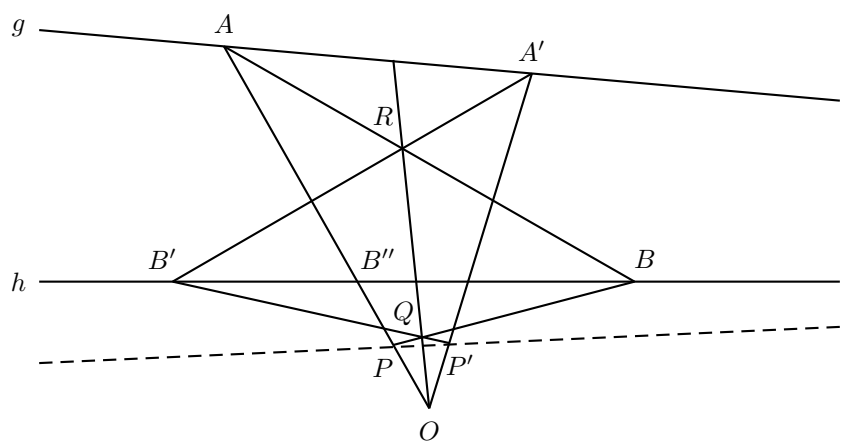

   Let $A$ and $A'$ be any points on $g$, $B$ and $B'$ any points on $h$, $AB \cap A'B' = R$, $O$ arbitrary on $AP$, $PB \cap RO = Q$, $B'Q \cap OA' = P'$.

   Then $\Delta A'B'P'$ is axial with $\Delta ABP$, with axis $RQ = RO$, and hence by $D_0$ they are in perspective, that is, $PP'$ goes through $g \cap h$.

   The construction is also correct when $B'$ becomes $B'' = h \cap AP$, in which case $D_2$ suffices for the proof.

3. **With the help of congruence axioms**, using Hjelmslev's method of reflection in lines.

   It is assumed that any point can be reflected in any line.

Then, for example, to connect an accessible point $P$ to the inaccessible intersection of lines $g$ and $h$, one reflects $P$ in $h$ to $P_h$, and $P_h$ in $g$ to $Q_1$. On the other hand, $P$ is reflected in $g$ to $P_g$, and $P_g$ in $h$ to $Q_2$. Let $s = Q_1Q_2$ and let $P'$ be the reflection of $P$ in $s$. Then $PP'$ is the desired connection. The proof is omitted here.

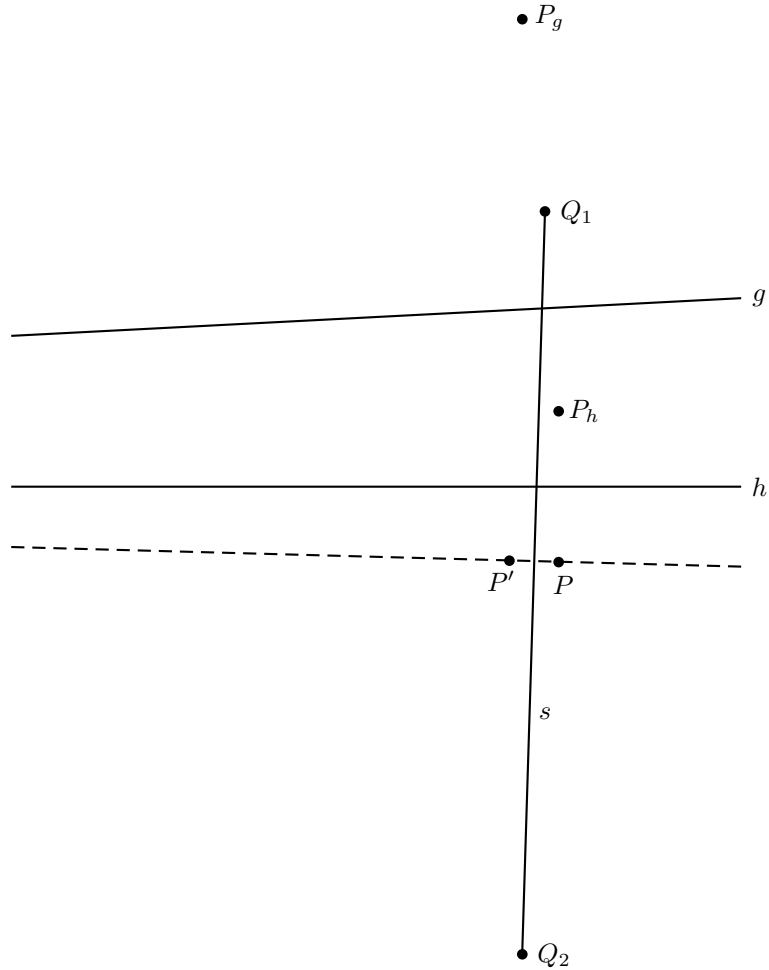

[Using $\xrightarrow{l}$ to denote reflection in the line $l$]

$$P \xrightarrow{h} P_h \xrightarrow{g} Q_1$$
$$P \xrightarrow{g} P_g \xrightarrow{h} Q_2$$
$$Q_1Q_2 = s$$
$$P \xrightarrow{s} P' \qquad PP' \text{ is the desired connection.}$$

# Appendix 2

### The role of the order axioms

In Hilbert's segment calculus, the vanishing of a polynomial indicates a configuration theorem (with the line of intersection points passing through $O$). Since Pappus' theorem determines an analytic geometry over a field, all plane configuration theorems must be analytically provable from it. That is, all plane configuration theorems follow from Pappus' theorem, as mentioned earlier.

The *problem of Dehn* is the (as yet incompletely answered) question whether there are incidence theorems $S$ "between" the theorems of Pappus and Desargues—which therefore do not imply Pappus and do not follow from Desargues. In the language of segment calculus, this asks whether there are computation rules which do not imply commutativity but do not hold in all skew fields.

Wagner (*Math. Annalen* 1937) has shown that the identical vanishing of a polynomial $P(a, b, \ldots, k)$ implies commutativity when the linear order axioms hold. Linear order is essential, as is shown by the example of $2 \times 2$ matrices over a field, for which we have the identity

$$\mathcal{L}(\mathcal{AB} - \mathcal{BA})^2 = (\mathcal{AB} - \mathcal{BA})^2 \mathcal{L},$$

which does *not* imply $\mathcal{AB} = \mathcal{BA}$.

The investigations of Wagner therefore answer the question in the negative for computation rules with *integer* coefficients.